\documentclass[12pt,twoside]{article}
\usepackage{graphicx, epsf, psfrag, epsfig}
\usepackage{amsmath, amsthm}
\RequirePackage{amsfonts,amssymb}

\def\const{\mbox{const}}

\def\mm#1{\begin{align}#1\end{align}}
\def\mmn#1{\begin{align*}#1\end{align*}}
\def\nn{\nonumber\\}
\def\Dt{\Delta t}

\def\dx{\partial_x}
\def\dy{\partial_y}
\def\vp{\varphi}
\def\htil{\tilde h}
\def\util{\tilde u}
\def\vtil{\tilde v}
\def\wtil{\tilde w}
\def\ctil{\tilde c}

\def\cn{\cos\theta}
\def\sn{\sin\theta}
\def\AvC#1{\,\frac1{2\pi}\int_0^{2\pi}#1d\theta} 
\def\AvS0#1{\,\frac1{2\pi}\int_0^{2\pi}#1d\theta} 
\def\AvD0#1{\,\big<#1\big>_{D_0}}

\def\vector#1{\left(\!\!\begin{array}{c}#1\end{array}\!\right)}
\def\vectorAvS0#1{\,\left<\!\!\begin{array}{c}#1\end{array}\!\!\right>_{\!\!\!S_0}}
\def\flexibleAvS0#1{\,\left<\!\!#1\!\!\right>_{\!\!\!S_0}}

\def\tt{{\cal O}(\Dt^2)}
\def\ttt{{\cal O}(\Dt^3)}

\newtheorem{lemm}{Lemma}[section]

\def\Bbb#1{\mathbb{#1}}

\def\softd{{\leavevmode\setbox1=\hbox{d}%
          \hbox to 1.05\wd1{d\kern-0.4ex{\char039}\hss}}}
\newcommand{\trivek}[3]{\left (
    \begin{array}{c}
    #1 \\ #2 \\ #3 \\
    \end{array} \right )}
\def\bold#1{\mbox{\boldmath $#1$}}

\newcommand{\umat}[1]{\bold{#1}}
\newcommand{\uu}[1]{\bold{#1}}
\newcommand{\D}{\Delta}

\newtheorem{theorem}{Theorem}[section]

\newcommand{\R}{\mathbb{R}}

\usepackage{fancyhdr}
\begin{document}

%
\def\leftrightscript#1#2#3{%
   \bgroup             
   \kern\scriptspace   
   {\vphantom{#2}}^{#1}
   \kern -\scriptspace 
   {#2\kern0pt}^{#3}
   \egroup             
}
%
%
\def\leftscript#1#2{%
   \kern\scriptspace
   {\vphantom{#2}}^{#1}%
   \kern -\scriptspace
   {#2\kern0pt}%
}
\def\bdot{\mbox{\tiny${^{\bullet}}$}}
\def\abs#1{\left|#1\right|} 
\def\norm#1{\left\|#1\right\|} 
\def\beq{\begin{equation}} \def\eeq{\end{equation}}
\def\bea{\bgroup\setlength\arraycolsep{2pt}\begin{eqnarray}} 
\def\eea{\end{eqnarray}\egroup}
\def\beas{\bgroup\setlength\arraycolsep{2pt}\begin{eqnarray*}} 
\def\eeas{\end{eqnarray*}\egroup}
\def\beg{\begin{eqngroup}} 
\def\eeg{\end{eqngroup}}
\def \sgn{\hbox{sgn}}
\def \rmand{\ {\rm and} \ }
\def \rmfor{\ {\rm for} \ }
\def \rmin{\ {\rm in} \ }
\def \rmif{\ {\rm if} \ }
\def \rmas{\ {\rm as} \ }
\def \rmon{\ {\rm on} \ }
\def \ii{\rm i}
\def \qed{\vrule height .5em width .5em depth 0em  \par \medskip}
\def \pd{\partial}
\def \<{\langle}
\def \>{\rangle}
\def \nab{\nabla}
\def \a{\alpha}
\def \b{\beta}
\def \ch{\chi}
\def \d{\delta}
\def \D{\Delta}
\def \L{\Lambda}
\def \g{\gamma}
\def \r{\rho}
\def \X{\Xi}
\def \G{\Gamma}
\def \o{\omega}
\def \O{\Omega}
\def \z{\zeta}
\def \k{\kappa}
\def \l{\lambda}
\def \m{\mu}
\def \n{\nu}
\def \s{\sigma}
\def \t{\tau}
\def \p{\pi}
\def \ph{\varphi}
\def \ps{\psi}
\def \qmin{q_{{\rm min}}}
\def \intw{\int_{-\infty}^{\infty}}
\def \inth{\int_0^{\infty}}
\def \us{u^{*}}
\def \ush{u^{\#}}
\def \fs{f^*}
\def \gs{g^*}
\def \hs{h^*}
\def \uh{\hat{u}}
\def \vh{\hat{v}}
\def \wh{\hat{w}}
\def \zh{\hat{z}}
\def \fh{\widehat{f}}
\def \gh{\hat{g}}
\def \hh{\hat{h}}
\def \Fh{\widehat{F}}
\def \Gh{\widehat{G}}
\def \Bh{\widehat{B}}
\def \Th{\widehat{T}}
\def \Wh{\widehat{W}}
\def \Uh{\widehat{U}}
\def \Nh{\widehat{N}}
\def \Rh{\widehat{R}}
\def \Sh{\widehat{S}}
\def \Eh{\widehat{E}}
\def \Fc{{\cal F}}
\def \tp{\tilde{p}}
\def \bt{\tilde{b}}
\def \ut{\tilde u}
\def \vt{\tilde v}
\def \Ut{\tilde U}
\def \Ht{\tilde H}
\def \Hc{{\cal H}}
\def \Mc{{\cal M}}
\def \Kc{{\cal K}}
\def \Wc{{\cal W}}
\def \Sb{\bar S}
\def \Ob{\bar{\O}}
\def \Lb{\bar L}
\def \Lt{\tilde L}
\def \Lp{\L~{\prime}}
\def \Lpp{\L~{\prime\prime}}
\def \bb{\bar \beta}
\def \cb{\bar c}
\def \xb{\bar x}
\def \yb{\bar y}
\def \ub{\bar u}
\def \vb{\bar v}
\def \wb{\bar w}
\def \zb{\bar z}
\def \fb{\bar f}
\def \gb{\bar g}
\def \Ub{\bar U}
\def \Nb{\bar N}
\def \Vb{\bar V}
\def \Wb{\bar W}
\def \Xb{\bar X}
\def \Yb{\bar Y}
\def \Zb{\bar Z}
\def \Lbar{\bar L}
\def \Bb{\bar B}
\def \BB{\bar\beta}
\def \xib{\bar \x}
\def \trip{\{f,g_0,g_1\}}
\def \triph{\{h,g_0,g_1\}}
\def \Trip{\{F,G_0,G_1\}}
\def \asy{{\rm asy}}
\def \pr{^{\prime}}
\def\half{\textstyle{\frac{1}{2}}\displaystyle}
\def\halfs{\scriptstyle{\frac{1}{2}}\displaystyle}
\newcommand{\jp}{_{j+1}}
\newcommand{\jm}{_{j-1}}
\newcommand{\mj}{^-_j}
\newcommand{\pj}{^+_j}
\newcommand{\ip}{_{i+1}}
\newcommand{\im}{_{i-1}}
\newcommand{\mi}{^-_i}
\newcommand{\e}{\mbox{e}}
\def\sx{\textstyle{\frac{1}{6}}\displaystyle}
\newcommand{\twoth}[1]{\textstyle{\frac{#1}{3}}\displaystyle}
\newcommand{\ep}{\epsilon}
\newcommand{\x}{(x)}
\newcommand{\y}{(y)}
\newcommand{\0}{(0)}
\def\ghj{G^h_j}
\newcommand{\sfrac}[1]{\textstyle{\frac{1}{#1}}\displaystyle}
\newcommand{\BM}[1]{\mbox{\boldmath{$#1$}}}
\newcommand{\DD}{\mbox{D}}
\newcommand{\dd}{\mbox{d}}

\bibliographystyle{plain}

\thispagestyle{empty}
\phantom{mm}
\vskip 1cm

\begin{center}
{\Large\bf Well-balanced finite volume evolution Galerkin methods for the shallow
water equations}
\\[0.3cm]
\vskip 1cm {\large M.~Luk\'a\v{c}ov\'a - Medvi\softd ov\'a\footnote[1]{Department of
Mathematics, University of Technology Hamburg, Schwarzenbergstra{\ss}e 95, 21 079
Hamburg, Germany, emails: lukacova@tu-harburg.de, kraft@tu-harburg.de},
S.~Noelle\footnote[2]{Division of Numerical Mathematics IGPM, RWTH Aachen,
Templergraben~55, 52062 Aachen, Germany, email: noelle@igpm.rwth-aachen.de} and
M.~Kraft$^1$}

\bigskip

\end{center}

\vskip 1cm
\begin{abstract}\noindent
We present a new well-balanced finite volume method within the framework of the
finite volume evolution Galerkin (FVEG) schemes. The methodology will be illustrated
for the shallow water equations with source terms modelling the bottom topography
and Coriolis forces. Results can be generalized to more complex systems of balance
laws. The FVEG methods couple a finite volume formulation with approximate evolution
operators. The latter are constructed using the bicharacteristics of
multidimensional hyperbolic systems, such that all of the infinitely many directions
of wave propagation are taken into account explicitly. We derive a well-balanced
approximation of the integral equations and prove that the FVEG scheme is
well-balanced for the stationary steady states as well as for the steady jets in the
rotational frame. Several numerical experiments for stationary and quasi-stationary
states as well as for steady jets confirm the reliability of the well-balanced FVEG
scheme. \end{abstract}



\section{Introduction}
\label{s1} \setcounter{equation}{0}
%
%

Consider the balance law in two space dimensions
\begin{equation}
\label{eq1} \uu{u}_t + \uu{f}_1(\uu{u})_x + \uu{f}_2(\uu{u})_y = \uu{b}(\uu{u},x,y),
\end{equation}
where $\uu{u}$ is the vector of  conservative variables, $\uu{f}_1,\, \uu{f}_2$ are
flux functions and $\uu{b}(\uu{u},x,y)$ is a source term.  In this paper we are
concerned with the finite volume evolution Galerkin (FVEG) method of
Luk\'a\v{c}ov\'a, Morton and Warnecke, cf.\ \cite{fvca3}-\cite{3d} and \cite{smhd}.
The FVEG methods couple a finite volume formulation with approximate evolution
operators which are based on the theory of bicharacteristics for the first order
systems \cite{mathcom}. As a result exact integral representations for solutions of
linear or linearized hyperbolic conservation laws can be derived, which take into
account all of the infinitely many directions of wave propagation.

In the finite volume framework the approximate evolution operators are used to
evolve the solution along the cell interfaces up to an intermediate time level
$t_{n+1/2}$ in order to compute fluxes. This step can be considered as a predictor
step. In the corrector step the finite volume update is done. The FVEG schemes have
been studied theoretically as well as experimentally with respect to their stability
and accuracy. Extensive numerical experiments confirm  robustness, good
multidimensional behaviour, high accuracy, stability, and efficiency of the FVEG
schemes, see Section \ref{section:num_exp} and also references \cite{jcp,sisc}. We
refer the reader to \cite{abgrall,bito,deconinck,roe,fey2,KroegerNoelle,lv,noelle}  for
other recent multidimensional schemes.

\medskip \noindent
For balance laws with source terms, the simplest approach is to use the operator
splitting method which alternates between the homogeneous conservation laws
$$
\uu{u}_t + \uu{f}_1(\uu{u})_x + \uu{f}_2(\uu{u})_y = 0
$$
and the ordinary differential equation
$$
\uu{u}_t  = \uu{b}(\uu{u},x,y)
$$
at each time step. For many situations this would be  effective and successful.
However, the original problem (\ref{eq1}) has an interesting structure, which is due
to the interplay between the differential terms and the right-hand-side source term
during the time evolution. For many flows which are of interest in geophysics, the
terms are nearly perfect balanced. If these terms are treated separately in a
numerical algorithm, the fundamental balance may be destroyed, resulting in spurious
oscillations. In particular, we will be interested in approximating correctly
equilibrium states or steady states, for which
$$
\uu{f}_1(\uu{u})_x + \uu{f}_2(\uu{u})_y = \uu{b}(\uu{u},x,y),
$$
and we want to approximate perturbations of such equilibrium states. Equilibrium
solutions play an important role because they are obtained usually as a limit when
time tends to infinity.

\medskip \noindent In this paper we present an approach which allows to incorporate
treatment of the source in the framework of the FVEG schemes without using the
operator splitting approach. The key ingredient is a new approximate representation
of the multi-dimensional solution which contains the full balance of the hydrostatic
pressure and the source terms, see Lemma \ref{lema32}. This new representation
allows to apply a recent non-standard quadrature rule (see \eqref{eq:sgn} and also
\cite{sisc}) to the mantle integrals of the bicharacteristic cone. This leads to a
surprisingly simple, accurate and efficient approximate evolution operator and hence
to an efficient, accurate and well-balanced finite volume scheme. We call our scheme
the {\sl well-balanced finite volume evolution Galerkin scheme} (FVEG). We refer the
reader to
\cite{bouchut1,botta,bouchut2,gallouet,greenberg,kurganov,lv1,noelle1,shijin,XingShu2005}
and the references therein for other related approaches for well-balanced finite
volume and finite difference schemes.

\medskip \noindent Our paper is organized as follows. In Section~2 we introduce the
shallow water equations in a rotating frame, derive a class of equilibria which
contains the lake at rest, jets in the rotating frame and combinations of these two
solutions. Then we introduce the class of two-step finite volume schemes used
throughout the paper. In Theorem \ref{theorem:well-balanced-fv} we give sufficient
conditions which guarantee well-balancing. Section \ref{section:EGoperator} is
devoted to the EG (evolution Galerkin) predictor step. Applying the theory of
bicharacteristics for multidimensional first order systems of hyperbolic type the
exact evolution operator is derived. A well-balanced approximation of the exact
evolution operator, which preserves some interesting steady states exactly and also
works well for their perturbations, is given in Lemma \ref{lema32} and equations
\eqref{approx_pwc} -- \eqref{approx_pwb}. In Theorem \ref{theorem:well-balanced-eg}
we prove that the FVEG scheme is well-balanced for the stationary steady states
(e.g.\ lake at rest), for the steady jets on the rotating plane as well as for
combinations of these two flows. In Section~4 we summarize the main steps of the
FVEG method and present its algorithm.   Numerical experiments for one and
two-dimensional stationary and quasi-stationary problems as well as for steady jets
presented in Section~5 confirm  reliability of the well-balanced FVEG scheme. We
also compare the accuracy and runtime of the new FVEG scheme with that of the recent
high order well-balanced WENO FV schemes of Audusse et a. \cite{bouchut1} and
Noelle et al. \cite{noelle1}. The
question of positivity preserving property of the scheme, i.e.~$ h > 0$, is not yet
considered here and will be addressed in our future paper.

\vskip 1cm

\section{Geophysical equilibria and well-balanced two-step finite volume schemes}
\label{s2} \setcounter{equation}{0}

\subsection{The shallow water equations}

There are many practical applications where the balance laws and the correct
approximation of their quasi-steady states are necessary. Some example include
shallow water equations with the source term modelling the bottom topography, which
arise in oceano\-graphy and atmospheric science, gas dynamic equations with
geometrical source terms, e.g. a duct with variable cross-section, or fluid dynamics
with gravitational terms. In what follows we illustrate the methodology on the
example of the shallow water equations with the source terms modelling the bottom
topography and the Coriolis forces. This system reads
\begin{equation}
\label{eq2} \uu{u}_t + \uu{f}_1(\uu{u})_x + \uu{f}_2(\uu{u})_y = \uu{b}(\uu{u}),
\end{equation}
where
\begin{eqnarray}
\nonumber
\uu{u}&=&\left(\begin{array}{c} h \\ hu \\
hv\end{array}\right),\
\uu{f}_1(\uu{u})=\left(\begin{array}{c} hu \\
hu^2+\frac12gh^2\\ huv\end{array}\right),
\\ \nonumber
\uu{f}_2(\uu{u})&=&\left(\begin{array}{c} hv \\ huv \\
hv^2+\frac12gh^2\end{array}\right),\
\uu{b}(\uu{u})=\left(\begin{array}{c} 0 \\
-ghb_x + f h v \\ -ghb_y - f h u\end{array}\right).
\end{eqnarray}

Here $h$ denotes the  water depth, $u,v$ are vertically averaged velocity components
in $x-$ and $y-$ direction, $g$ stands for the gravitational constant, $f$ is the
Coriolis parameter, and $b(x,y)$ denotes the bottom topography.

Note that these equations are also used in climate modelling and meteorology for
geostrophic flow, see, e.g., \cite{botta,klein}. For simulation of river or
oceanographic flows some additional terms modelling the bottom friction need to be
considered as well.

\subsection{Equilibria}

Many geophysical flows are close to equilibrium, or stationary state. It is easiest
to identify these states when the system is written in primitive variables,

\begin{equation}
\label{sw_primitive} \uu{w}_t +   \umat{A}_1(\uu{w}) \uu{w}_{x} + \umat{A}_2(\uu{w})
\uu{w}_{y} = \uu{t}(\uu{w}),
\end{equation}
\mm{ \label{sw_jacobians}
 \uu{w} =  \left (\begin{array}{c}
                      h \\
                      u    \\
                      v    \\
                      \end{array} \right ),
  \umat{A}_1 = \left ( \begin{array}{ccc}
                      u & h & 0 \\
                      g & u & 0 \\
                      0 & 0 & u \\
                      \end{array} \right ),
  \umat{A}_2 = \left ( \begin{array}{cccc}
                      v & 0 & h  \\
                      0 & v & 0     \\
                      g & 0 & v  \\
                      \end{array} \right ),
 \uu{t} =  \left ( \begin{array}{c}
                        0 \\
                       -g b_x + f v \\
                       -g b_y - f u
               \end{array} \right ).
}
Here we consider states which are both stationary,
\mm{(h,u,v)_t=0,}
and constant along streamlines,
\mm{(\dot h,\dot u,\dot v)=0,}
where $\;\dot{} =\partial_t+u\partial_x+v\partial_y$ is the material derivative. For
such states, we obtain
\mm{uh_x+vh_y = uu_x+vu_y = uv_x+vv_y = 0.}
Since the velocity vector is now constant along the streamlines, these become
straight lines. It is then natural to align the coordinates with the streamlines.
The desired solution has to satisfy the conditions
\mm{u &= 0 \\ v_y &= 0 \\ vh_y &= 0 \\ g(h+b)_x &= fv \\ g(h+b)_y &= 0.}
In the region $\{(x,y)|\, v(x)=0\}$ we obtain the lake at rest solution, where the
water level  $h+b$ is flat. When $v(x)\neq0$, we must have $h_y$ = 0 and hence
$b_y=0$. Hence the topography is locally one-dimensional, and along the rise of the
bottom we have a one-dimensional flow. This solution, which is well-known to
oceano\-graphers, is called the {\em jet in the rotational frame}. Due to the earths
rotation the jet exerts a sidewards pressure $fv$ onto the water, which is balanced
by a raise in the water level $g(h+b)_x$. In meteorological literature this state is
also called the geostrophic equilibrium.

For future reference we also define the primitive $(U,V)$ of the Coriolis force, as
introduced by Bouchut et al. in \cite{bouchut1}, via
\mm{\label{UV} V_x=\frac{f}{g}v \quad\mbox{and}\quad U_y=\frac{f}{g} u}
and the potential energies
\mm{\label{KL} K:=g(h+b-V) \quad\mbox{and}\quad L:=g(h+b+U).}

\subsection{Two-step finite volume schemes}

The FVEG (Finite Volume Evolution Galerkin) scheme which we propose for the balance
laws are time-explicit two-step schemes, similarly as Richtmyer's two-step version
of the Lax-Wendroff scheme \cite{RichtmyerMorton1967} and Colella's CTU (corner
transport upwind) scheme \cite{Colella1990}.

The first step, called predictor step, evolves the point value at a quadrature node
to the half-timestep. This can be done by a simple finite difference operator as in
\cite{RichtmyerMorton1967}, by one-dimensional characteristic theory as in
\cite{Colella1990} or by fully multidimensional, bicharacteristic theory as in
\cite{icfd} and related works of Luk\'a\v{c}ov\'a, Morton, Warnecke et al.. Our
predictor step is based on this bicharacteristic theory.

The second step is the standard finite volume update. It approximates the flux
integral across the interfaces by a quadrature of the fluxes evaluated at the
predicted states at the half-timestep.

We proceed as follows: in the present section we study the finite volume step (i.e.
the second of the two steps). This will give us sufficient conditions for
well-balancing which should be satisfied by the values computed in the predictor
step (the first step). In Section~\ref{s3} we will introduce the evolution Galerkin
predictor step and prove that it satisfies the sufficient conditions derived in the
present section.

In order to define the class of two step finite volume schemes, let us divide a
computational domain $\Omega$ into a finite number of regular finite volumes
$\Omega_{ij} = [x_{i-\frac12}, x_{i+\frac12}]\times[y_{j-\frac12}, y_{j+\frac12}] =
[x_i - \hbar/2, x_i + \hbar/2] \times [y_j -\hbar/2, y_j + \hbar/2]$, $i,j \in \Bbb
Z$, where $\hbar$ is the mesh size. Denote by $\uu{U}_{ij}^n$ the piecewise constant
approximate solution on a mesh cell $\Omega_{ij}$ at time $t_n$ and start with
initial approximations obtained by the integral averages $\uu{U}_{ij}^0 =
\int_{\Omega_{ij}} \uu{U}(\cdot, 0)$. Integrating the balance law (\ref{eq2}) and
applying the Gauss theorem on any mesh cell $\Omega_{ij}$ yields the following
finite volume update formula

\begin{equation}
\label{so_6} \uu{U}_{ij}^{n+1} = \uu{U}_{ij}^n - \lambda \sum_{k=1}^2
\delta^{ij}_{x_k} \bar{\uu{f}}_{k}^{n + 1/2} + \lambda \uu{B}_{ij}^{n+1/2},
\end{equation}

where $\lambda = \Delta t / \hbar$, $\Delta t$ is a time step, $\delta^{ij}_{x_k}$
stands for the central difference operator in the $x_k$-direction, $k=1,2$ and
$\bar{\uu{f}}_{k}^{n+1/2}$ represents an approximation to the edge flux at the
intermediate time level $t_n + \Delta t/2$. Further $\uu{B}_{ij}^{n+1/2}$ stands for
the approximation of the source term multiplied with the mesh size, $\hbar\uu{b}$.
The cell interface fluxes $\bar{\uu{f}}_{k}^{n+1/2}$ are evolved using an
approximate evolution operator denoted by $E_{\Delta t/2}$  to $t_n + \Delta t/2$
and averaged along the cell interface edge denoted by ${\cal E}$,

\begin{equation}
\label{so_7} \bar{\uu{f}}_k^{n + 1/2} := \sum_{j} \o_j \uu{f}_k(E_{\Dt/2}
\uu{U}^{n}(\uu{x}^j({\cal E}))).
\end{equation}
Here $\uu{x}^j({\cal E})$ are the nodes and $\o_j$ the weights of the quadrature for
the flux integration along the edges.

For simplicity, we introduce the following notation. Along the edges, we have
quadrature nodes $(x_{i\pm\frac12},y_{j+j'})$ resp.\ $(x_{i+i'},y_{j\pm\frac12})$,
where $i',j'\in\{0,\pm\frac12\}$. These nodes are already sufficient for the
midpoint, the trapezoidal and Simpson's rule. We denote the values given at the
predictor step by
\mm{ \vector{\hat h\\ \hat u\\ \hat v}_{i\pm\frac12,j+j'} \!\!\!\! &:=
\vector{h\\u\\v}_{i\pm\frac12,j+j'}^{n+\frac12} \quad \mathrm{and} \quad
\vector{\hat h\\ \hat u\\ \hat v}_{i+i',j\pm\frac12} \!\!\!\! :=
\vector{h\\u\\v}_{i+i',j\pm\frac12}^{n+\frac12} }
and the corresponding fluxes in $x$ resp. $y$-direction by
\mm{ \uu{\hat f}^1_{i\pm\frac12,j+j'} &:= \uu{f}_1((\hat h,\hat u,\hat v)_{i\pm\frac12,j+j'}) \\
\uu{\hat f}^2_{i+i',j\pm\frac12} &:= \uu{f}_2((\hat h,\hat u,\hat
v)_{i+i',j\pm\frac12}) .}
With this notation, we obtain
\mm{ \delta^{ij}_{x_1} \bar{\uu{f}}_{1}^{n + 1/2} &= \sum_{j'}\o_{j'}
\delta^{i,j+j'}_{x_1}\uu{\hat f}^1_{i,j+j'} \\
\delta^{ij}_{x_2} \bar{\uu{f}}_{2}^{n + 1/2} &= \sum_{i'}\o_{i'}
 \delta^{i+i',j}_{x_2}\uu{\hat f}^2_{i+i',j}.}
Finally we discretize the source term by
\begin{equation}
\label{source}
\uu{B}_{ij}^{n+\frac12} = -\;g\vector{0\\
\sum_{j'}\o_{j'}\;(\mu^{i,j+j'}_{x_1}\;\hat h_{i,j+j'})
\;(\delta^{i,j+j'}_{x_1}\;(\hat b-\hat V)_{i,j+j'}) \\[1ex]
\sum_{i'}\o_{i'}\;(\mu^{i+i',j}_{x_2}\;\hat h_{i+i',j})
\;(\delta^{i+i',j}_{x_2}\;(\hat b+\hat U)_{i+i',j})}.
\end{equation}
Here $U$ and $V$ are the discrete primitives of the Coriolis forces (see
\eqref{UV}), defined by
\mm{
\delta^{i,j+j'}_{x_1}\hat V_{i,j+j'}&=\hbar\frac{f}{g}\;\mu^{ij}_{x_1}\hat v_{i,j+j'}\\
\delta^{i+i',j}_{x_2}\hat U_{i+i',j}&=\hbar\frac{f}{g}\;\mu^{i+i',j}_{x_2}\hat
u_{i+i',j},}
and we have used the average operators
\mmn{
\mu^{ij}_{x_1}a &= (a_{i+1/2,j}+a_{i-1/2,j})/2\\
\mu^{ij}_{x_2}a &= (a_{i,j+1/2}+a_{i,j-1/2})/2.}
The following theorem states conditions which guarantee that the two-step finite
volume scheme \eqref{so_6} -- \eqref{so_7} is well-balanced for the lake at rest as
well as for the jet in the rotating frame.

\begin{theorem}\label{theorem:well-balanced-fv}
Suppose that the values $(\hat h,\hat u, \hat v)$ given by predictor step satisfy
for all $i,j,i',j'$
\begin{eqnarray}
\hat u_{i,j+j'} &=&0  \label{c1}\\
\delta_{y}^{i+i',j}\hat v_{i+i',j} &=&0 \label{c2}\\
\hat v_{i+i',j}\;\delta_{y}^{i+i',j}\hat h_{i+i',j} &=&0 \label{c3}\\
\delta_{x}^{i,j+j'}\hat K_{i,j+j'} &=&0 \label{c4}\\
\delta_{y}^{i+i',j}(\hat h_{i+i',j} + b_{i+i',j}) &=&0 \label{c5}
\end{eqnarray}
where $K$ is defined in \eqref{KL}. Then the finite volume scheme preserves the lake
at rest and the jet in the rotating frame.
\end{theorem}

\proof Since this argument is already standard for the lake at rest, we only sketch
it briefly for  the jet in the rotating frame. Let us study conservation of momentum
in the $y$-direction over cell $\O_{ij}$. Using (\ref{c1}) -- (\ref{c3}),
(\ref{c5}), and the discrete product rule
\mm{ \delta(\hat a\hat b) &= \delta(\hat a)\mu(\hat b) + \mu(\hat a)\delta(\hat b) }
it is straightforward to show that the sum of the flux differences and the source
term vanishes:
\mm{\nonumber &\; (hv)_{ij}^{n+1}-(hv)_{ij}^n\\
\label{eq:wb11} =&\; \sum_{j'}\o_{j'}\;\delta_{x}^{i,j+j'}(\hat h\hat u\hat v) +
\sum_{i'}\o_{i'}\left(\delta_{y}^{i+i',j}(\hat h\hat v^2+\frac{g}{2} \hat h^2) +
\;g\;(\mu_{y}^{i+i',j}\hat h)\;\delta_{y}^{i+i',j}(b+\hat U)\right).}
Now $\hat u_{i\pm\frac12,j+j'}=0$, and
\mmn{ \quad \delta_{y}^{i+i',j}(\hat h\hat v^2) &= \delta_{y}^{i+i',j}(\hat h\hat v)
\; \mu_{y}^{i+i',j}\hat v
+ \mu_{y}^{i+i',j}(\hat h\hat v) \; \delta_{y}^{i+i',j}\hat v \\
&= \delta_{y}^{i+i',j}\hat h \; \mu_{y}^{i+i',j}\hat v
+ \mu_{y}^{i+i',j}\hat \; \delta_{y}^{i+i',j}\hat v \\
&= 0.}
Dropping the corresponding terms in \eqref{eq:wb11}, we obtain
\mmn{(hv)_{ij}^{n+1}-(hv)_{ij}^n &= \; \sum_{i'}\o_{i'}\;(\mu_{y}^{i+i',j}\hat
h)\;\delta_{y}^{i+i',j}
g(\hat h+b+\hat U) \\
&=\; \sum_{i'}\o_{i'}\;(\mu_{y}^{i+i',j}\hat h)\;\delta_{y}^{i+i',j}\hat L \\
&=\;\; 0. }
This is the desired well-balanced property for the $y$-momentum. The $x$-momentum
and the balance of mass can be treated analogously. \hfill \qed

\section{The well-balanced approximate evolution operators}
\label{section:EGoperator} \setcounter{equation}{0}

The predictor step in the FVEG scheme will be based on exact and approximate
integral representations of solutions to the linearized shallow water equations. We
begin this section by formulating the exact integral representation. A few
clarifying remarks should help the reader to understand the structure of this
representation. Details of the derivation are given in Appendix
\ref{section:appendix_EG}. We proceed to derive two approximate integral
representations. For the first order scheme the approximate evolution operator
$E_{\Delta t/2}^{const}$ for the piecewise constant data is used. For the second
order method the continuous bilinear recovery $R$ is applied. In the case of
discontinuous solutions slopes in $R$ are limited yielding a discontinuous piecewise
bilinear recovery $R$, cf.\ Section~4 as well as \cite{jcp}. The predicted solution
at the quadrature nodes on the cell interfaces at the half timestep is obtained by a
suitable combination of $E_{\Delta}^{const}$ and $E_{\Delta}^{bilin}$,
\begin{equation}
\label{so_77}
E_{\Delta t/2}\uu{U}^{n} := E^{bilin}_{\Delta t/2} R \uu{U}^{n} +
E^{const}_{\Delta t/2} (1 - \mu_x^2 \mu_y^2) \uu{U}^{n},
\end{equation}
where $\mu_x^2 U_{ij} = 1/4(U_{i+1,j} + 2 U_{ij} + U_{i-1,j})$; an analogous
notation is used for the $y-$direction. It has been shown in \cite{sisc} that the
combination (\ref{so_77}) yields the best results with respect to accuracy as well
as stability among other possible second order FVEG schemes. It is particularly
important that the constant evolution term $E^{const}_{\Delta t/2} (1 - \mu_x^2
\mu_y^2) \uu{U}^{n}$  corrects the conservativity of the bilinear recovery and hence
of the intermediate solutions along cell-interfaces. If it is not used the scheme is
second order formally, but unconditionally unstable, cf.\ the FVEG-B scheme
\cite{sisc}.

Finally we show that approximate evolution operators lead to well-balanced two-step
finite volume schemes for the lake at rest and the jet in the rotational frame.

\subsection{Exact integral representation}
\label{s3}

We believe that the most satisfying methods for evolutionary problems are based on
the approximation of evolution operator or at least its dominant part. For the
two-step FVEG method, we use two fundamental evolution operators. One of the steps
is the classical finite volume update for the cell averages and uses the integral
form of the conservation law. Its well-balanced properties have been established in
Theorem \ref{theorem:well-balanced-fv}. The other step, which precedes the finite
volume update, is needed to predict point values in order to evaluate fluxes on cell
interfaces. It is here that the classical bicharacteristic theory comes into play.
It provides exact integral formulae for point values of solutions to
multidimensional hyperbolic systems.


Let $P:=( x, y,t_{n+1/2})$ be one of the quadrature points where the finite volume
fluxes will be evaluated, and let $\wtil=(\htil,\util,\vtil)$ be a suitable local
average of the solution around $P$. We will derive an exact integral representation
of the solution of the linearized shallow water equation at $P$. Similarly as in
\eqref{sw_primitive}, the linearized system in primitive variables reads

\begin{equation}\label{sw_linearized}
\uu{w}_t +   \umat{A}_1(\uu{\wtil}) \uu{w}_{x} + \umat{A}_2(\uu{\wtil})
\uu{w}_{y} = \uu{t}(\uu{w}),
\end{equation}

where the Jacobian matrices $A_1$ and $A_2$ are defined in \eqref{sw_jacobians}.


The homogeneous part of (\ref{sw_linearized}) yields a hyperbolic system. Fix a direction angle
$\theta$ with corresponding unit normal vector $(\cn,\sn)$. The matrix pencil
$\umat{A} \equiv  \umat{A}(\uu{\wtil})= \cn \umat{A}_1 + \sn \umat{A}_2$,
has three eigenvalues \begin{eqnarray} \nonumber \lambda_1 &=&
\cn\, \util + \sn\, \vtil - \ctil, \\ \nonumber \lambda_2 &=& \cn\, \util + \sn\, \vtil, \\
\nonumber \lambda_3 &=& \cn\, \util + \sn\, \vtil + \ctil, \end{eqnarray} and a full
set of right eigenvectors $$ \uu{r}_1 = \trivek{-1}{g\cn/\ctil}{g\sn/\ctil}, \: \:
\uu{r}_2 = \trivek{0}{\sn}{-\cn}, \:\:  \uu{r}_3 = \trivek{1}{g\cn/\ctil}{g\sn/\ctil
}, $$ where $c= \sqrt{gh}$ denotes the wave celerity. The eigenvalues
$\lambda_{1,3}$ correspond to fast waves, the so-called inertia-gravity waves,
whereas slow modes are related to $\lambda_2$. Analogously to the gas dynamics the
Froude number $Fr= | \uu{u} | / c$ plays an important role in the classification of
shallow flows. The shallow flow is called supercritical, critical or subcritical for
$Fr >1, Fr=1, $ and $Fr < 1$, respectively.

Applying the theory of bicharacteristics to the linearized system (\ref{sw_linearized})
yields an exact integral representation of the solution. Since the computations are
closely related to \cite{sisc}, we summarize the key steps only briefly and refer to
Appendix~A for further details.
\begin{itemize}
\item
Fix a point $P=(x,y,t_{n} + \tau)$, \ $\tau = \frac{\Delta t}{2}$.

For each spatial direction $(\cos(\theta),\sin(\theta)), \theta\in[0,2\pi)$ apply
the corresponding one-dimensional characteristic decomposition to two-dimensional
system \eqref{sw_linearized}.
\item
Integrate the resulting equations along each bicharacteristic curve from time $t_n$
to time $t_n + \tau$.
\item
Integrate the resulting equations over all direction angles $\theta$.
This gives a representation formula for the solution at the point $P$.
\end{itemize}

\begin{figure}[h]
\label{fig_11}
   \psfrag{p}{\small{$P=(x,y,t_n + \tau)$}}
   \psfrag{pp}{\small{$Q_0$}}
   \psfrag{qq}{\small{$Q(\theta)$}}
   \psfrag{x}{$x$}
   \psfrag{y}{$y$}
   \psfrag{z}{$t$}
   \begin{center}
      \epsfig{file=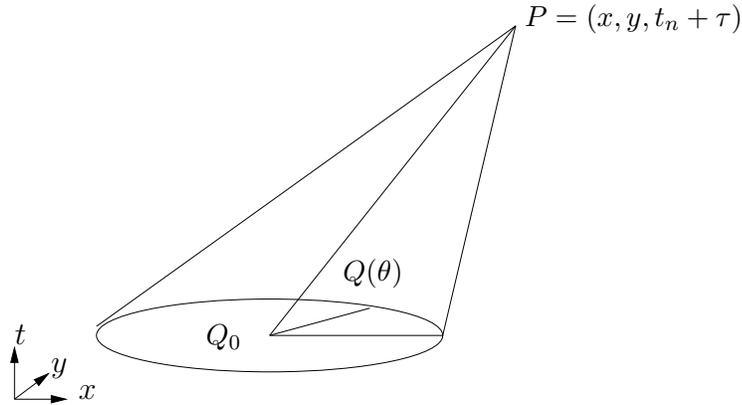,width=7cm}%
   \end{center}
\caption{Bicharacterestics  cone. }
  \hfill
\end{figure}

The EG integral representation derived in Appendix A then reads,
cf.~(\ref{eq3-ap})-(\ref{eq5-ap})
\begin{eqnarray}
\label{eq3} h\left( P\right) &=& \frac{1}{2\pi}\int_{0}^{2\pi}
   h\left(Q\right) - \frac{\tilde{c}}{g} \left( u\left(Q\right)\cos\theta+
    v\left(Q\right)\sin\theta \right) \dd\theta
    \phantom{mmmmmmmmmmmmmmmm}
\nonumber \\
&& - \frac{1}{2\pi} \int_{t_n}^{t_n + \tau} \frac{1}{t_n + \tau - \tilde{t}}
     \int_{0}^{2 \pi}\frac{\tilde c}{g}\left( u(\tilde Q) \cos \theta + v(\tilde Q) \sin \theta\right)
     \dd \theta \dd \tilde{t}
\\ \nonumber
&& + \frac{1}{2\pi} \tilde{c} \int_{t_n}^{t_n + \tau}  \int_{0}^{2 \pi} \left(
b_x(\tilde Q) \cos \theta + b_y (\tilde Q ) \sin \theta \right)     \dd \theta \dd
\tilde{t} \nonumber
\\ \nonumber
&& - \frac{1}{2\pi} \frac{\tilde{c} f}{g} \int_{t_n}^{t_n + \tau}  \int_{0}^{2 \pi}
\left( v(\tilde Q) \cos \theta - u (\tilde Q ) \sin \theta \right)     \dd \theta
\dd \tilde{t},
\end{eqnarray}
\begin{eqnarray}
\nonumber \label{eq4}
 u\left(P \right) &=&\frac{1}{2}u\left( Q_0\right) + \frac{1}{2 \pi}\int_{0}^{2\pi}
   - \frac{g}{\tilde{c}}  h\left(Q\right) \cos\theta
   + u \left(Q \right)\cos^2\theta  + v\left(Q \right)\sin\theta\cos\theta \, \dd\theta   \nonumber \\
&& - \frac{g}{2}  \int_{t_n}^{t_n + \tau} \left( h_x(\tilde Q_0) + b_x (\tilde
Q_0)\right)  \dd \tilde{t}
\\ \nonumber
&& - \frac{g}{2\pi}  \int_{t_n}^{t_n + \tau} \int_0^{2 \pi} \left(b_x(\tilde Q)
\cos^2 \theta + b_y(\tilde Q) \sin \theta \cos \theta \right) \dd \theta \dd
\tilde{t}
\\  \nonumber
&&+\frac{1}{2\pi} \int_{t_n}^{t_n + \tau} \frac{1}{t_n + \tau - \tilde{t}}
 \int_{0}^{2 \pi} \left( u(\tilde Q) \cos 2 \theta + v(\tilde Q) \sin  2\theta \right)
 \dd \theta \dd \tilde{t}
\\ \nonumber
&& + \frac{f}{2}  \int_{t_n}^{t_n + \tau}  v(\tilde Q_0) \dd \tilde{t} +
\frac{f}{2\pi} \int_{t_n}^{t_n + \tau} \int_0^{2 \pi} \left(v(\tilde Q) \cos^2
\theta - u(\tilde Q) \sin \theta \cos \theta \right) \dd \theta \dd \tilde{t},
\end{eqnarray}
\begin{eqnarray}
\label{eq5}
\nonumber  v\left(P \right) &=&\frac{1}{2}v\left( Q_0\right) + \frac{1}{2
\pi}\int_{0}^{2\pi}
   - \frac{g}{\tilde{c}} h\left(Q\right) \sin\theta
   + u\left(Q \right)\sin\theta\cos\theta + v \left(Q \right)\sin^2\theta
   \, \dd\theta   \nonumber \\
&& - \frac{g}{2}  \int_{t_n}^{t_n + \tau} \left( h_y(\tilde Q_0) + b_y (\tilde
Q_0)\right)  \dd \tilde{t}
\\ \nonumber
&& - \frac{g}{2\pi}  \int_{t_n}^{t_n + \tau} \int_0^{2 \pi} \left(b_x(\tilde Q) \sin
\theta \cos \theta + b_y(\tilde Q) \sin^2 \theta \right) \dd \theta \dd \tilde{t}
\\  \nonumber
&&+\frac{1}{2\pi} \int_{t_n}^{t_n + \tau} \frac{1}{t_n + \tau - \tilde{t}}
 \int_{0}^{2 \pi} \left( u(\tilde Q) \sin 2 \theta - v(\tilde Q) \cos  2\theta \right)
 \dd \theta \dd \tilde{t}
\\ \nonumber
&& - \frac{f}{2}  \int_{t_n}^{t_n + \tau}  u(\tilde Q_0) \dd \tilde{t} +
\frac{f}{2\pi} \int_{t_n}^{t_n + \tau} \int_0^{2 \pi} \left( v(\tilde Q) \sin \theta
\cos\theta - u(\tilde Q) \sin^2 \theta  \right) \dd \theta \dd \tilde{t}.
\end{eqnarray}

Evolution takes place along the bicharacteristic cone, see Fig.~1,  where
$P=(x,y,t_n + \tau)$ is the peak of the bicharacteristic cone, $Q_0=(x-\tilde{u}
\tau, y- \tilde{v} \tau, t_n)$ denotes the center of the sonic circle at time $t_n$,
$\tilde Q_0 = (x- \tilde{u} (t_n + \tau - \tilde t), y- \tilde{v} (t_n + \tau -
\tilde t), \tilde t)$, $\tilde Q  = (x - \tilde{u} (t_n + \tau- \tilde t) + c (t_n +
\tau - \tilde t) \cos \theta, y- \tilde{v}(t_n + \tau - \tilde t) + c (t_n + \tau -
\tilde t )\sin \theta, \tilde t)$ stays for arbitrary point on the mantle and $Q =
Q(\tilde t)\Big |_{\tilde t = t_n} $ denotes a point at the perimeter of the sonic
circle at time $t_n.$

At first view, it may seem difficult to interpret the terms in (\ref{eq3})-(\ref{eq5}). However, even
if one is not familiar with bicharacteristic theory, there is a simple explanation of all
terms, which we sketch in the following paragraph.

\subsection{A simple interpretation of the EG integral representation }

Our interpretation of the EG integral representation is based on a comparison with the approximate representation of the solution which results from a Taylor expansion along the streamlines.

The streamlines of the linearized shallow water equations are given by
$\{(x(t),y(t),t)|\dot x=\util,\dot y=\vtil\}$. As in Section~2.2, let $ \dot \vp :=
\vp_t+\util\vp_x+\vtil\vp_y $ be the material derivative of a function $\vp: \bar D
\to \Bbb R$.

Then the linearized shallow water equations (3.1) reduce to
\mm{
\dot h &= - \htil(u_x+v_y) \\
\dot u &= - K_x \\
\dot v &= - L_y.
}
This implies that
\mm{
\ddot h &= -\htil(\dot u_x+\dot v_y)
         = \htil (K_{xx} + L_{yy}) \\
\ddot u &= -\dot K_x
         = g \htil(u_x+v_y)_x- g(\util b_x+\vtil b_y)_x+fL_y \\
\ddot v &= -\dot L_y
         = g\htil(u_x+v_y)_y-g(\util b_x+\vtil b_y)_y-fK_x,
}
so the exact solution derived in (\ref{eq3})-(\ref{eq5}) can be approximated (up
to $\ttt$) by
\mm{ \vector{h_1\\u_1\\v_1} &:= \vector{h\\u\\v}_{\!\!\!\!\! Q_0} \!\!\! - \tau
\vector{\htil(u_x+v_y)\\K_x\\L_y}_{\!\!\!\!\! Q_0} \!\!\!  +\frac{\tau^2}2
\vector{\htil(K_{xx}+L_{yy})\\ g\htil(u_x+v_y)_x-g(\util b_x+\vtil b_y)_x+fL_y\\
g\htil(u_x+v_y)_y-g(\util b_x+\vtil b_y)_y-fK_x}_{\!\!\!\!\! Q_0} \!\!\!\! .
\label{eq:huv_1} }
Now we will indicate how to compare the RHS of the integral representation
(\ref{eq3})-(\ref{eq5}) with that of \eqref{eq:huv_1}. We will show that they agree
up to terms of $\ttt$. For this we need the following Taylor expansions on the sonic
circle and the bicharacteristic cone. For simplicity we introduce the notation
\[ \a := \cn, \quad \b := \sn.\]
\begin{lemm}\label{aux_lemma:taylor_expansions}
(Taylor expansions on the sonic circle and the bicharacteristic cone) Let
$\vp:C^{m+1}(\bar D;\R)$ and $\vp_0:=\vp(Q_0)$, $\tilde\vp_0:=\vp(\tilde Q_0)$.
Then
\mm{ \AvS0{\!\!\a^m\b^n\vp} &= \sum\limits_{k,l=0}^{k+l \leq m
}\frac{\AvC{\a^{k+m}\b^{l+n}}}{k!l!}(\ctil\tau)^{k+l}\,\dx^k\dy^l\vp_0 +
\mathcal{O}(\Dt^{m+1})
\label{eq:taylor_S0} \\ \nonumber
\frac{1}{2 \pi \tau} \int_{t_n}^{t_n + \tau} \!\!\!\! \int_0^{2 \pi}{\!\!\a^m\b^n\vp}\dd \theta
\dd \tilde t &= \sum\limits_{k,l=0}^{k+l \leq
m}\frac{\AvC{\a^{k+m}\b^{l+n}}}{k!l!}\frac{1}{\tau}\int_{t_n}^{t_n + \tau}
\hspace*{-0.5cm} \left(\tilde c(t_n + \tau - \tilde t )\right)^{k+l}\,\dx^k\dy^l
\tilde\vp_0 \dd \tilde t
\\ & \quad +\mathcal{O}(\Dt^{m+1}). \label{eq:taylor_Sb}}
Moreover, $\AvC{\a^{k+m}\b^{l+n}}=0$ if either $k$ or $l$ are odd integers, and
\mm{
\AvC{\a^2}=\AvC{\b^2}&=\frac12,\\
\quad\AvC{\a^4}=\AvC{\b^4}&=\frac38,\\
\quad\AvC{\a^2\b^2}&=\frac18.
}
\end{lemm}
\proof can be obtained by a direct evaluation. \hfill{\qed}
Using this Taylor expansion together with appropriate smoothness assumptions it is an
elementary exercise to prove that
\mm{\vector{h_1\\u_1\\v_1} &= \vector{h\\u\\v}\!\!(P)+\ttt.
}

\subsection{Approximate evolution operators}

In the previous section, we could give the integral representation \eqref{eq3} --
\eqref{eq5} a straightforward interpretation by deriving it from a Taylor expansion
along the streamlines.

However, it would be far too expensive to evaluate \eqref{eq3} -- \eqref{eq5} at
each quadrature node of the two-step finite volume scheme. In the present section we
derive the crucial approximations of  \eqref{eq3} -- \eqref{eq5}, leading to an
efficient and accurate algorithm.
This approximation comes in two steps. First we derive in Lemma~3.2 a suitable
approximation, which contains all terms necessary for the balance between the
pressure terms and the sources, i.e.\ $K_x=0=L_y$. This approximation is still
continuous. Afterwards, we apply a special numerical quadrature to approximate
the mantle integrals (i.e.~time dependent integrals), in order to obtain approximate
evolution operators which are explicit in time.

For the present paper, we are interested in second order schemes. It is therefore
sufficient that the predictor steps is first order accurate, i.e. accurate up to
terms of order $\tt$. In order to obtain a fully explicit first order approximation
of $(h,u,v)(P)$, we would like to convert the mantle integrals on the LHS of
(\ref{eq3})-(\ref{eq5}) into integrals over the sonic circle $S_0$. One possibility
would be, analogously to \cite{mathcom}, to use the simple rectangle rule
\mm{\label{eq:3} \frac{1}{2\pi \tau} \int_{t_n}^{t_n + \tau} \int_0^{2
\pi}{\vp(\tilde t,\theta)} \dd \theta \dd \tilde{t} &= \AvS0{\vp(t_n,
\theta)}+{\mathcal O}(\Dt). }

In the second step we can further eliminate derivatives over the sonic circle by
means of the per-partes formula, cf.~Lemma~2.1 in \cite{mathcom}
\mm{\label{eq:4} \frac{\ctil \tau}2\AvS0{(\vp_x+\psi_y)} &= \AvS0{(\alpha \vp + \beta
\psi)}+\tt. }

However, it has been shown in \cite{mathcom,sisc,stab} that the
application of classical quadrature rules, such as the rectangle rule in
(\ref{eq:3}), are not well suited for approximation of discontinuous waves, which
may propagate along the mantle of the bicharacteristic cone. It resulted in a
reduced stability range of the FVEG.  In particular,  if the mantle integrals are
approximated by the rectangle rule the CFL stability number was 0.63 and 0.56 for
the first and second order FVEG scheme, respectively; it is the so-called FVEG3
scheme, cf.\ \cite{eccomas}. In the recent paper \cite{sisc} new quadrature rules
have been proposed for the mantle integrals.  For example, if $f = f(x)$ is a
piecewise constant function, then it was shown in \cite{sisc} that
\begin{equation}
\label{eq:sgn} \frac{1}{2 \pi} \int_{0}^{2 \pi} f(Q) \cos \theta \dd \theta +
\frac{1}{2 \pi} \int_{t_n}^{t_n + \tau}\frac{1}{t_n + \tau - \tilde t} \int_0^{2
\pi} f(\tilde Q) \cos \theta \dd \tilde t = \frac{1}{2 \pi} \int_{0}^{2 \pi} f(Q)
\mbox{ sgn} \cos \theta \dd \theta,
\end{equation}
an analogous relation holds for $f(Q) \sin \theta$. Similarly, the quadrature rules
for bilinear data have been derived, cf.\ Lemma~A.1 in the Appendix of \cite{sisc}.

As a result new approximate evolution operators evaluate exactly
each planar wave propagating either in $x-$ or $y-$
directions and increase the stability range of the FVEG scheme substantially yielding the CFL number close to~1.

Before we apply the quadrature rules proposed in \cite{sisc}, we approximate and simplify the
exact integral equations (\ref{eq3})-(\ref{eq5}). This is done in Lemma~\ref{lema32}. Our
strategy is to drop as many of the second order terms as possible, but to keep all those
terms which enter the balance of convective fluxes and source terms, i.e. $K_x = 0 = L_y$,
and are therefore needed for well-balancing. The remaining terms will be reformulated or
approximated up to the order ${\mathcal O}(\Dt^2)$ in such a way, that the above mentioned
quadrature rules from \cite{sisc} can be applied. Thus, we keep the balance conditions for
source terms and at the same time we will be able to approximate all resulting mantle
integrals in a stable way.
\begin{lemm}
\label{lema32} The following operator is a first order approximation of the exact
integral equations \eqref{eq3}-\eqref{eq5}
\begin{eqnarray}
\nonumber
  h\left( P\right) &=& - b(P) + \frac{1}{2\pi}\int_{0}^{2\pi}
   \left( h\left(Q\right)  + b\left(Q\right) \right)
   - \frac{\tilde{c}}{g} \left( u\left(Q\right)\cos\theta +
    v\left(Q\right)\sin\theta \right)\dd\theta    \\
      \label{eq6}
&&- \frac{1}{2\pi} \int_{t_n}^{t_n + \tau} \frac{1}{t_n + \tau - \tilde{t}}
\int_{0}^{2 \pi}\frac{\tilde c}{g} \left( u(\tilde Q) \cos \theta + v(\tilde Q) \sin
\theta \right)  \dd \theta \dd \tilde{t}
\\ \nonumber
&& + \frac{1}{2 \pi}  \int_{t_n}^{t_n + \tau}  \int_0^{2 \pi}  \left( \tilde u b_x
(\tilde Q) + \tilde v b_y (\tilde Q)  \right) \dd \theta \dd \tilde t + {\mathcal
O}\left(\Delta t^2 \right),
\end{eqnarray}
\begin{eqnarray}
\nonumber u\left(P \right) &=&    \frac{1}{2 \pi} \int_{0}^{2\pi}
-\frac{1}{\tilde{c}} K\left(Q\right) \cos\theta
   + u \left(Q \right)\cos^2\theta  + v\left(Q \right)\sin\theta\cos\theta \, \dd\theta   \\
&& +\frac{1}{2}u\left( Q_0\right) - \frac{1}{2\pi} \int_{t_n}^{t_n + \tau}
\frac{1}{t_n + \tau - \tilde{t}} \int_{0}^{2 \pi} \frac{1}{\tilde{c}} K (\tilde Q)
\cos \theta \, \dd \theta \dd \tilde{t}  \label{eq7}
\\ \nonumber
&& + \frac{1}{2\pi} \int_{t_n}^{t_n + \tau} \frac{1}{t_n + \tau - \tilde{t}}
 \int_{0}^{2 \pi} \left( u(\tilde Q) \cos 2 \theta + v(\tilde Q) \sin  2\theta \right)
 \dd \theta \dd \tilde{t}   + {\mathcal O}\left(\Delta t^2 \right),
\end{eqnarray}
\begin{eqnarray}
\nonumber v\left(P \right) &=&  \frac{1}{2 \pi} \int_{0}^{2\pi}
   -\frac{1}{\tilde{c}} L\left( Q \right) \sin\theta
   + u \left(Q \right)\sin \theta \cos\theta  + v\left(Q \right)\sin^2\theta \, \dd\theta   \\
&& +\frac{1}{2}v\left( Q_0\right)  - \frac{1}{2\pi} \int_{t_n}^{t_n + \tau}
\frac{1}{t_n + \tau - \tilde{t}} \int_{0}^{2 \pi} \frac{1}{\tilde{c}} L(\tilde Q)
\sin \theta \,   \dd \theta \dd \tilde{t} \label{eq8}
\\ \nonumber && +
\frac{1}{2\pi} \int_{t_n}^{t_n + \tau} \frac{1}{t_n + \tau - \tilde{t}}
 \int_{0}^{2 \pi} \left( u(\tilde Q) \sin 2 \theta -  v(\tilde Q) \cos  2\theta \right)
 \dd \theta \dd \tilde{t}   + {\mathcal O}\left(\Delta t^2 \right).
\end{eqnarray}
\end{lemm}
The proof of this lemma is postponed to Appendix~B.  Now, in order to obtain time
explicit approximate evolution operators we approximate time integrals in
(\ref{eq6})-(\ref{eq8}).

The only integral which is not of the form (\ref{eq:sgn}) is the last term in
(\ref{eq6}). Here we apply the rectangle rule (\ref{eq:3}) and get
$$
\frac{1}{2 \pi} \int_{t_n}^{t_n + \tau} \int_0^{2 \pi} \left( \tilde u b_x (\tilde
Q) + \tilde v b_y (\tilde Q) \right) \dd \theta \dd \tilde t =
 \frac{\tau}{2 \pi} \int_0^{2 \pi} \left(\tilde u b_x (Q) + \tilde v b_y
(Q)  \right) \dd \theta + {\mathcal O} \left(\Delta t^2 \right).
$$
Moreover, for the special case when the bottom topography  slopes $b_x, b_y$ are
approximated by a piecewise constant functions, which is the case of our bilinear
recovery, for example, we can evaluate this term exactly. Note, that $\tilde u =
\const.,$ $\tilde v = \const.,$ and $b = b(x,y)$ does not change in time. Thus, we
have
$$
\frac{1}{2 \pi} \int_{t_n}^{t_n + \tau} \int_0^{2 \pi} \left( \tilde u b_x (\tilde
Q) + \tilde v b_y (\tilde Q) \right) \dd \theta \dd \tilde t = \frac{\tau}{2 \pi}
\int_0^{2 \pi} \left( \tilde u b_x(Q) + \tilde v b_y (Q) \right) \dd \theta.
$$

All other integrals are of the form (\ref{eq:sgn}), so we can apply numerical
quadratures from \cite{sisc}. They will be used separately for
constant and bilinear approximations. Thus, using (\ref{eq:sgn}) we get, analogously
to \cite{sisc}, the approximate evolution operator $E_{\Delta}^{const}$ using the
{piecewise constant approximate functions}
\begin{eqnarray}
\nonumber
   h\left( P\right) &=& - b(P) + \frac{1}{2\pi}\int_{0}^{2\pi}\left[
   ( h\left(Q\right) + b(Q) ) - \frac{\tilde{c}}{g}\left( u\left(Q\right)\sgn(\cos\theta)
   + v\left(Q\right)\sgn(\sin\theta)
   \right) \right] \dd\theta  \\
   && + \frac{\tau}{2 \pi} \int_0^{2 \pi} \left( \tilde u b_x (Q) +
\tilde v b_y (Q)  \right) \dd \theta \nonumber\\ \nonumber
   u\left(P \right) &=&\frac{1}{2\pi}\int_{0}^{2\pi}\left[
   -\frac{1}{\tilde{c}} K \left(Q\right) \sgn(\cos\theta)
      +u\left(Q\right)\left(\cos^2\theta + \frac12 \right)
       +v\left(Q\right)\sin\theta\cos\theta \right] d\theta \\
  \nonumber
   v\left(P \right) &=&\frac{1}{2\pi}\int_{0}^{2\pi}\left[
   -\frac{1}{\tilde{c}} L\left(Q\right) \sgn(\sin\theta)
      + u\left(Q\right)\sin\theta\cos\theta
      +v\left(Q\right)\left(\sin^2\theta + \frac12 \right)\right] d\theta.
\\ \label{approx_pwc}
\end{eqnarray}
For the {piecewise bilinear ansatz functions} the justification of the mantle
integrals approximation is more involved. The reader is referred to \cite{sisc} for
more details.  Applying the approximations from \cite{sisc} the approximate
evolution operator $E_{\Delta}^{bilin}$ reads
 \begin{eqnarray}
 \nonumber
    h\left( P\right) &=& -b(P) + h(Q_0) + b(Q_0)
    + \frac14\int_0^{2\pi} (h(Q) - h(Q_0)) + (b(Q) - b(Q_0))
    \dd \theta \\ \nonumber
    && - \frac{1}{\pi}\int\limits_0^{2\pi}
    \left[\frac{\tilde{c}}{g}u(Q)\cos\theta + \frac{\tilde{c}}{g}v(Q)\sin\theta\right]
\dd \theta  + \frac{\tau}{2 \pi} \int_0^{2 \pi} \left( \tilde u b_x (Q)
 + \tilde v b_y (Q) \right) \dd \theta\\
 \nonumber
    u\left(P \right) &=& u(Q_0) - \frac{1}{\pi}\int_0^{2\pi}
    \frac{1}{\tilde{c}}  K(Q)  \cos\theta \dd \theta
 \\
\nonumber
    && +\frac{1}{4}\int_0^{2\pi}
    \left[3u(Q)\cos^2\theta + 3v(Q)\sin\theta\cos\theta - u(Q) - \frac{1}{2} u(Q_0)
    \right] \dd\theta \\ \nonumber
    v\left(P \right) &=& v(Q_0) - \frac{1}{\pi}\int_0^{2\pi}
    \frac{1}{\tilde{c}}  L(Q)  \sin\theta \dd \theta\\
\label{approx_pwb}
    &\ &
    +\frac{1}{4}\int_0^{2\pi}
    \left[3u(Q)\sin\theta \cos\theta + 3v(Q)\sin^2\theta - v(Q) - \frac{1}{2} v(Q_0)
    \right] \dd\theta.
\end{eqnarray}
The approximate evolution operators (\ref{approx_pwc}), (\ref{approx_pwb}) together
with the finite volume update (\ref{so_6}) define the FVEG schemes. We will
summarize the algorithm in the Section~4.

Let us pause for a moment and discuss the possible advantages of using the
approximate evolution operators instead of the Taylor expansion \eqref{eq:huv_1} in
the predictor step entering the finite volume update \eqref{so_6}.
\begin{itemize}
\item The above approximate evolution operators does not rely upon any derivative of
the unknowns $(h,u,v)$, while the Taylor expansion uses $h_x, h_y, u_x+v_y$.
Therefore (\ref{approx_pwc}), (\ref{approx_pwb}) are less dependent upon the
reconstruction of the piecewise constant solution than \eqref{eq:huv_1}.
\item The integrations along the cone in
(\ref{eq3})-(\ref{eq5}) take into account the whole domain of dependence.
This may result in a more robust algorithm, in particular for discontinuous solutions.
\end{itemize}

\subsection{Well-balanced property of the approximate evolution operators}

The aim of this subsection is to verify that the approximate evolution operators
\eqref{approx_pwc}, \eqref{approx_pwb} are well-balanced for the lake at rest as well as for
the jet in the rotational frame. The proof of Theorem \ref{theorem:well-balanced-eg} below is
based on the sufficient conditions for well-balancing formulated in Theorem
\ref{theorem:well-balanced-fv}.
\begin{theorem}\label{theorem:well-balanced-eg}
Suppose that the reconstructions at time $t^n$ satisfy for all $(x,y)$
\begin{eqnarray}
u^n(x,y) &\equiv& 0  \label{cc1}\\
\dy v^n(x,y) &\equiv& 0  \label{cc2} \\
\mu_x v^n(x,y) \dy  h^n(x,y) &\equiv& 0 \label{cc3}  \\
\dx K^n(x,y) &\equiv& 0 \label{cc4} \\
\dy L^n(x,y) &\equiv& 0.   \label{cc5}
\end{eqnarray}
Then the approximate EG predictor steps defined by \eqref{approx_pwc},
\eqref{approx_pwb} satisfy the conditions for well-balancing of Theorem {\rm
\ref{theorem:well-balanced-fv}.}

Therefore, the FVEG schemes based on the above approximate evolution operators are
well-balanced for the lake at rest and the jet in the rotational frame.

\end{theorem}

\proof We prove here that the approximate evolution operator for piecewise constant
data $E_{\Delta}^{const}$, see\ (\ref{approx_pwc}), satisfies conditions
(\ref{c1})-(\ref{c5}) of Theorem \ref{theorem:well-balanced-fv}. The proof for the
approximate evolution operator for piecewise bilinear data $E_{\Delta}^{bilin}$, see
(\ref{approx_pwb}), is analogous.

%
%

First we use conditions \eqref{cc1} -- \eqref{cc5} to simplify the approxmate evolution
operator \eqref{approx_pwc}. Due to \eqref{cc1} and \eqref{cc2}
\mmn{\frac{1}{2\pi}\int_{0}^{2\pi}
   u\left(Q\right)\sgn(\cos\theta)\dd\theta =
   \frac{1}{2\pi}\int_{0}^{2\pi}
   v\left(Q\right)\sgn(\sin\theta)\dd\theta = 0}
Similarly
\mmn{\frac{1}{2\pi}\int_{0}^{2\pi}K\left(Q\right)\sgn(\cos\theta)\dd\theta =
   \frac{1}{2\pi}\int_{0}^{2\pi}L\left(Q\right)\sgn(\sin\theta)\dd\theta =
   \frac{1}{2\pi}\int_{0}^{2\pi}v\left(Q\right)\sin\theta\cos\theta\dd\theta = 0,}
while
\mmn{\frac{1}{2\pi}\int_{0}^{2\pi}v\left(Q\right)
\left(\sin^2\theta+\frac12\right)\dd\theta = \frac12(v^L+v^R).}
Due to \eqref{cc5}, \eqref{c1} and \eqref{cc3}
\mm{\frac{1}{2\pi}\int_{0}^{2\pi}\tilde v b_y\left(Q\right)\dd\theta =
-\frac{1}{2\pi}\int_{0}^{2\pi}\tilde v h_y\left(Q\right)\dd\theta = 0.}
Using the above identities in \eqref{approx_pwc} gives the simplified approximate evolution
operator, valid for the jet in the rotational frame
\begin{eqnarray}
   \hat h\left( P\right) &=& - b(P) + \frac{1}{2\pi}\int_{0}^{2\pi}
   ( h\left(Q\right) + b(Q) )\dd\theta \nn%
   \hat u\left(P \right) &=& 0 \nn
   \hat v\left(P \right) &=& \frac12(v^L+v^R).
\label{approx_pwc_jet}
\end{eqnarray}
From this, \eqref{c1}, \eqref{c2} and \eqref{c5} follow immediately. To verify \eqref{c3},
we set $P_0:=P-(0,0,\tau)$ (the projection of $P$ onto the plane $t\equiv t_n$) and
compute
\mmn{ \hat v\left( P\right)\partial_y\hat h\left( P\right) &= \hat v\left(
P\right)\partial_y\left\{ - b(P) + \frac{1}{2\pi}\int_{0}^{2\pi} ( h\left(Q\right) +
b(Q) )\dd\theta \right\} \nn &= -\hat v\left( P\right)\partial_y b(P_0) \nn &= -\hat
v\left( P\right)\partial_y \left((h(P_0)+b(P_0))-h(P_0)\right) \nn &= \hat v\left(
P\right)\partial_y h(P_0) \nn &= 0.}

It remains to prove that \eqref{c4} holds. Since $K=g(h+b-V)$,
\begin{eqnarray}
\hat K(P) &=& g\left(\hat h(P) + b(P) -  V^{n+1/2}(P)\right) \nonumber \\
\nonumber &=& \frac{g}{2 \pi} \int_0^{2 \pi} ( h^n(Q) + b(Q) - V^n(Q)) + \frac{g}{2
\pi} \int_0^{2 \pi} V^n(Q) \dd \theta  -  g\;V^{n+1/2}(P), \nonumber  \\ \nonumber
&=& \frac{1}{2 \pi} \int_0^{2 \pi} K^n(Q) \dd \theta + \frac{g}{2 \pi} \int_0^{2
\pi} V^n(Q) \dd \theta - g\;V^{n+1/2}(P).
\end{eqnarray}
Differentiating this equation with respect to $x$ and applying (\ref{UV}) implies
\begin{eqnarray}
\nonumber
\dx \hat K(P) &=& \frac{1}{2 \pi} \int_0^{2 \pi} \dx K^n(Q) \dd \theta +
\frac{g}{2 \pi} \int_0^{2 \pi} \dx V^n(Q) \dd \theta - g\;\dx V^{n+1/2}(P)
\\ \nonumber
&=& \frac{f}{2 \pi} \int_0^{2 \pi} v^n(Q) \dd \theta - f \hat v(P)
\\ \nonumber
&=& f \frac{v^R + v^L}{2} - f \frac{v^R + v^L}{2} = 0,
\end{eqnarray}
which is the well-balanced condition (\ref{c4}) and concludes the proof. \hfill
\qed

\section{Summary of the FVEG algorithm}
In this section we summarize the main steps of the FVEG method by presenting the
algorithm for the first and second order scheme including the effects of bottom
topography as well as the Coriolis forces.\\

{\bf Algorithm}
\begin{itemize}
\item[1] Given are piecewise constant approximations at time $t_n$:
$h_{ij}^n, u_{ij}^n, v_{ij}^n, i,j \in \Bbb Z$, the bottom topography  $b(x,y)$,
mesh and time steps $ \hbar, \Delta t$ and constants $g,\, f$; compute
\begin{eqnarray}
\nonumber &&b_{ij}^n = b(x_i, y_j, t^n),\\ \nonumber
&& V_{ij}^n =  \frac{f}{g} \hbar
\sum_{i'=i_0}^i \frac{v_{i'-1,j}^n + v_{i'j}^n}{2},
\\ \nonumber
&& U_{ij}^n = \frac{f}{g} \hbar \sum_{j'=j_0}^j \frac{u_{i,j'-1}^n + u_{ij'}^n}{2}.
\end{eqnarray}
\item[2] {\sl recovery step:} \\
If the scheme is second order, do the recovery step. For smooth parts of solution
apply the continuous bilinear recovery, cf.\ \cite{jcp}. Possible overshoots on
discontinuities are limited, e.g.\ by the minmod limiter; cf.\ \cite{jcp}. This
yields the piecewise bilinear approximations $R_{\hbar} h^n,\, R_\hbar u^n,\,
R_\hbar v^n,\, R_\hbar b^n$, $R_\hbar U^n$, $R_\hbar V^n$.
\item[3]
{\sl predictor step / approximate evolution:} \\
Compute the intermediate solutions at time level $t_{n+1/2}$ on the cell interfaces
by the approximate evolution operators. For the first order scheme use the
approximate evolution operator $E_\Delta^{const}$ (\ref{approx_pwc}); the second
order scheme is computed using both approximate evolution operators
$E_\Delta^{const}$ (\ref{approx_pwc}) as well as $E_\Delta^{bilin}$
(\ref{approx_pwb}), cf.\ (\ref{so_77}). Integration along the cell interfaces is
realized numerically by the Simpson rule.
\item[4]
{\sl corrector step / FV-update:} \\
Compute the Coriolis forces and  the bottom topography at the intermediate time
level $t_{n+1/2}$ and at each integration points on cell interfaces, i.e.\ at
vertices and midpoints:
\begin{eqnarray}
\nonumber
&&b_{k \ell}^{n+1/2} = b(x_k,y_{\ell}), \quad k=i, i\pm 1/2,\ \ell=j, j\pm 1/2; \\
&& V_{i+1/2, {\ell}}^{n+1/2} =  \frac{f}{g} \hbar \sum_{i'=i_0}^i
\frac{v_{i-1/2,{\ell}}^{n+1/2} + v_{i+1/2, {\ell}}^{n+1/2}}{2}, \quad \ell=j, j\pm
1/2; \nonumber \\ \nonumber
&& U_{k,j+1/2}^{n+1/2} = \frac{f}{g} \hbar \sum_{j'=j_0}^j
\frac{u_{k,j-1/2}^{n +1/2} + u_{k j+1/2}^{n+1/2}}{2}, \quad k=i, i\pm 1/2.
\end{eqnarray}
Do the FV-update (\ref{so_6}) using the well-balanced approximation of the source
terms (\ref{source}).
\end{itemize}

\section{Numerical experiments}\label{section:num_exp}
One interesting steady state, which should be correctly resolved by a well-balanced
scheme,  is the stationary steady state, i.e. $h+b=\const.$ and $u=0=v.$ In this
section we demonstrate well-balanced behaviour of the proposed FVEG schemes through
several benchmark problems for stationary and quasi-stationary states, i.e.\ $h+b
\approx \mbox{ const.}$ and $u \approx 0 \approx v$; see \cite{kurganov,lv1}
for related results in literature. Further, we present results for steady jets
including effects of the Coriolis forces and show that the FVEG scheme is
well-balanced also for this nontrivial steady state. At the end of this section we
compare accuracy and computational time of the well-balanced FVEG method and the
well-balanced  second and fourth order  FVM of Audusse et al. \cite{bouchut1}
and Noelle et al.\ \cite{noelle1}.

\subsection{One-dimensional stationary and quasi-stationary states}
In this experiment we have tested  the preservation of a stationary steady state as
well as the approximation of small perturbations of this steady state. The bottom
topography consists of one hump
\begin{eqnarray}\nonumber
  \label{QuasiStat_1d_a}
  b(x)& = & \left\{
  \begin{array}  {ll}
    0.25(\cos(10\pi(x-0.5)) + 1) & \mathrm{if}\  |x-0.5|<0.1\\
    0 & \mathrm{otherwise}
  \end{array}  \right. \\
  \nonumber
  \mathrm{and\ the\ initial\ data\ are  }\  u(x,0) &=& 0,\\
  \nonumber
  h(x, 0)& = &\left\{
  \begin{array}  {ll}
    1 - b(x) + \varepsilon & \mathrm{if}\  0.1 < x < 0.2\\
    1 - b(x) & \mathrm{otherwise.}
  \end{array}\right. 
  \label{QuasiSteady_large}
\end{eqnarray}
The parameter $\varepsilon$ is chosen to be $0$, $0.2$ or $0.01$. The computational
domain is the interval $[0,1]$ and absorbing boundary conditions have been
implemented by extrapolating all variables. The gravitational constant $g$ was set
to 1 analogously as in \cite{kurganov,lv1}.   It should be pointed out that
the one-dimensional problems are actually computed by a two-dimensional code by
imposing zero tangential velocity $v=0.$

Firstly we test the ability of the FVEG scheme to  preserve the stationary steady
state, i.e. the lake at rest case, by taking $\varepsilon = 0$. In Table~1 the
$L^1$-errors for different times computed with the first order FVEG method, cf.\
(\ref{approx_pwc}), and with the second order FVEG method, cf.\ (\ref{approx_pwb}),
are presented. Although we have used a rather coarse mesh consisting of $20 \times
20$ mesh cells, it can be seen clearly that the FVEG scheme balances up to the
machine accuracy also for long time computations.

\begin{table}[ht]
\caption{ The $L^1$-error of the well-balance FVEG scheme using $20 \times 20$ mesh
cells. }

\begin{center} \footnotesize
\begin{tabular}{|c|c|c|c|c|} \hline
 Method&  $t=0.2$   & $ t=1 $& $t=10$ \\
 \hline
 first order FVEG &$ 0$ &$ 0$&$2.22\times 10^{-16}$\\
second order FVEG &$ 1.67 \times 10^{-17}$ &$ 1.11\times 10^{-17}$& $4.27\times 10^{-16}$ \\
\hline
\end{tabular}
\end{center}
\end{table}

\begin{figure}[ht]
\begin{center}
     \epsfig{file=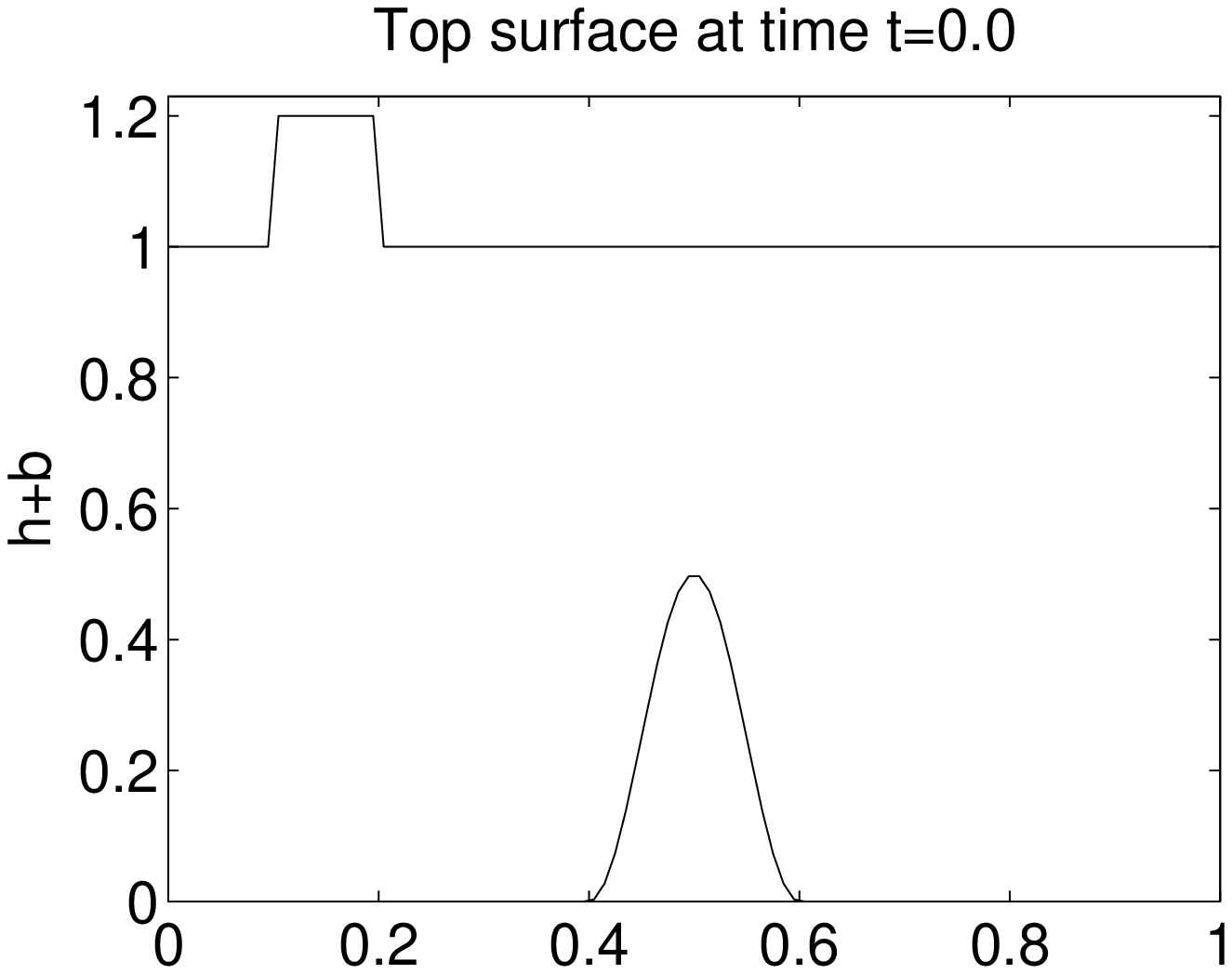,width=0.49\textwidth}
     \epsfig{file=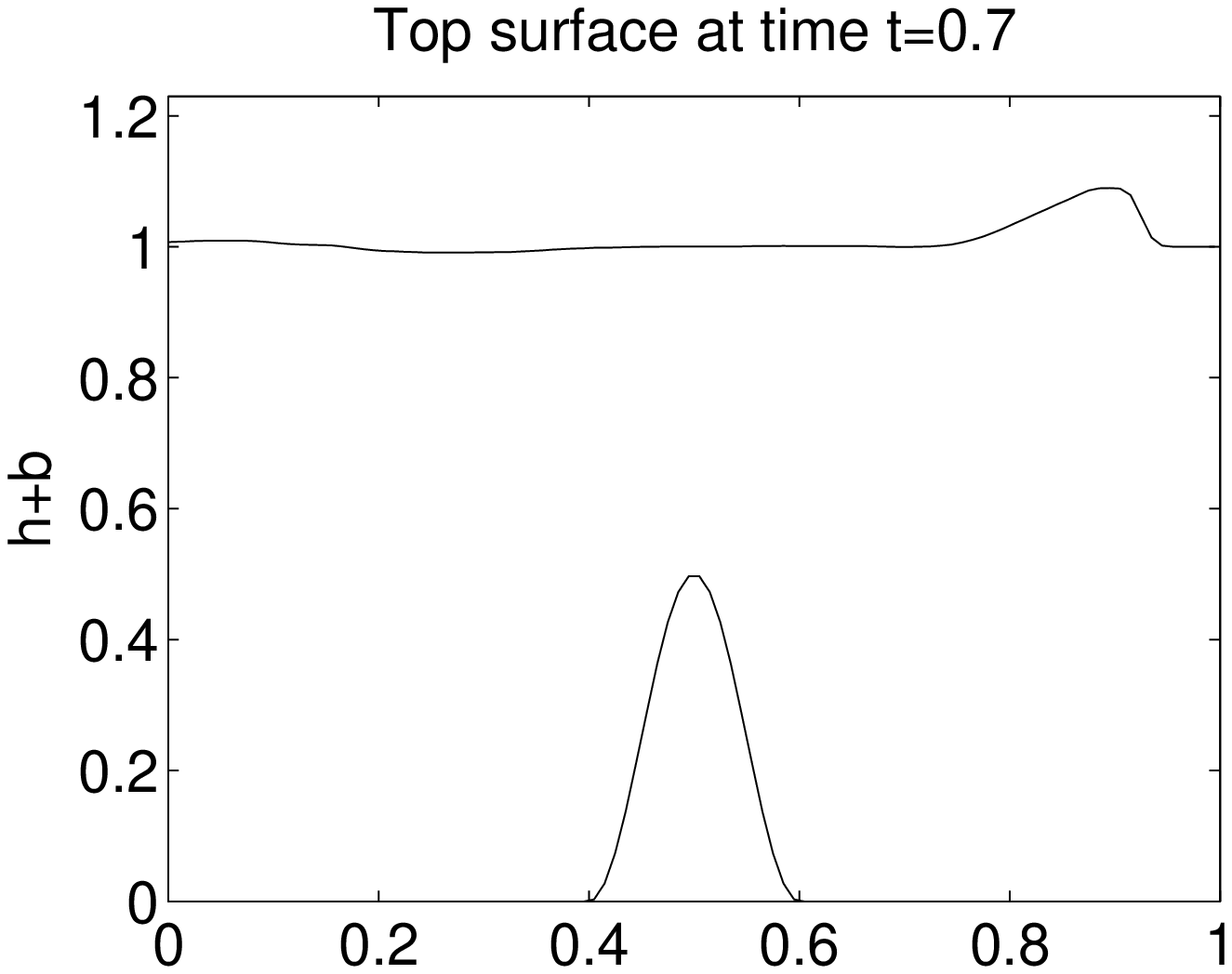,width=0.49\textwidth}
\end{center}
\caption{Propagation of small perturbations, $\varepsilon=0.2$.}
\end{figure}

\begin{figure}[ht]
\begin{center}
   \epsfig{file=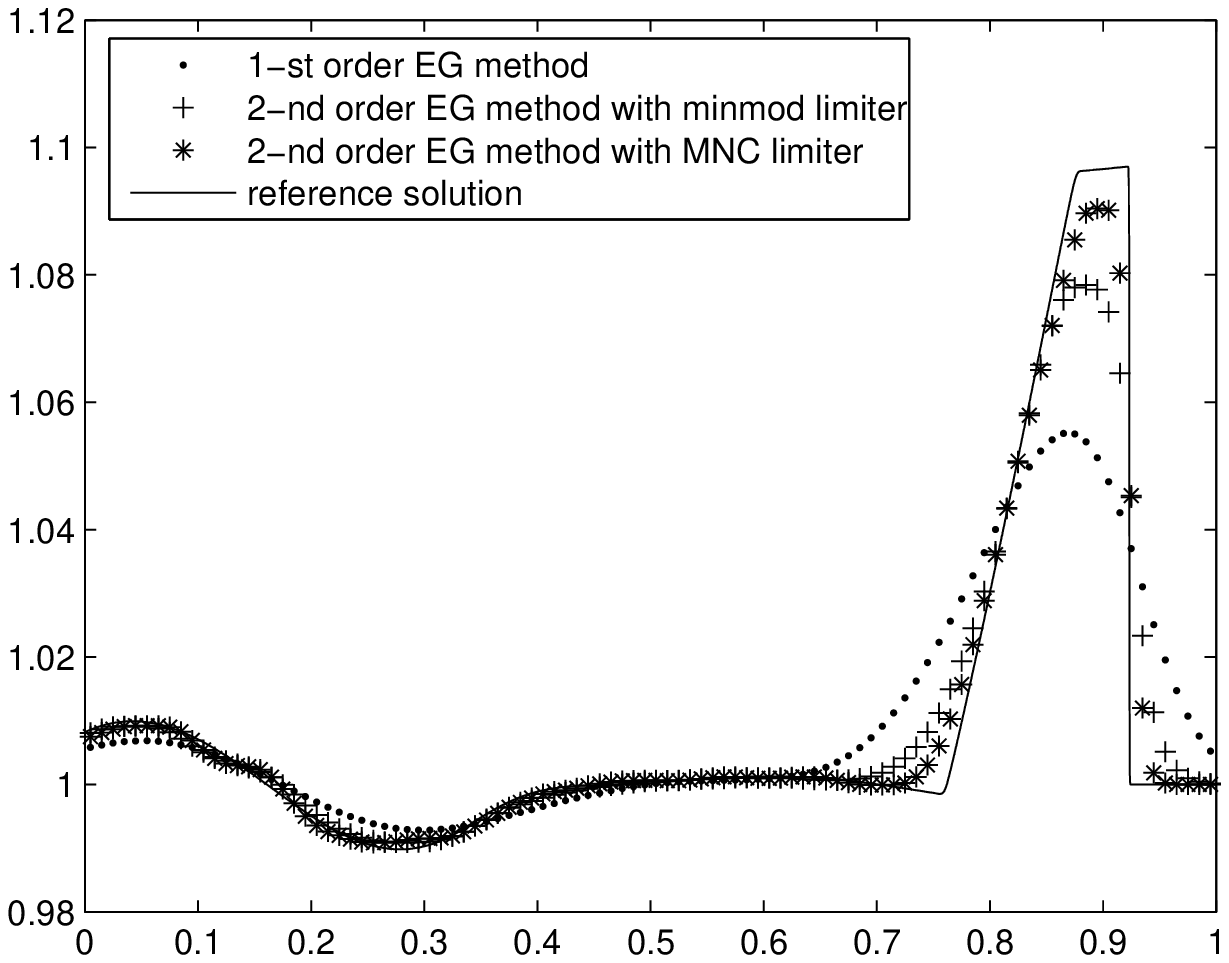, width=0.49\textwidth}
   \epsfig{file=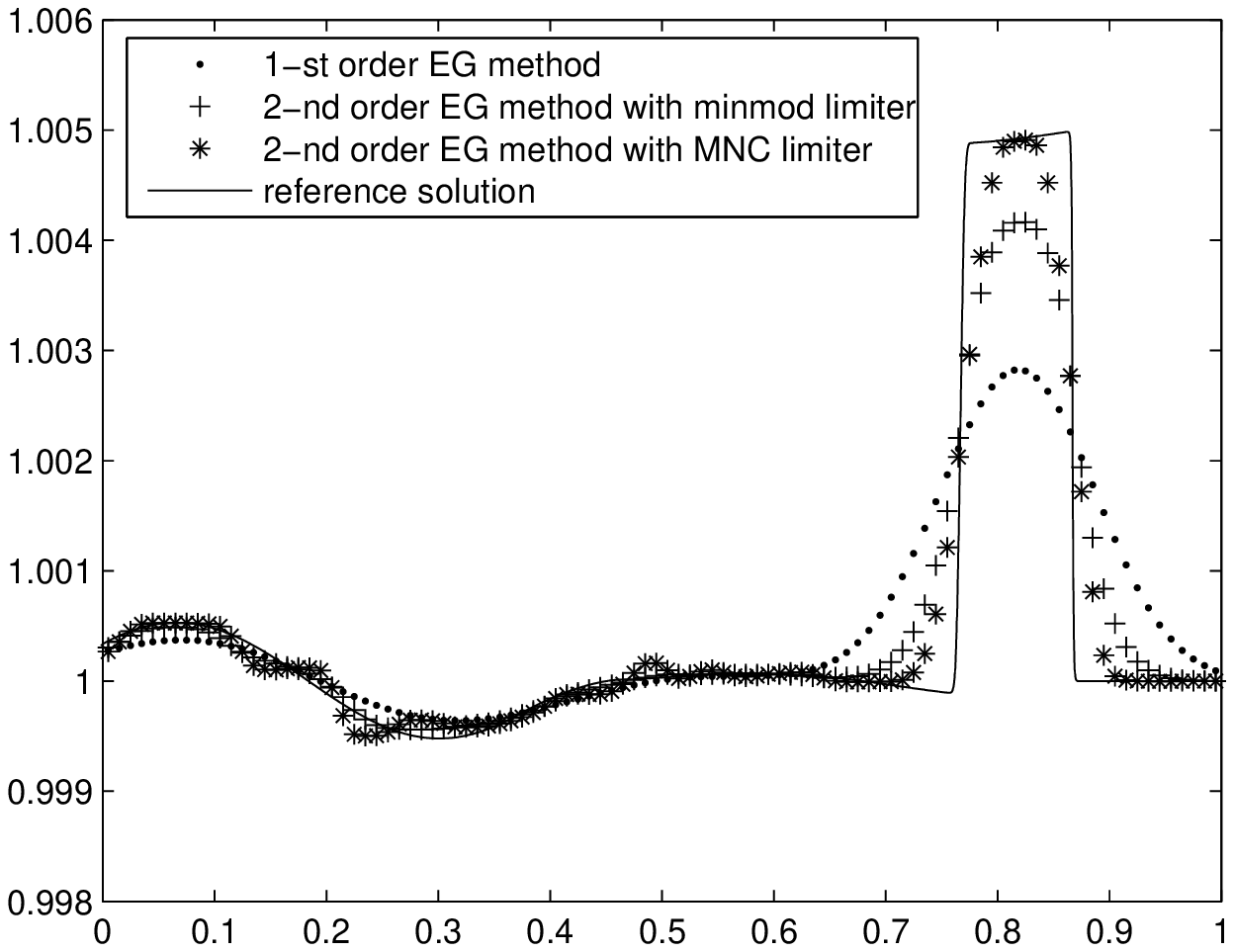, width=0.49\textwidth}
\end{center}
\label{QuasiSteady} \caption{Propagation of small perturbations, magnified view;
$\varepsilon = 0.2$ (left) and  $\varepsilon = 0.01$ (right).}
\end{figure}

In  Figure~2
the typical propagation of small height perturbations is shown at time $t=0.7$.
The solution is computed on a mesh with $100 \times 5$
cells and the height of the initial perturbation was $\varepsilon=0.2$.  The initial
disturbance generates two waves, the left-going wave runs out of the computational
domain, and the right-going wave  passes the bottom elevation obstacle. It is known
that if the perturbations are relatively large in comparison to the discretization
error a ``naive"  approximation of the source term, i.e.\ not well-balanced scheme,
e.g. the fractional step method, can still yield reasonable approximations. However,
for small perturbations, i.e.\ $\varepsilon$ of order of the discretization errors,
such a scheme would yield strong oscillations over the bottom hump and the wave of
interest will be lost in the noise, see \cite{lv1}.

In Figure~3
we compare results for water depth $h$ at time $t=0.7$ obtained by the first and
second order FVEG methods using the minmod limiter and the monotonized centered
limiter (denoted as MNC), respectively. In the left picture  $\varepsilon =0.2$, the
right picture shows results for $\varepsilon=0.01.$  The reference solutions was
obtained by the second order FVEG method with the minmod limiter on a mesh with
10000 cells. For the first order scheme and the second order scheme with minmod
limiter we can notice correct resolution of small perturbations of the stationary
steady state even if the perturbation is of the order of the truncation error. The
MNC limiter resolves the peak much more sharply, but overcompress the left-going
wave. This is a well-known feature of compressive limiters, see e.g.
the discussion in \cite{sweby,noelle-lie,kurganov_2}.

\subsection{Two-dimensional quasi-stationary problem}
The second example is a two-dimensional analogue of the previous one. The bottom
topography is given by the function
\begin{eqnarray}
  \label{QuasiStat_2d_a}
  b(x,y) &=& 0.8\,\mathrm{exp}\left(-5\left(x-0.9\right)^2 - 50\left(y-0.5\right)^2\right)
\end{eqnarray}
and\ the\ initial\ data\ are
\begin{eqnarray}
  \nonumber
  h(x,y, 0) &=& \left\{
  \begin{array}  {ll}
    1 - b(x,y) + \varepsilon & \mathrm{if}\  0.05 < x < 0.15\\
    1 - b(x,y) & \mathrm{otherwise}
  \end{array}\right. \\
  \label{QuasiStat_2d_b}
  u(x,y,0) &=& v(x,y,0) = 0.
\end{eqnarray}
The parameter $\varepsilon$ is set to $0$ and $0.01$. The computational domain is
$[0,2]\times[0,1]$ and the absorbing extrapolation boundary condition are used.

First, we take $\varepsilon=0$ and test the preservation of a two-dimensional lake
at rest on a mesh with $20 \times 20$ mesh cells, see~Table~2. Analogously to the
one-dimensional case this steady state is preserved up to the machine accuracy.  In
the Figure~4
we present two solutions of a perturbed problem, $\varepsilon=0.01$, which are
computed on a $200\times100$ grid (left) and on a $600 \times 300$ grid (right) by
the second order FVEG scheme with the minmod limiter. Notice that the FVEG method
correctly approximates small perturbed waves, the perturbation propagates over the
bottom hump without any oscillations. Note that the wave speed is slower over the
hump, which leads to a distortion of the initially planar perturbation. The
perturbed wave  runs out of the computational domain and the flat surface is
obtained at the end. Our results are in a good agreement with other results
presented in literature, cf., e.g., \cite{kurganov,lv1,noelle1,XingShu2005}.


\begin{table}
\caption{ The $L^1$-error of the well-balance FVEG scheme using $20 \times 20$ mesh
cells. }
\begin{center} \footnotesize
\begin{tabular}{|c|c|c|c|c|} \hline
 Method&  $t=0.2$   & $ t=1 $& $t=10$ \\
 \hline
 first order FVEG &$ 2.35 \times 10^{-17}$ &$ 5.09 \times 10^{-17}$&$5.02\times 10^{-17}$\\
second order FVEG &$ 4.97 \times 10^{-17}$ &$ 6.74\times 10^{-17}$& $1.53\times 10^{-16}$ \\
\hline
\end{tabular}
\end{center}
\end{table}

\begin{figure}[Ht]
  \begin{center}
    \epsfig{file=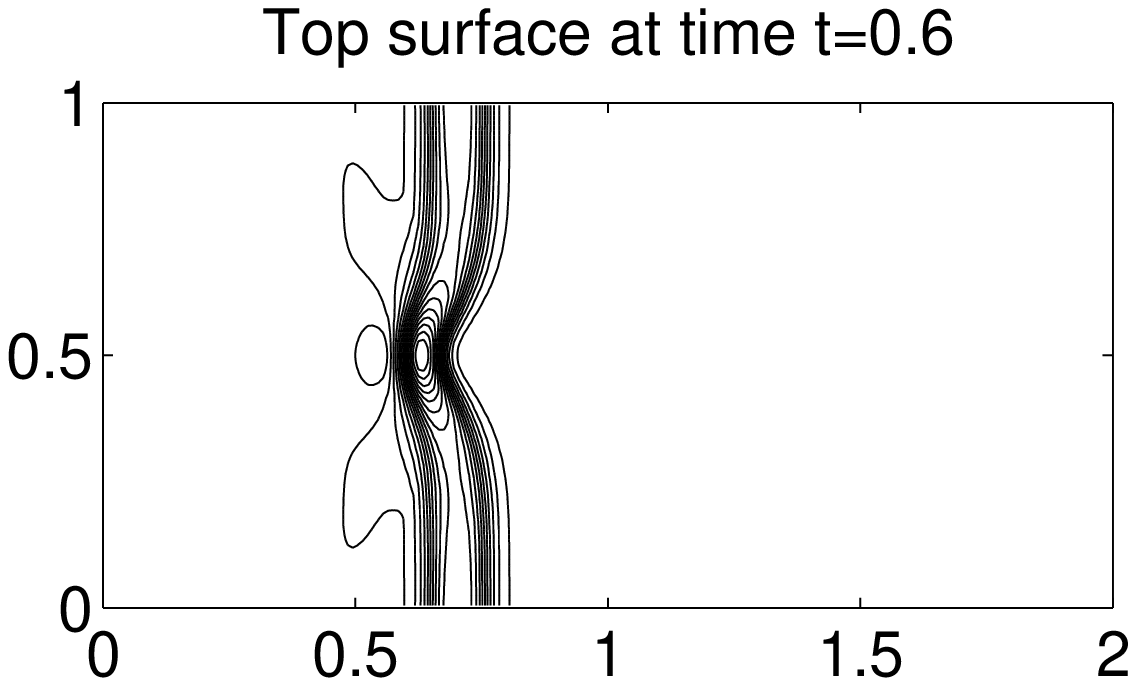, width=0.33\textwidth}
    \epsfig{file=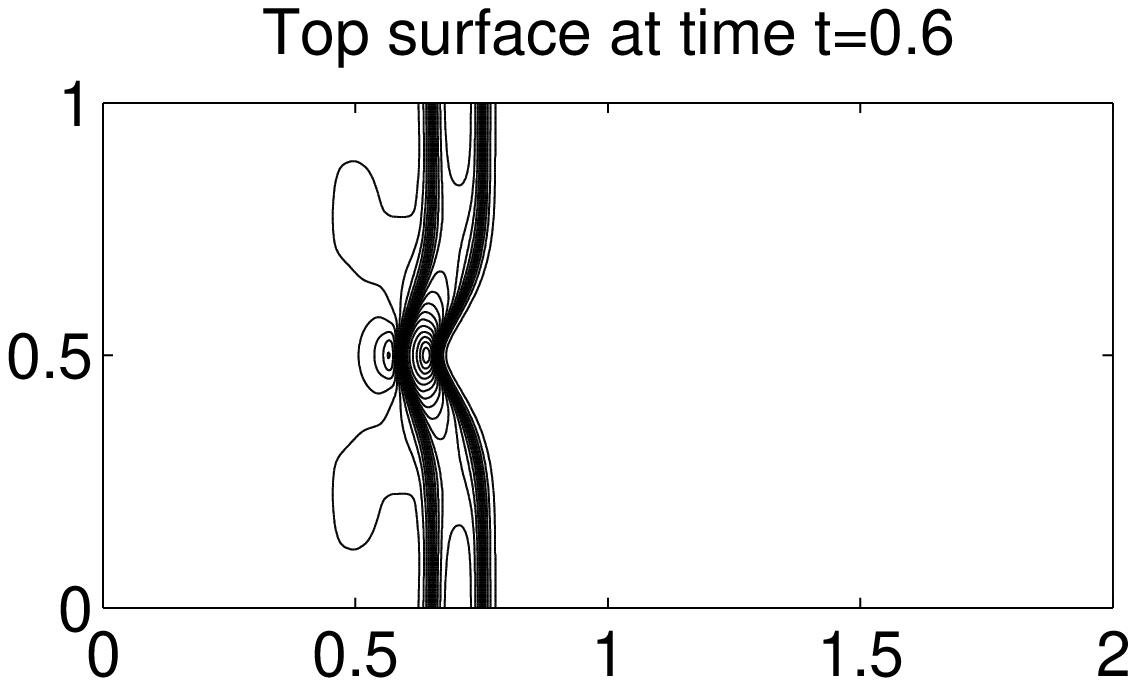,width=0.33\textwidth}\\
    \epsfig{file=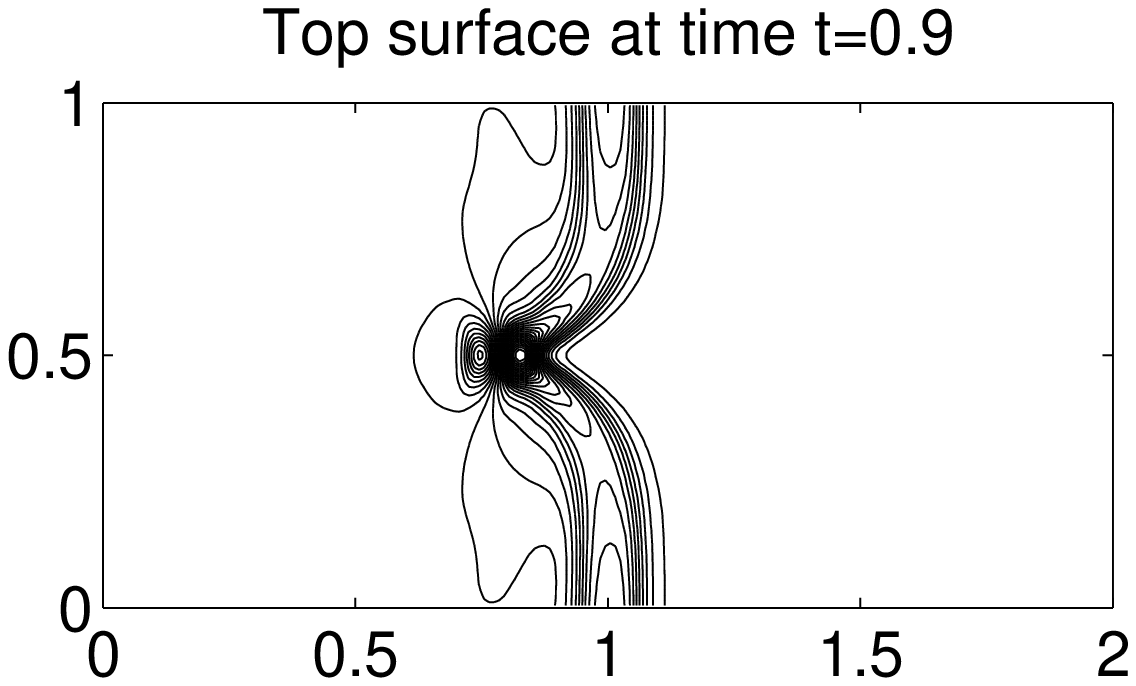, width=0.33\textwidth}
    \epsfig{file=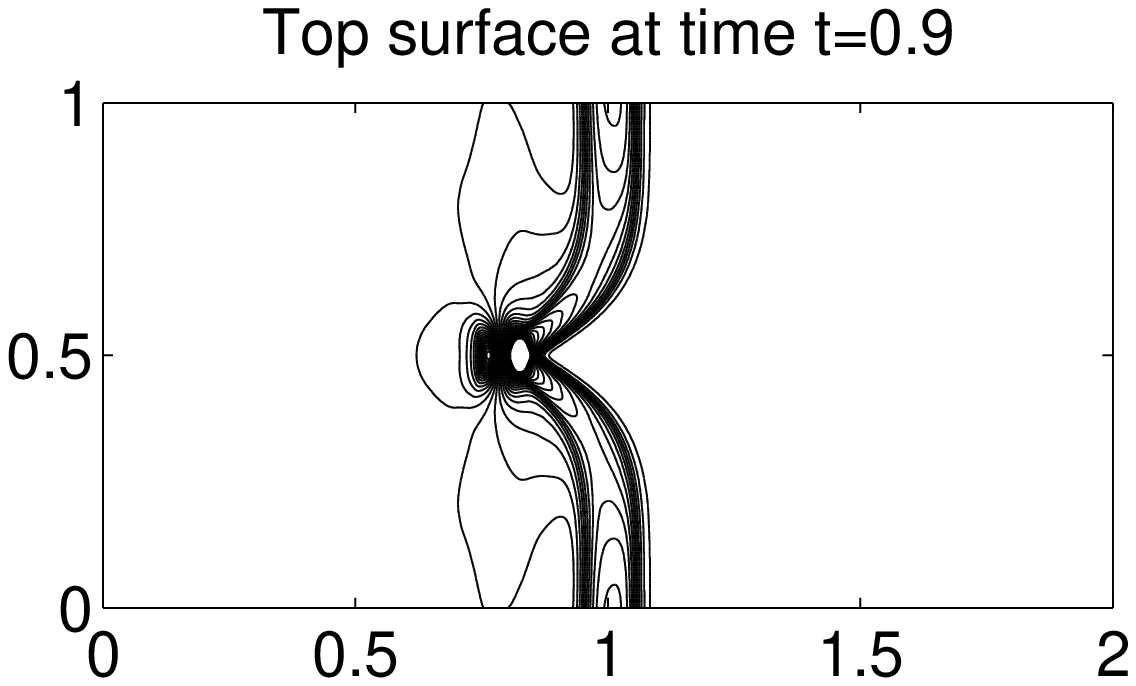,width=0.33\textwidth}\\
    \epsfig{file=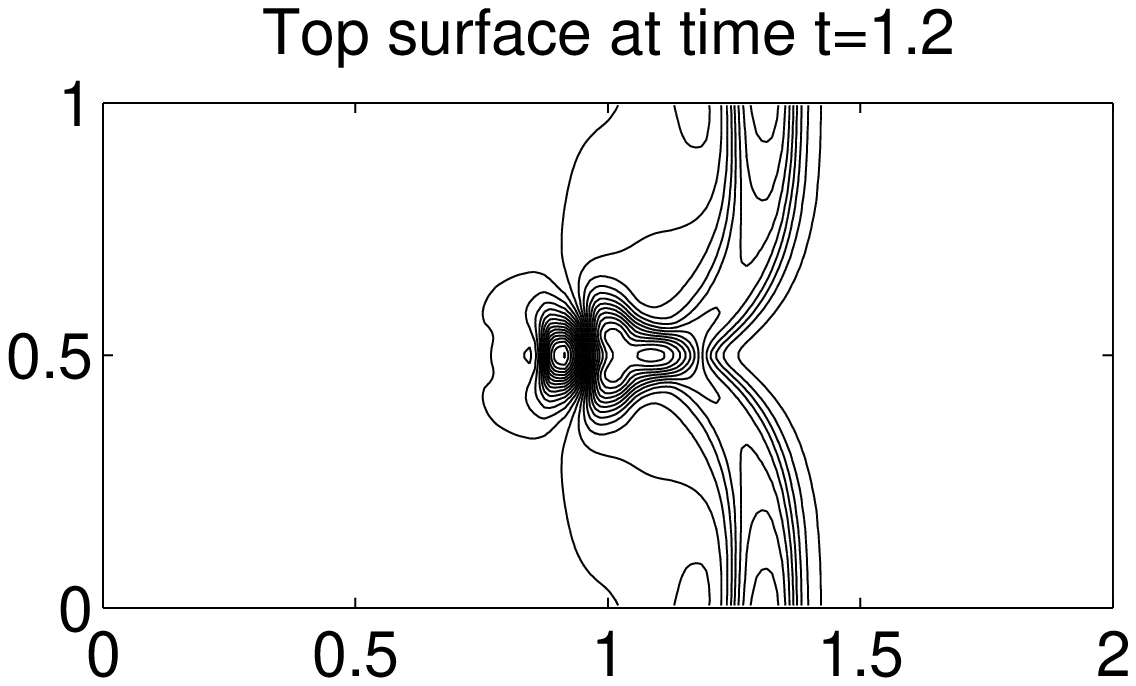, width=0.33\textwidth}
    \epsfig{file=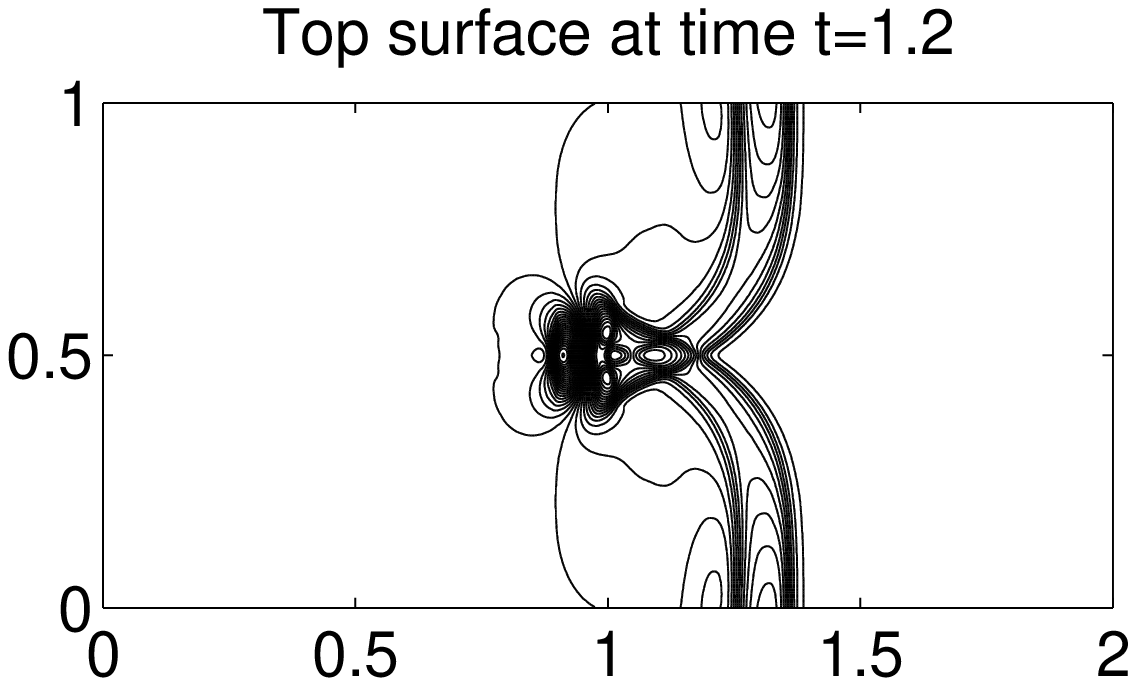,width=0.33\textwidth}\\
    \epsfig{file=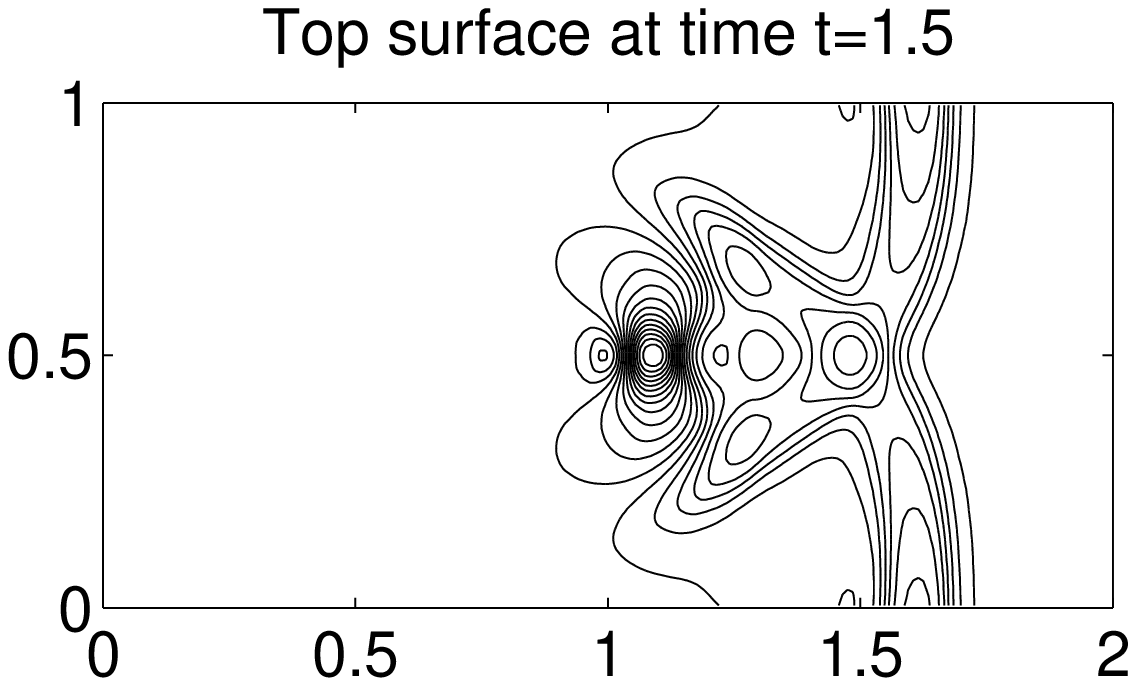, width=0.33\textwidth}
    \epsfig{file=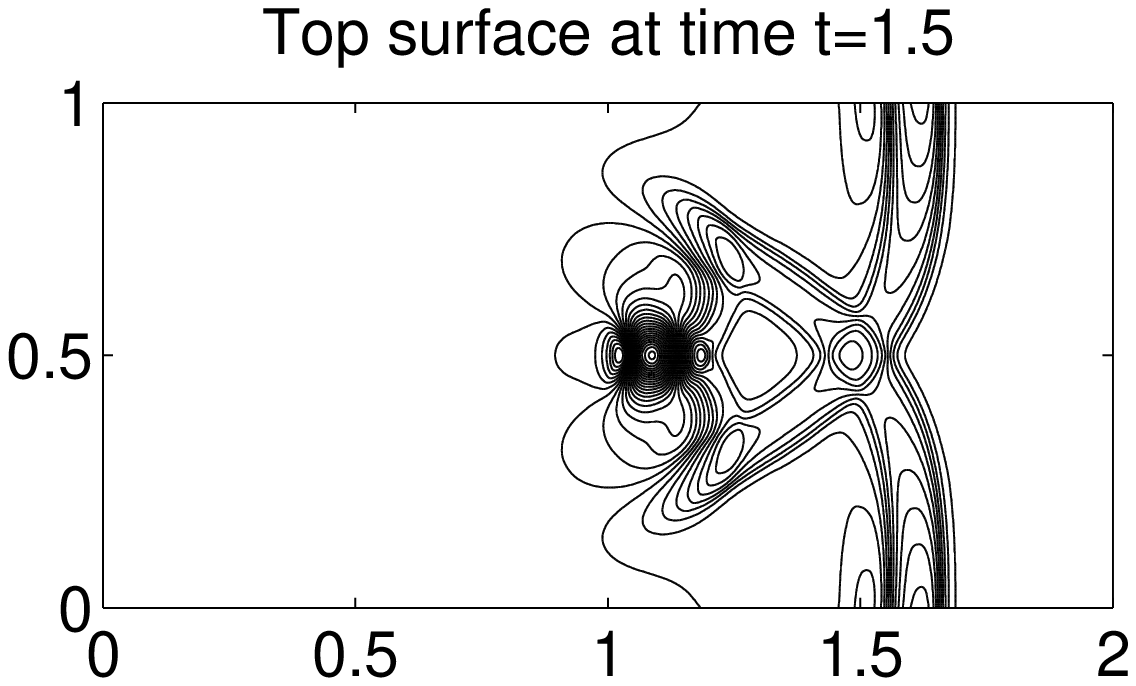,width=0.33\textwidth}\\
    \epsfig{file=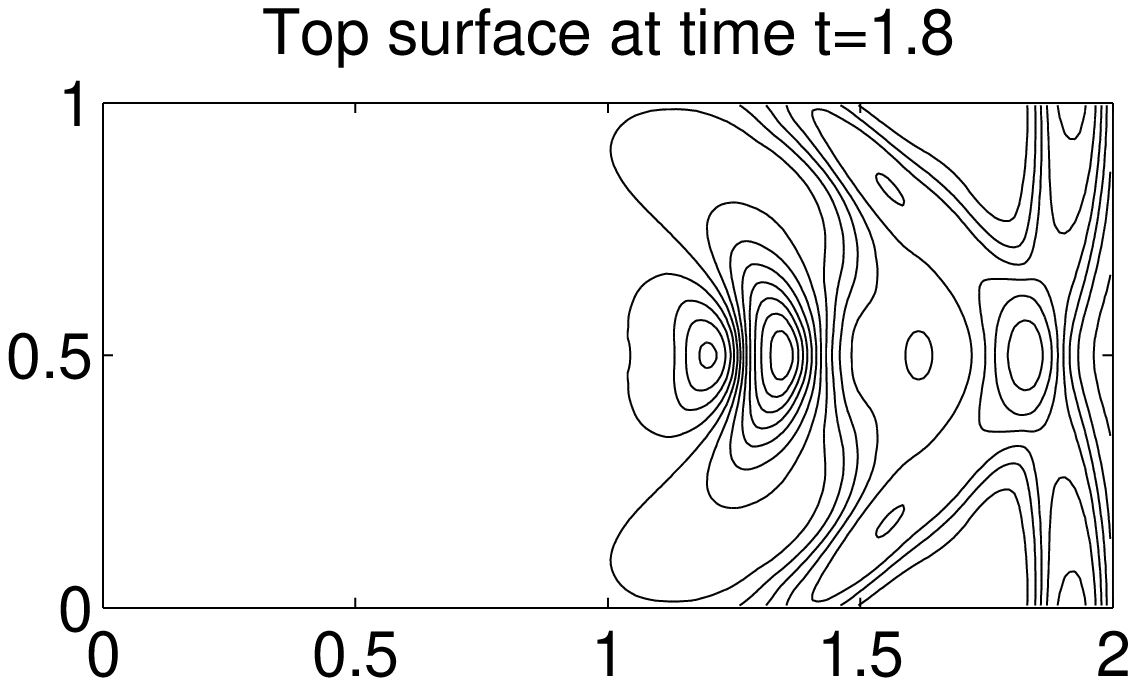, width=0.33\textwidth}
    \epsfig{file=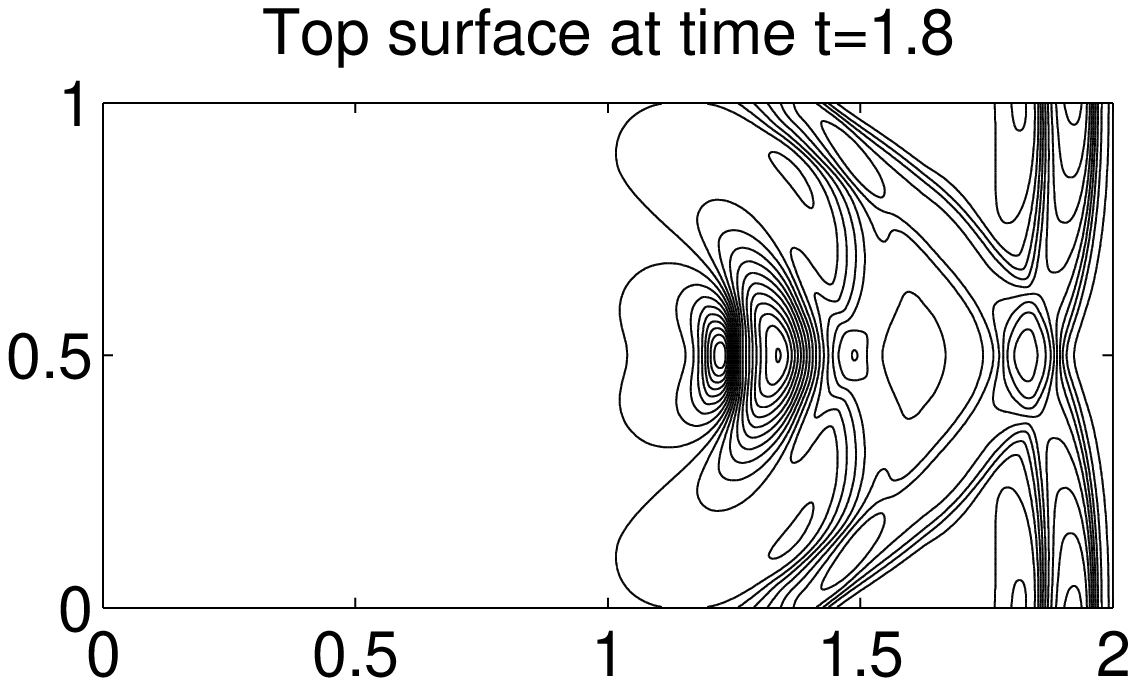,width=0.33\textwidth}\\
    \caption{Two dimensional quasi-stationary problem
    (\ref{QuasiStat_2d_a}), (\ref{QuasiStat_2d_b}).}
    \label{2DQuasiStationary}
  \end{center}
\end{figure}

\subsection{Steady jet in the rotational frame}
This is a classical Rossby adjustment of an unbalanced jet in an open domain, see
e.g.~\cite{bouchut2}. The initial data are a rest state superimposed by a one-dimensioal jet,
$$
h(x,y, 0) = 1.0, \quad u(x,y,0)  = 0, \quad v(x,y,0) = 2 N_L(x),
$$
where  the shape of the velocity $v$ is given by a smooth profile
$$
N_L(x) = \frac{(1+ tanh(4x/L+2)) \,(1 - tanh(4x/L-2))}{(1+tanh(2))^2}
$$
with $L=2$. We have used flat bottom topography $b(x)=0$, the parameter of the
Coriolis forces $f$ and the gravitational acceleration $g$ are set to 1. The
nondimensional parameter representing the effects of Coriolis forces, the Rossby
number $Ro=\frac{| v(x,y,0) |}{f L} = 1$ and the Burgers number reflecting the
nonlinear effects is $Bu = \frac{g |h(x,y,0)|}{f^2 L^2} = 0.25.$ The initial jet
adjusts a momentum unbalance, which emits the waves, the so-called  gravity waves,
propagating out from the jet. The formation of shocks can be noticed within the jet
core approximately at $\pi / f,$ which is a half of a natural time scale $ T_f =
2\pi /f $, see~Figure~5, 6. As time is evolved the solution tends to the equilibrium
state $ f v = gh_x$, which is a geostrophic balance as demonstrated in Figure~7. We
can notice that even for  long time simulations there are still small oscillations
around the geostrophic equilibrium. As pointed out by Bouchut et al.\
\cite{bouchut2}
 some wave modes with the frequencies close to $f$ remain for a longer time in the
core of the jet. Their analysis for a linearized situation shows that they
correspond to the gravity wave modes having almost zero group velocity, and thus
are almost not propagating. For another extensive study of the stability of jets,
which gives interesting eigenfunctions similar to those in Figure~5 we refer
to \cite{gjevik}.

\begin{figure}[ht]
\begin{center}
{\hbox{\includegraphics[scale=0.50]{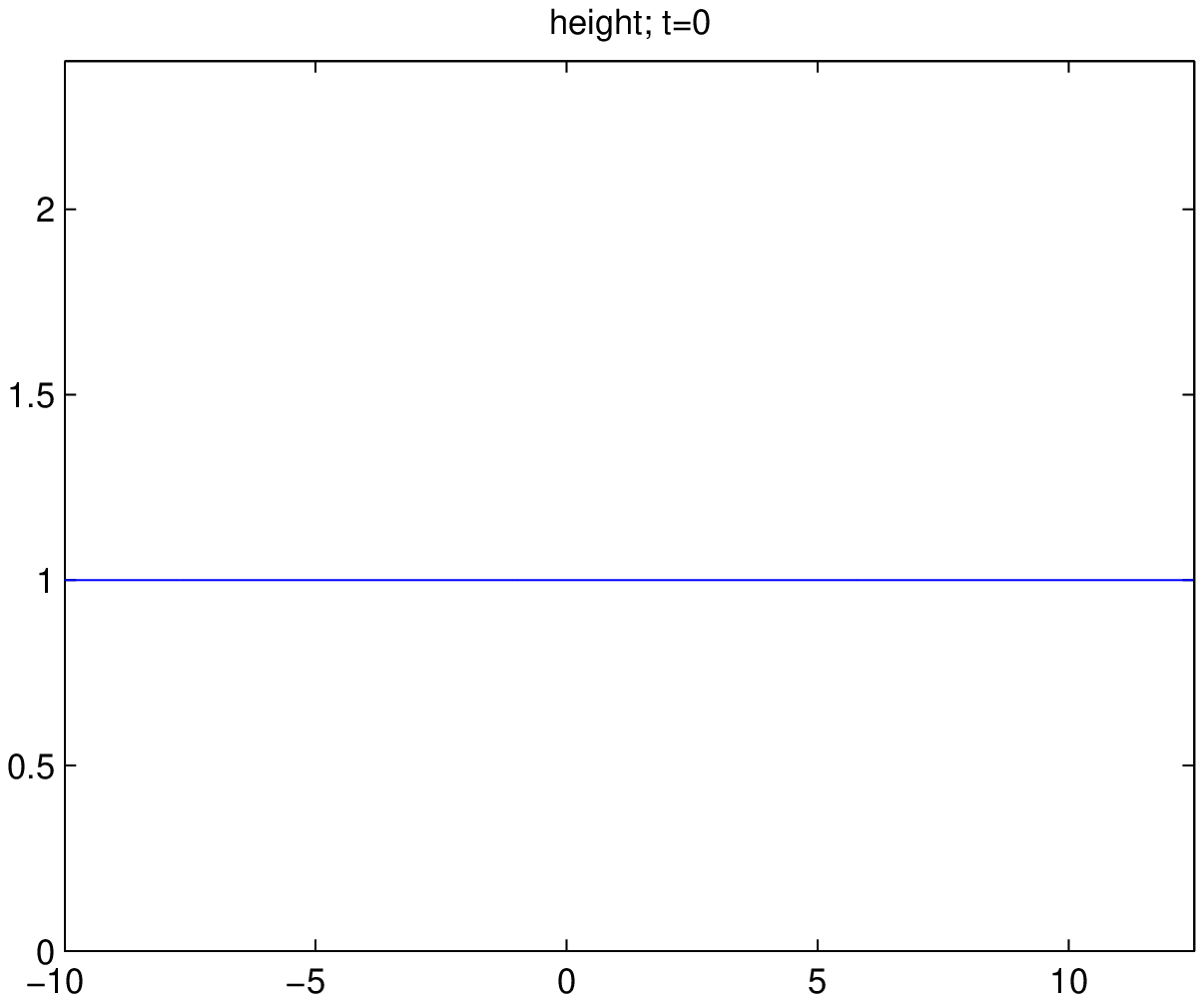}
\includegraphics[scale=0.50]{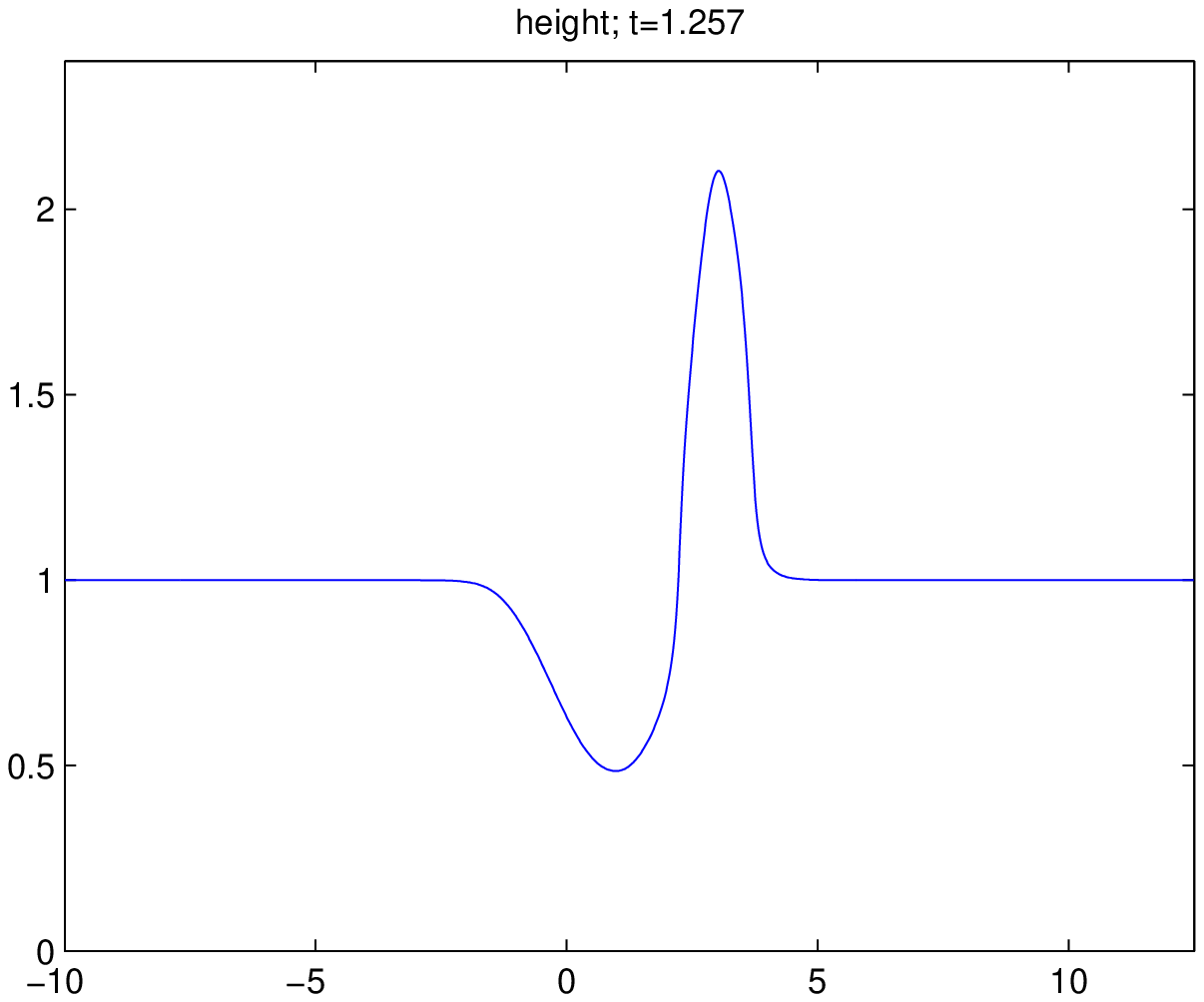}}}
{\hbox{
\includegraphics[scale=0.50]{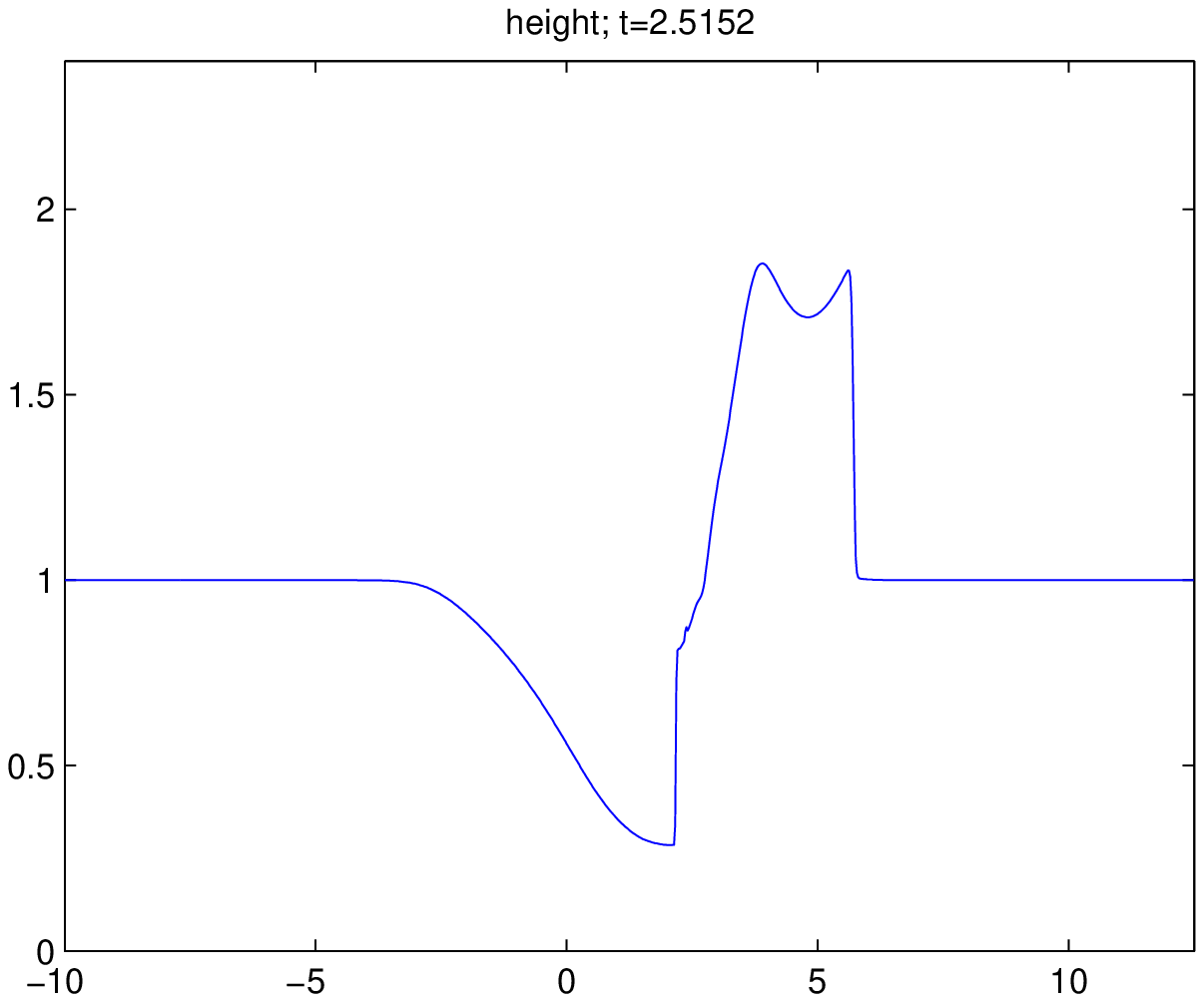}
\includegraphics[scale=0.50]{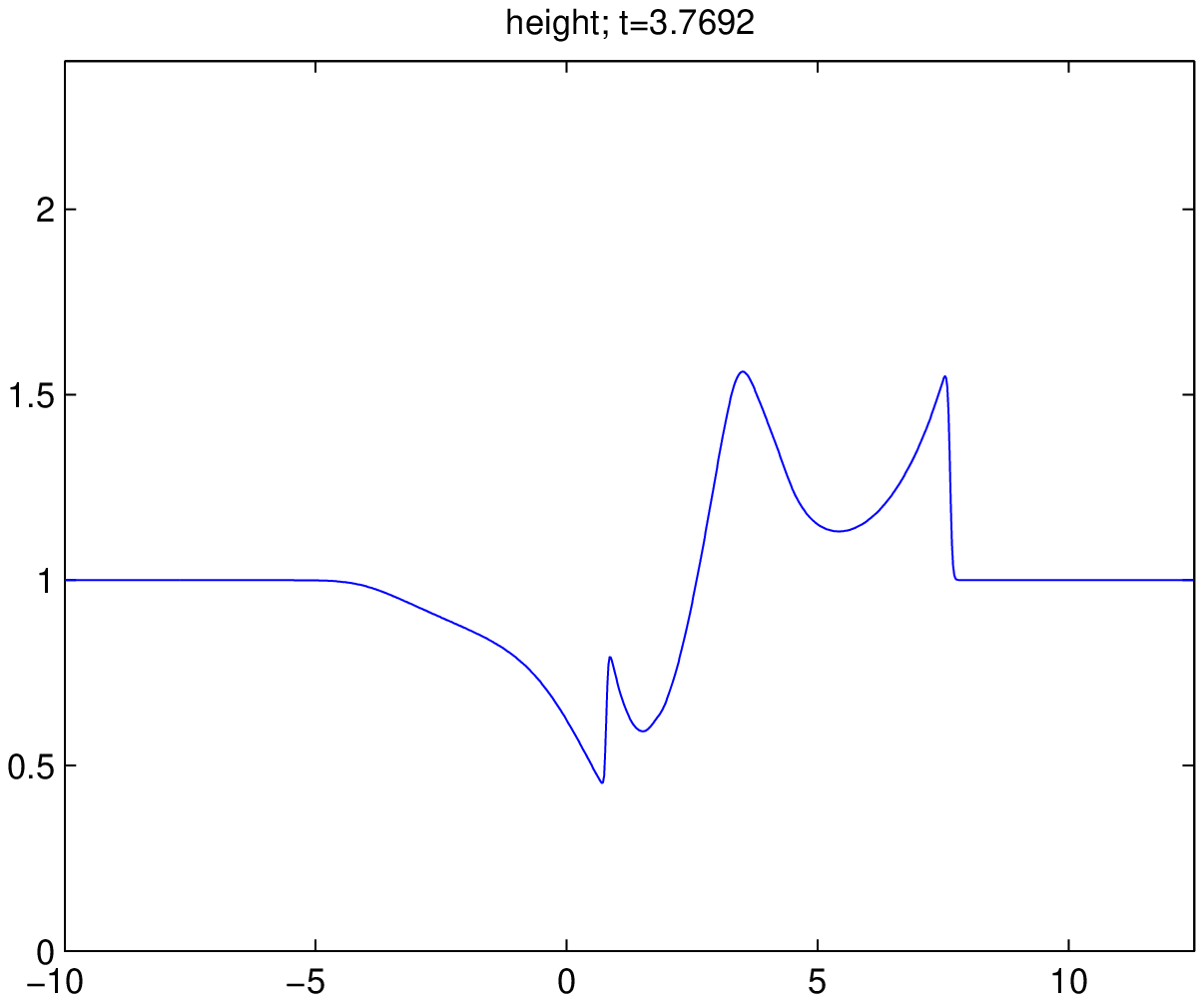}}}
{\hbox{\includegraphics[scale = 0.50]{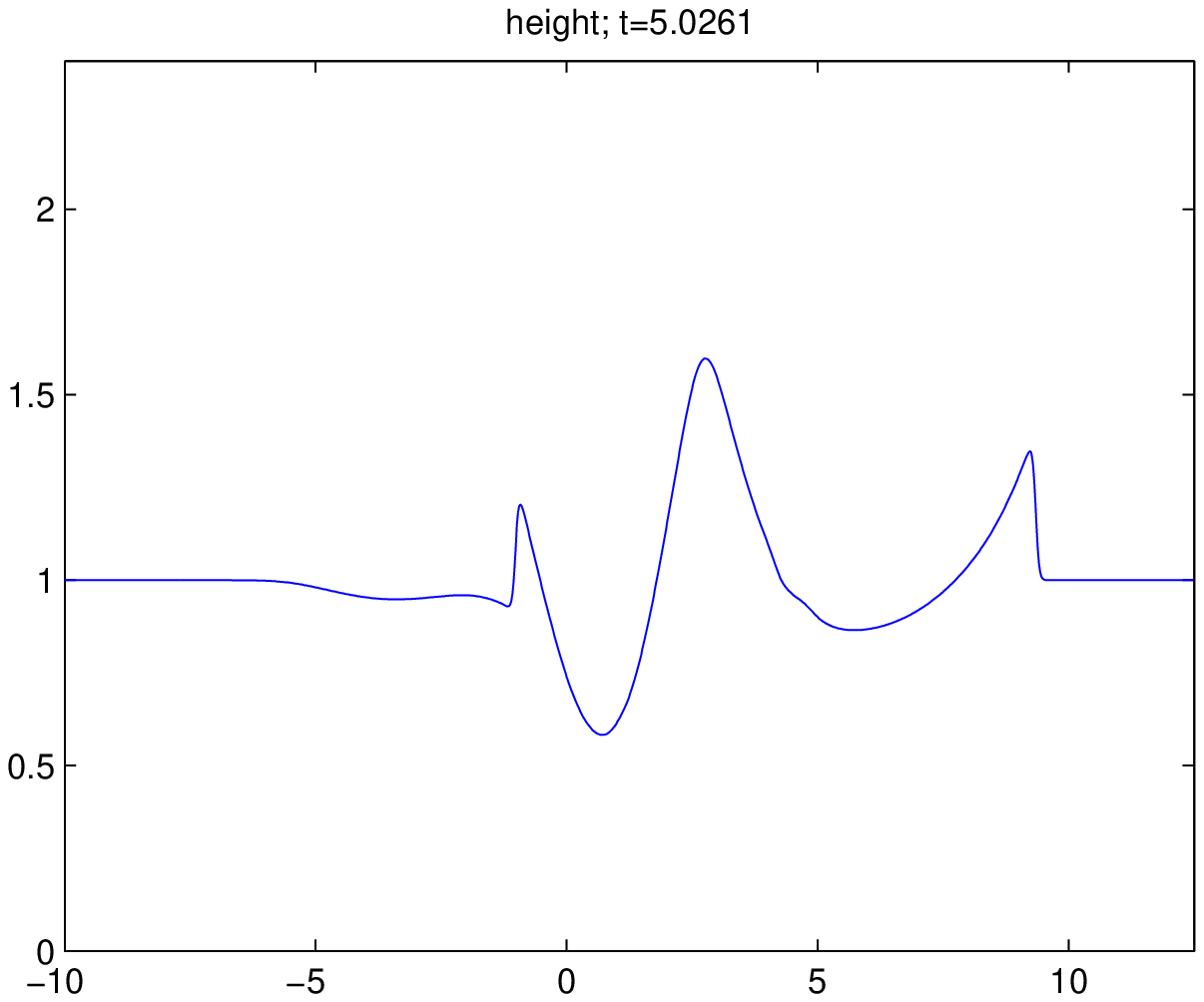}
\includegraphics[scale=0.50]{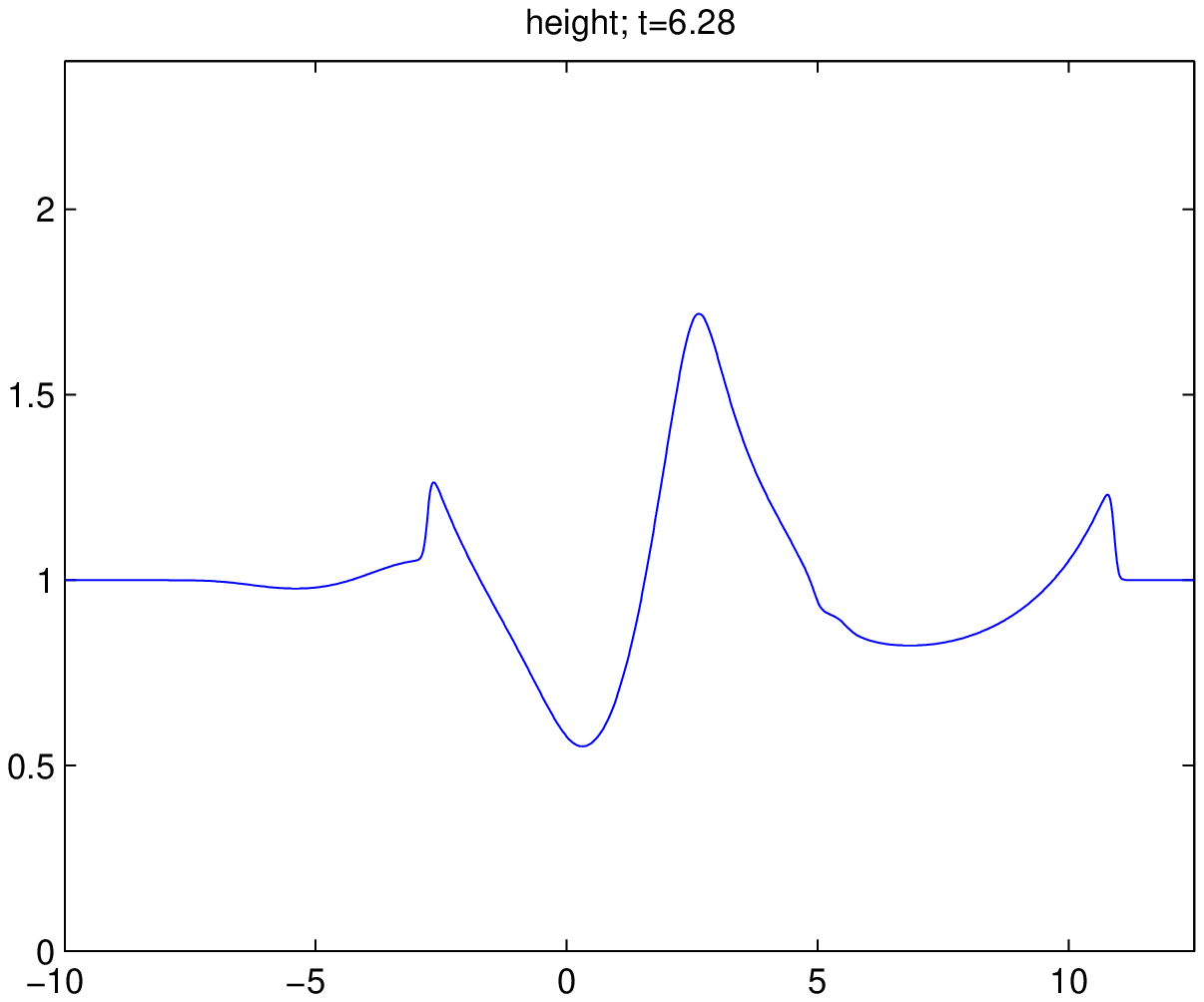}}}
\caption{One-dimensional Rossby adjustment problem, time evolution of water height.}
\end{center}
\end{figure}

\begin{figure}[ht]
\begin{center}
{\hbox{\includegraphics[scale=0.50] {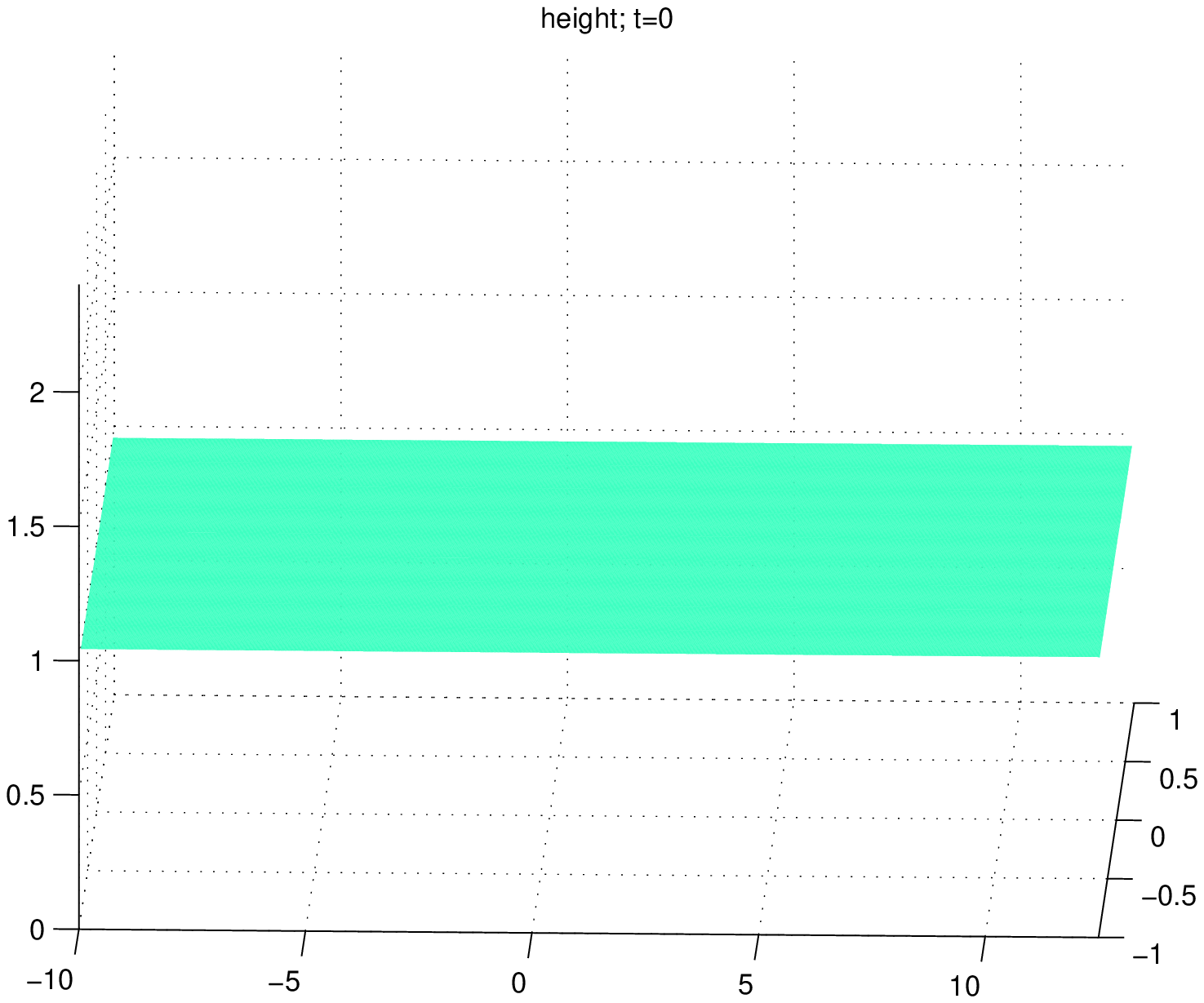}
\includegraphics[scale=0.50]
{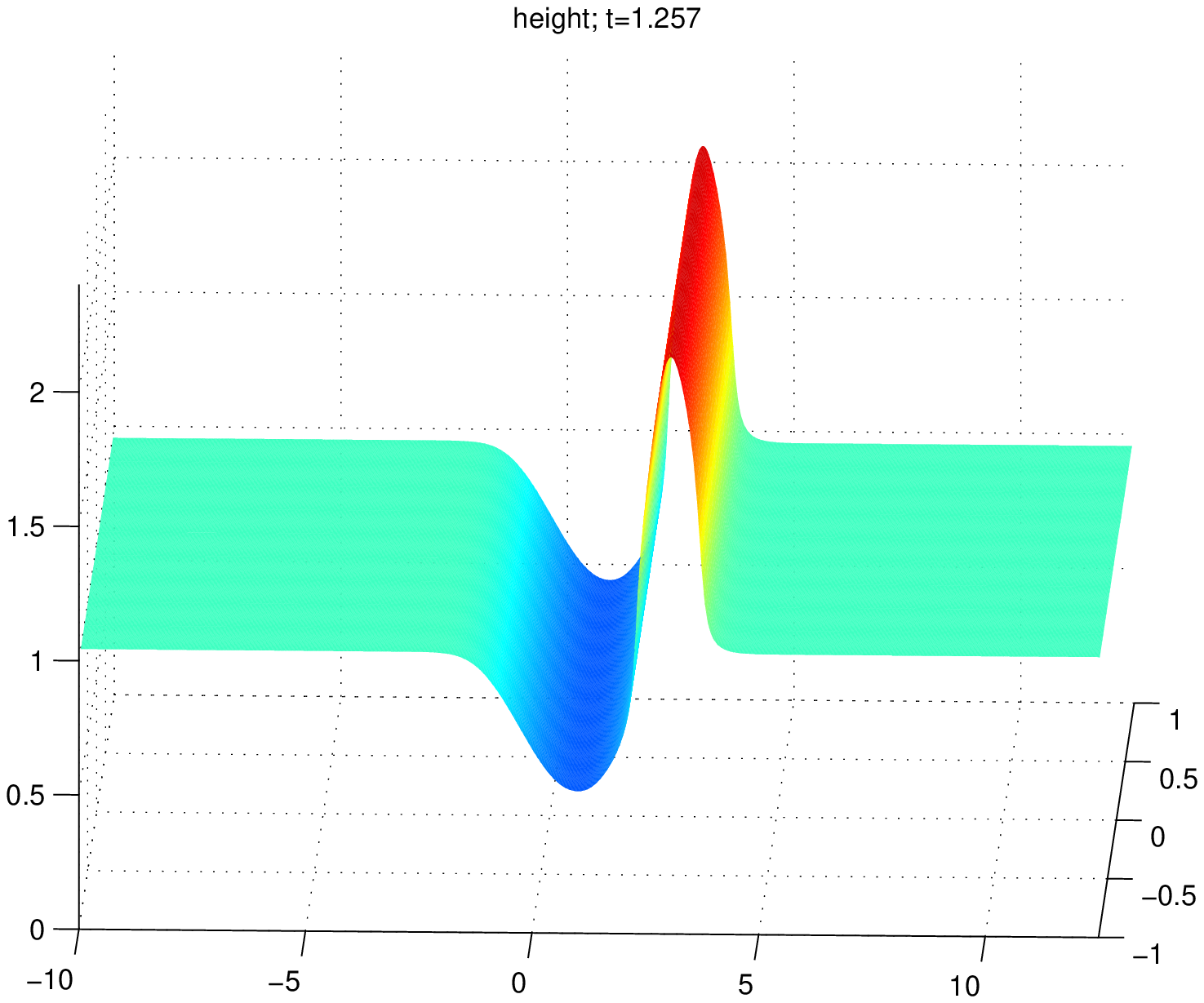}}}

{\hbox{
\includegraphics[scale=0.50]
{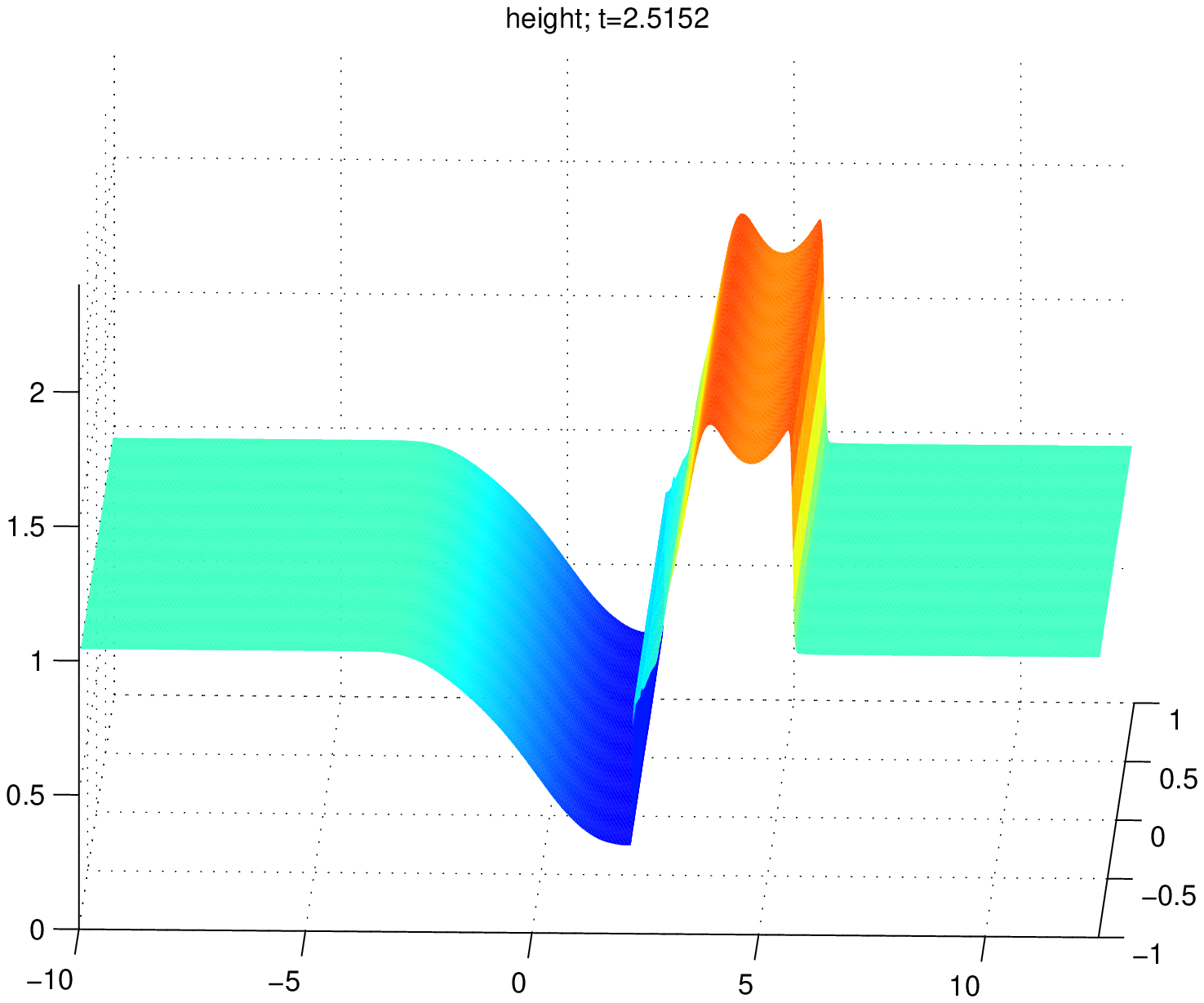}
\includegraphics[scale=0.50]
{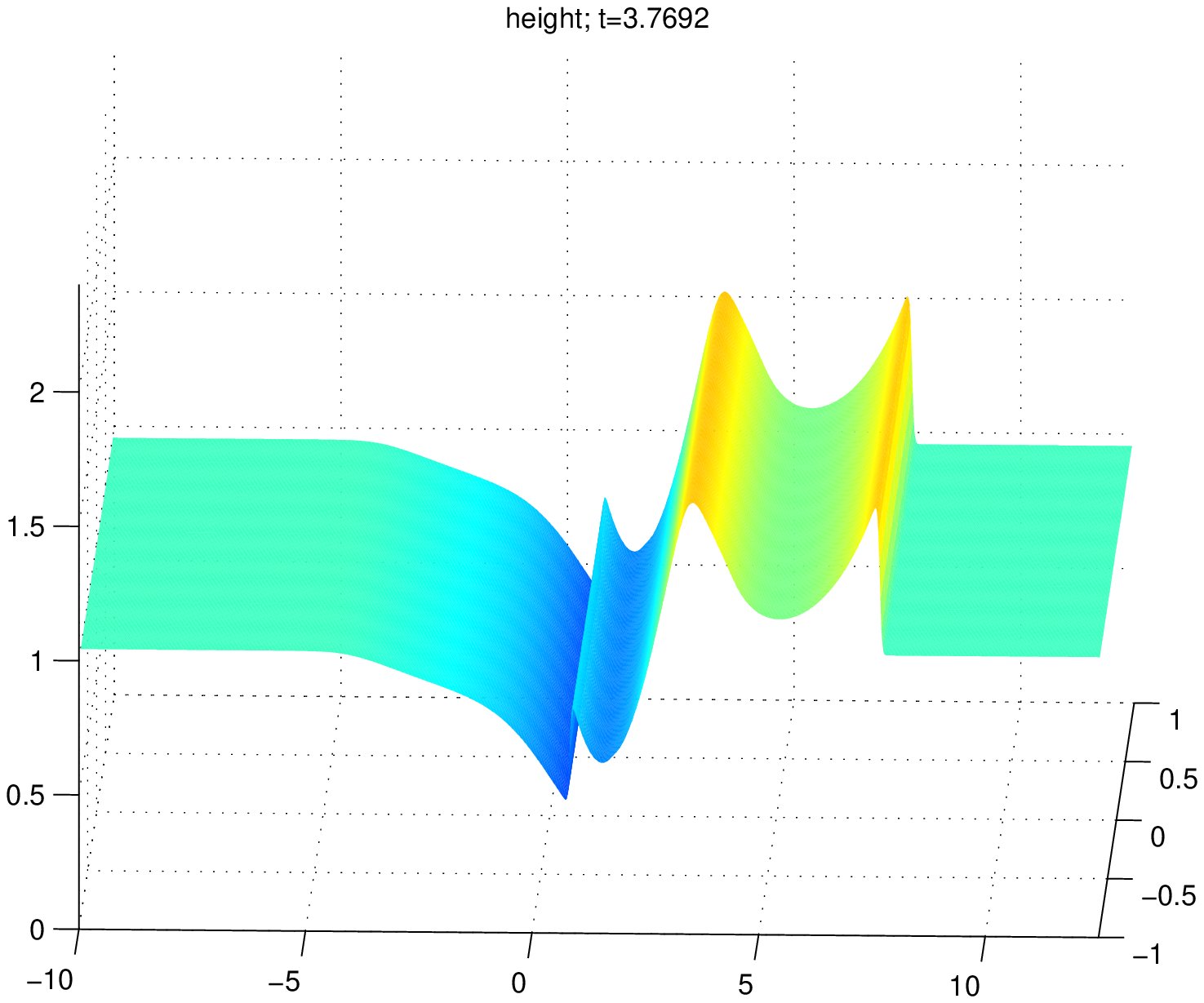}}}

{\hbox{\includegraphics[scale = 0.50] {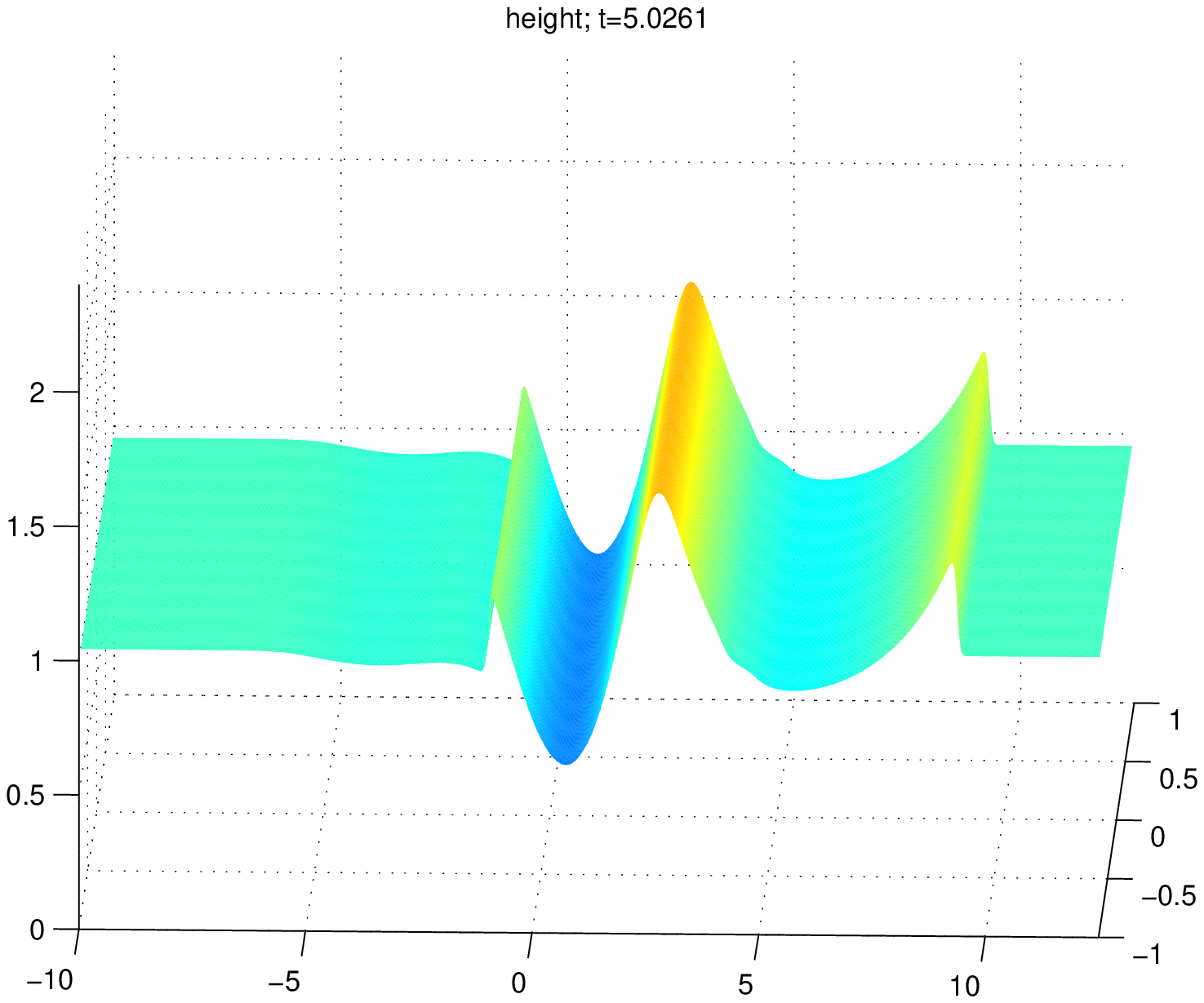}
\includegraphics[scale=0.50]
{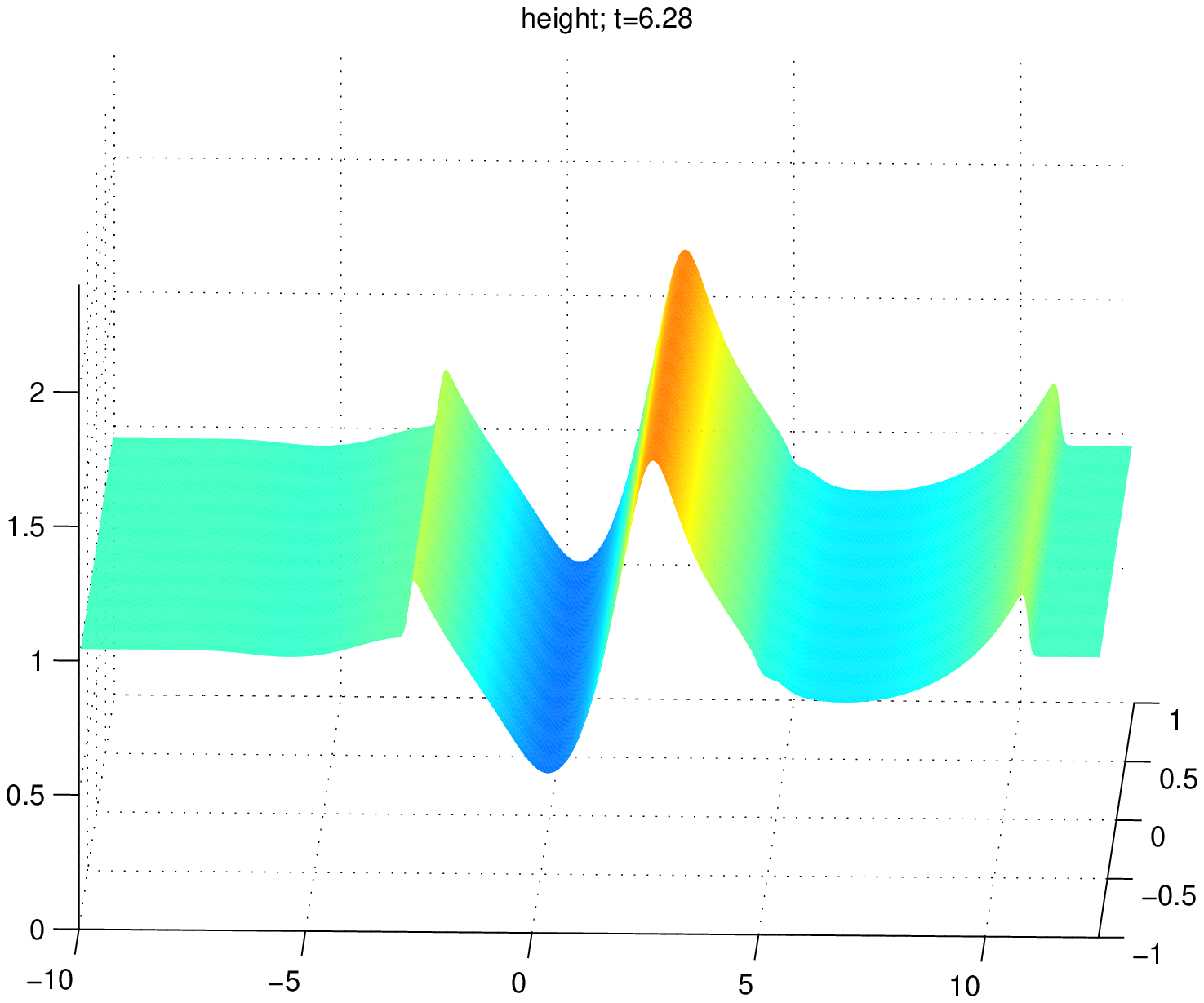}}}

\caption{Rossby adjustment problem, time evolution of water height, two-dimensional graphs.}
\end{center}
\end{figure}

\begin{figure}[hb]
\begin{center}
 \hbox{  \epsfig{file=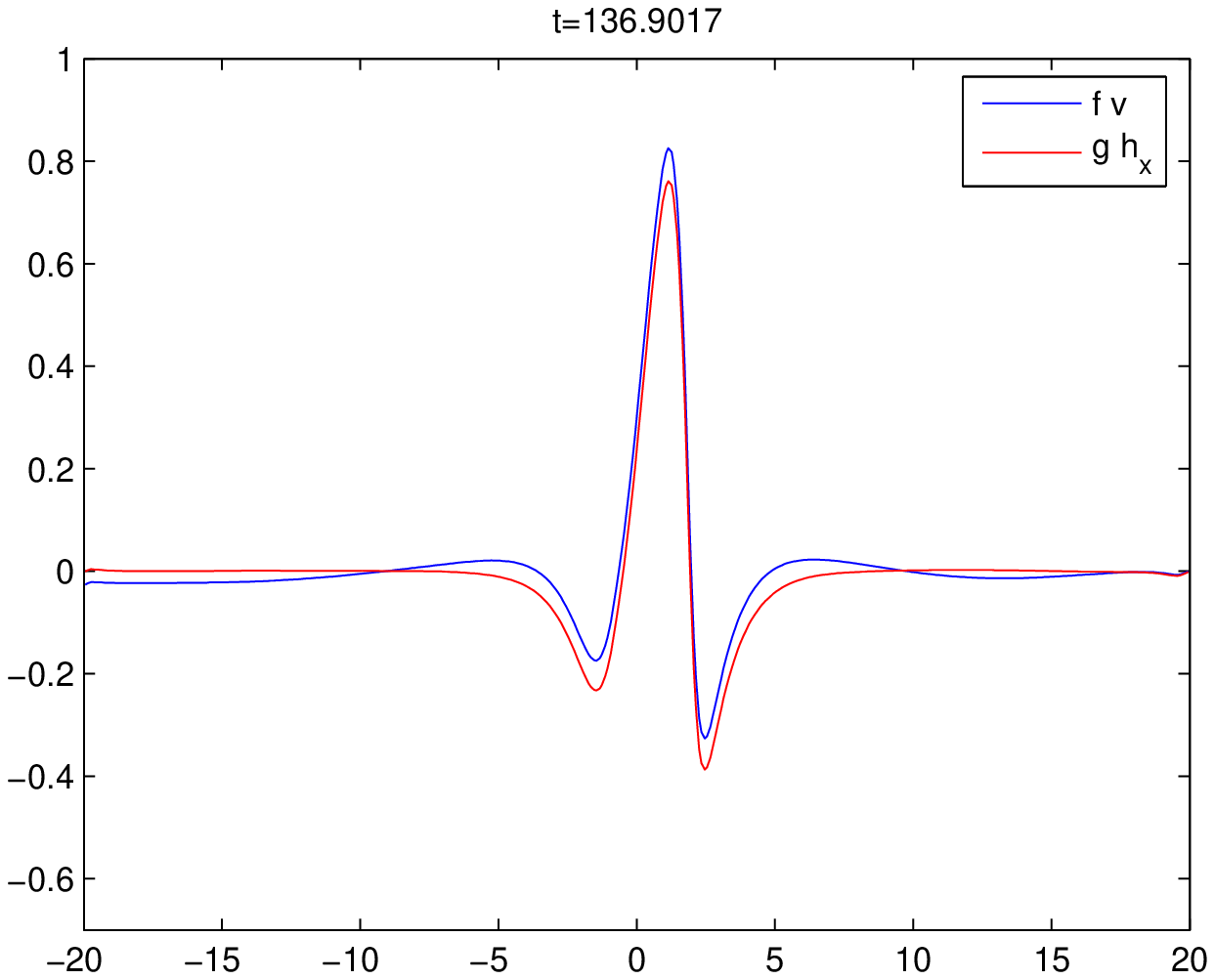, height=6cm}
  \epsfig{file=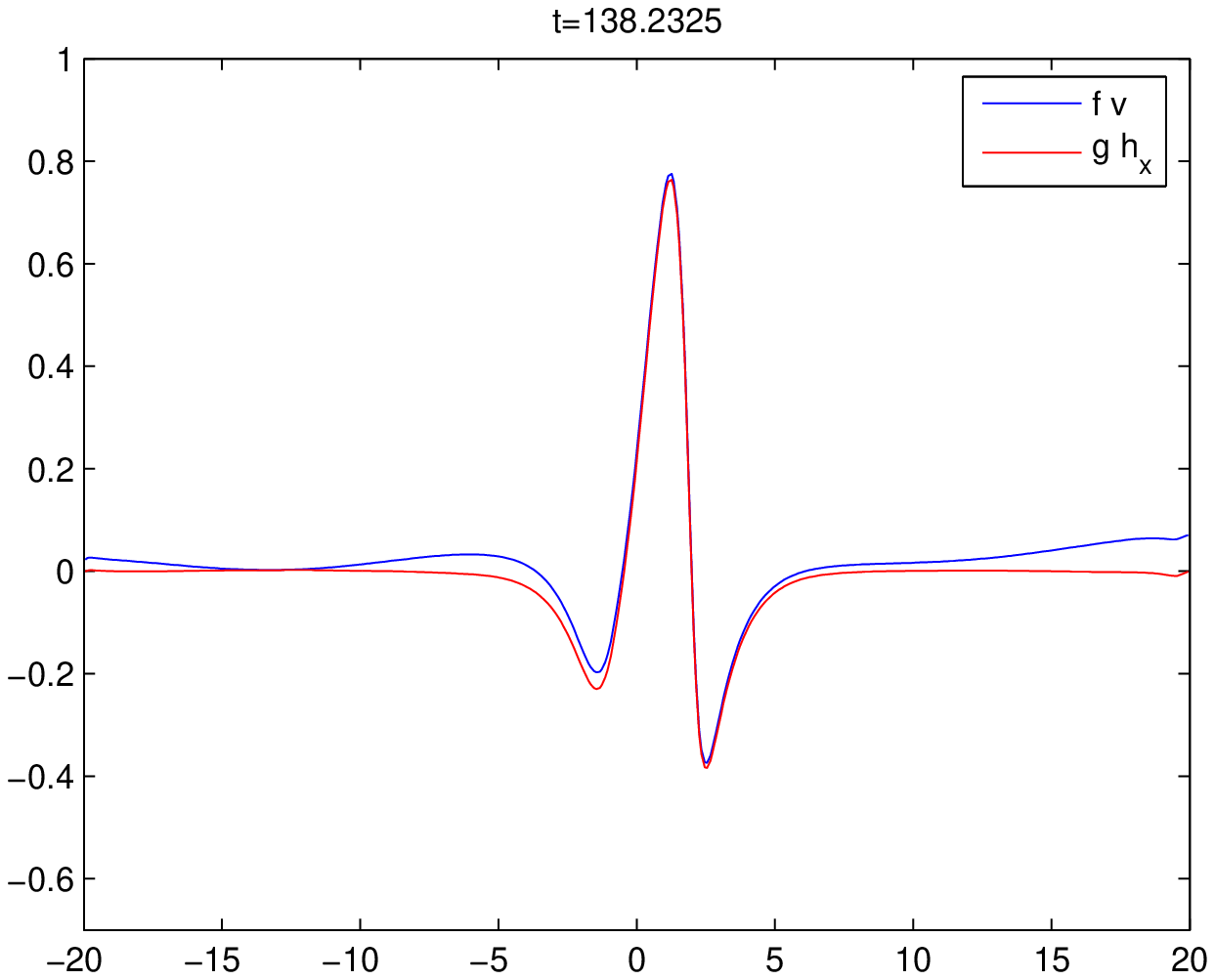, height=6cm}}
\epsfig{file=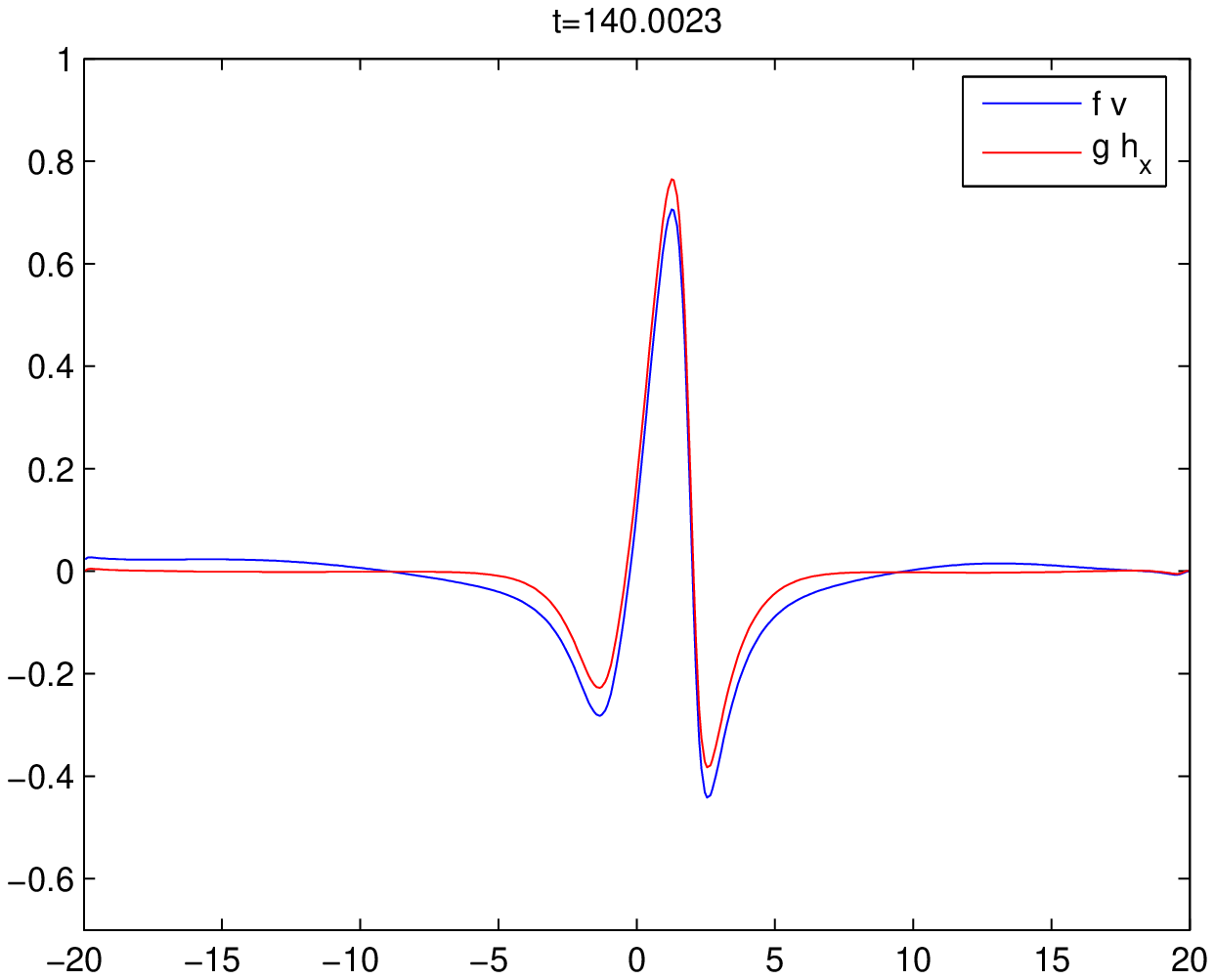, height=6cm} \caption{One-dimensional
Rossby adjustment problem at different times, geostrophic balance.}
\end{center}

\end{figure}

\subsection{Accuracy and performance}
In this experiment we compare  accuracy and computational time of the well-balanced
FVEG, the second order well-balanced FV method of Audusse et al. \cite{bouchut1} and
its fourth order extension due to Noelle et al.\ \cite{noelle1}.
We choose a fully two-dimensional experiment analogous to that of Xing and Shu
\cite{XingShu2005}, but include moreover the effect of Coriolis forces by setting $f=10$.
The gravitational constant
was set to $g=9.812.$ The bottom topography and the initial data are given as
follows
\begin{eqnarray}
\nonumber b(x,y) &=& \sin(2 \pi x) + \cos(2 \pi y), \\
 \nonumber h(x,y,0) &=& 10 + \mbox{exp}({\sin(2 \pi x)}) \cos(2 \pi y),
\\ \nonumber
hu(x,y,0) &=& \sin(\cos(2 \pi x)) \sin(2 \pi y),
\\ \nonumber
h v(x,y,0) &=& \cos(2 \pi x) \cos(\sin(2 \pi y)).
\end{eqnarray}
The computational domain $[0,1] \times [0,1]$  was consecutively divided into $25$,
$50, \dots, 800$ mesh cells in each direction. We have compared solutions obtained
by the second order FVEG scheme as well as by the second order and fourth order
well-balanced FVM at time $T=0.05$. For the second order well-balanced FVM of Audusse
et al.\ the second order Runge-Kutta method was used for time integration, the third
order Gaussian quadrature was used for cell-interface integrals of fluxes and the
second order WENO recovery was applied. The reference solution was obtained by the
fourth order well-balanced FV method of Noelle et al.\ \cite{noelle1}.

Tables~3 and 4 contain the $L^1$ errors and experimental order of convergence (EOC)
for the FVEG for both CFL numbers 0.8 as well as 0.5, respectively. The
well-balanced higher order directional splitting FVM is in general stable up to
CFL=0.5. The $L^1$ errors for its second order version are presented in Table~5. We
can indeed see that both methods are second order accurate in all components.
Moreover,  the second order FVEG scheme is almost 10 times more accurate than the
second order FVM, see~Table~5 as well as the left picture of the Figure~8.

Figure~8 illustrates the CPU/accuracy behaviour graphically. We use the logarithmic
scale on $x-,y-$ axis. On the $y-$ axis the $L^1$ errors in first component $h$ is
depicted. Errors in other components yield analogous results. On the left of
Figure~8 a comparison between second order FVEG and FV methods are presented,
whereas on the right we show the comparison between the fourth order well-balanced
FVM of Noelle \cite{noelle1} and  the second order FVEG scheme. The FVEG schemes
yields on coarse meshes still more accurate solutions. In fact, for meshes up to
approximately $100 \times 100$ cells, which are actually often used for practical
computations, it is more efficient to use the second order FVEG scheme than the
fourth order FVM. The superiority of the fourth order scheme is showing up on fine
grids, see~the right graph of Figure~8.

We should point out that no attempt has been made in order to optimize the codes
with respect to their CPU performance. Our extensive numerical treatment indicates
that both well-balanced second order methods, the FVEG as well as the FVM are
actually comparable with respect to their computational time.

\begin{table}[h]
\caption{FVEG scheme: Convergence in the $L^1$ norm, CFL=0.8}
\begin{center}
\begin{tabular}[h]{|c|c|c|c|c|c|c|} \hline
N&  $L^1$ error in $ h $ &EOC& $L^1$ error in $ hu $& EOC& $L^1$ error in $ hv $&
EOC
\\ \hline \hline
25 & 1.04e-02&        &3.56e-02&    &8.52e-02&
\\ \hline
50 & 2.42e-03&  2.10   &8.71e-03&2.03&2.15e-02&1.99
\\ \hline
100& 6.01e-04& 2.01 &2.23e-03 &1.96&5.50e-03&1.96
\\ \hline
200& 1.54e-04& 1.96 &5.76e-04 &1.95&1.44e-03&1.93
\\ \hline
400& 3.97e-05 &1.96 &1.47e-04& 1.97&3.69e-04&1.96
\\ \hline
800& 1.02e-05&  1.97&3.71e-05& 1.98&9.40e-05&1.97
\\ \hline \hline
\end{tabular}
\end{center}
\end{table}

\begin{table}[h]
\caption{FVEG scheme: Convergence in the $L^1$ norm, CFL=0.5}
\begin{center}
\begin{tabular}[h]{|c|c|c|c|c|c|c|} \hline
N&  $L^1$ error in $ h $ &EOC& $L^1$ error in $ hu $& EOC& $L^1$ error in $ hv $&
EOC
\\ \hline \hline
25 & 1.37e-02&        &6.19e-02&    &1.18e-01&
\\ \hline
50 & 2.80e-03& 2.29   &1.05e-02&2.56&2.33e-02&2.34
\\ \hline
100& 5.23e-04& 2.42 & 1.80e-03 &2.54&4.25e-03&2.45
\\ \hline
200& 1.04e-04& 2.33 &3.63e-04 &2.31&8.12e-04&2.39
\\ \hline
400& 2.45e-05 & 2.09 &8.79e-05&2.05&1.80e-04&2.17
\\ \hline
800& 6.14e-06& 1.99&  2.20e-05& 2.00&4.36e-05&2.04
\\ \hline \hline
\end{tabular}
\end{center}
\end{table}

\begin{table}[h!]
\caption{FV scheme: Convergence in the $L^1$ norm, CFL=0.5}
\begin{center}
\begin{tabular}[h]{|c|c|c|c|c|c|c|} \hline
N&  $L^1$ error in $ h $ &EOC& $L^1$ error in $ hu $& EOC& $L^1$ error in $ hv $&
EOC
\\ \hline \hline
25 & 4.53e-02&        &2.13e-01&    &3.40e-01&
\\ \hline
50 & 1.32e-02&  1.77   &5.57e-02&1.94&9.51e-02&1.84
\\ \hline
100& 3.50e-03& 1.92 &1.42e-02&1.97&2.52e-02&1.92
\\ \hline
200& 8.95e-04& 1.97 &3.58e-03 &1.99&6.46e-03&1.96
\\ \hline
400& 2.26e-04 &1.99 &8.96e-04& 2.00&1.63e-03&1.99
\\ \hline
800& 5.67e-05&  1.99&2.24e-04& 2.00&4.10e-04&1.99
\\ \hline \hline
\end{tabular}
\end{center}
\end{table}

\begin{figure}[hbt] \label{fig-7}
\begin{center}
 \hbox{\epsfig{file=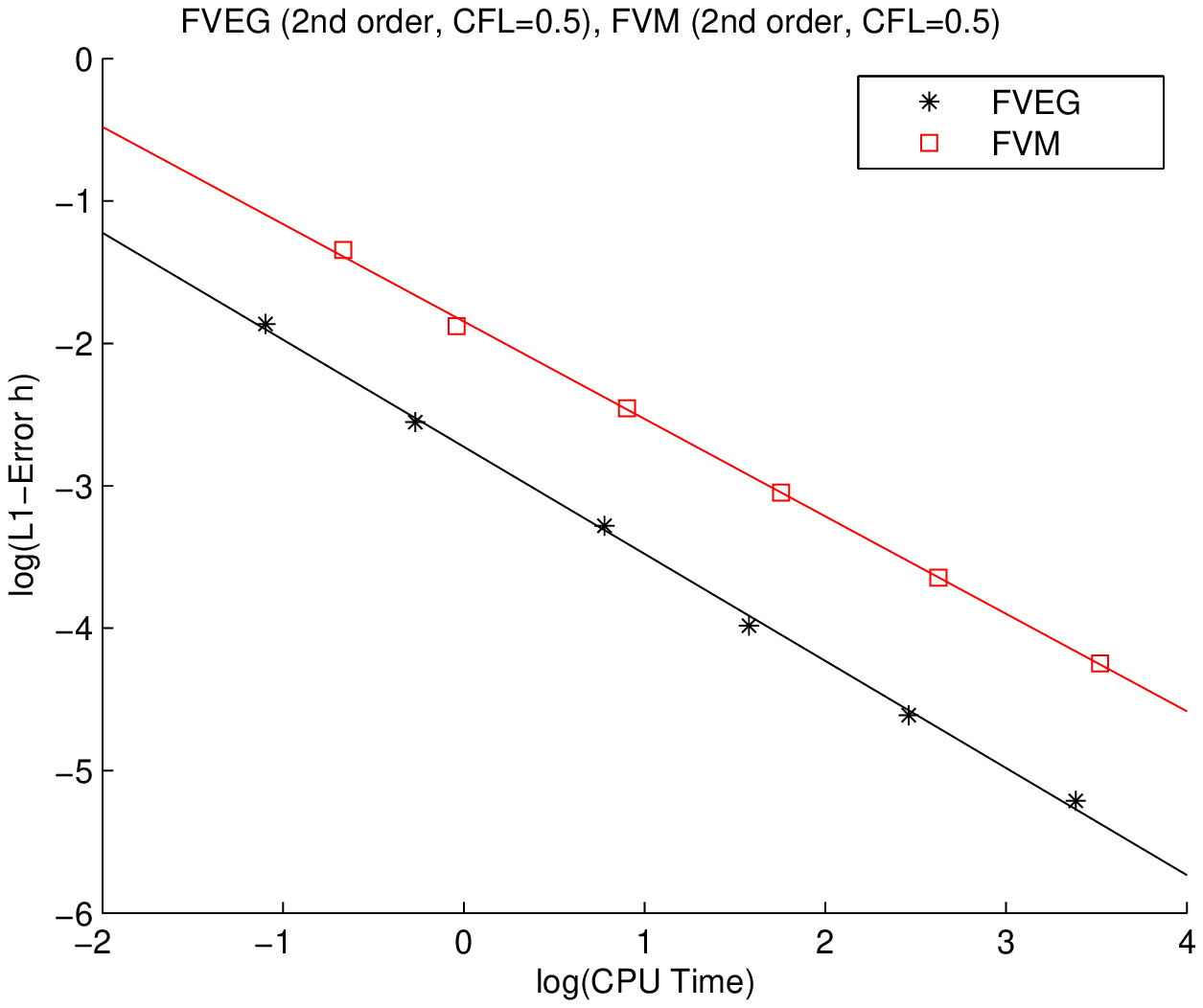, height=6cm}
  \epsfig{file=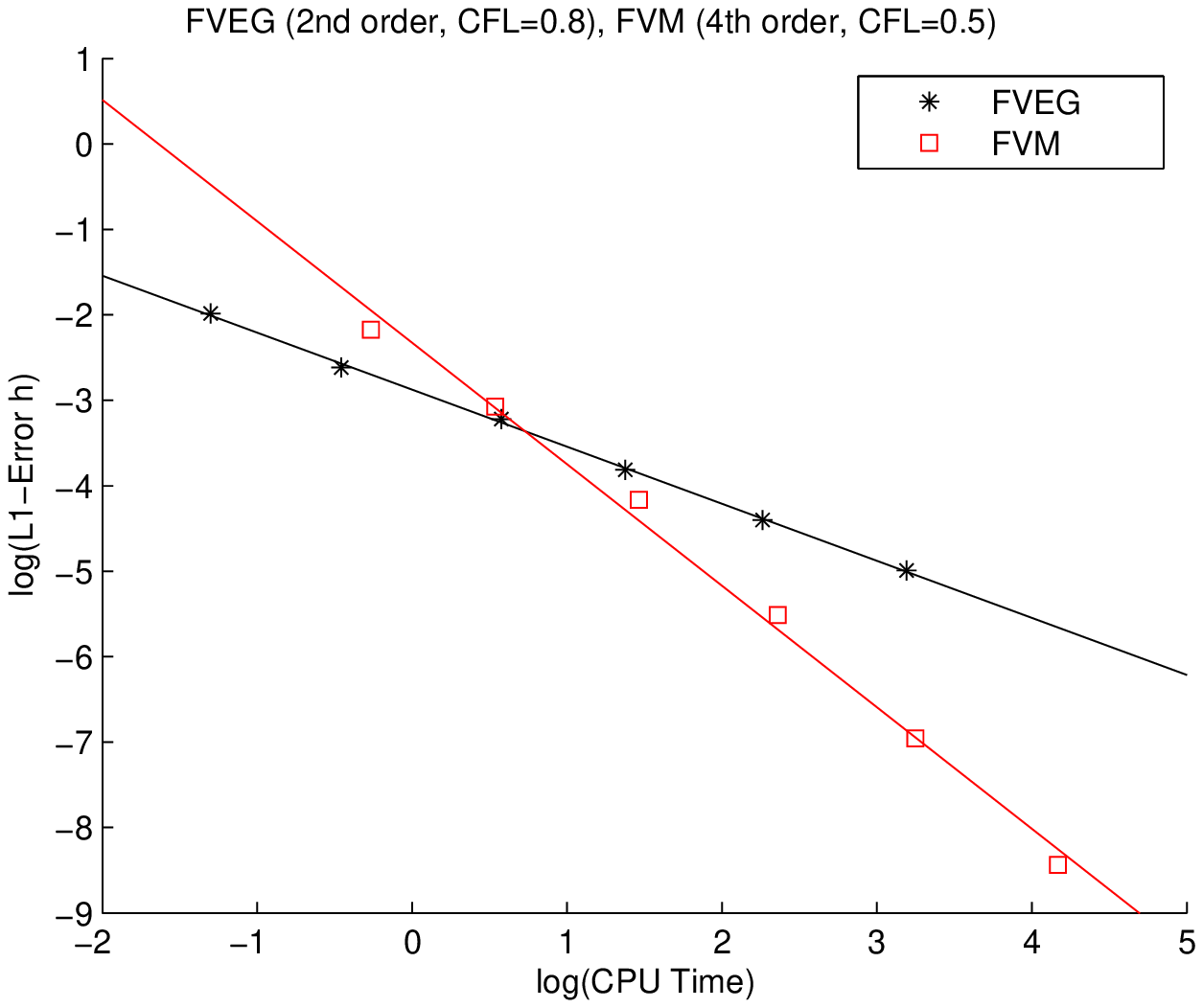, height=6cm}}
\end{center}
\caption{Efficiency test: $L^1$ error over the CPU time;  the second order FVEG and
second order FV schemes (left) as well as the fourth order FVM  and the second order
FVEG scheme (right). }
\end{figure}

\section{Conclusions}

In the present paper we have developed a new well-balanced scheme within the framework of
the finite volume evolution Galerkin (FVEG) scheme. The scheme is applied for the shallow
water equations with source terms modelling the bottom topography and the Coriolis forces.
The key ingredient of this FVEG scheme is a new well-balanced approximate representation of
the solution, cf.\ Lemma~3.2, which together with a recent quadrature rule from \cite{sisc}
leads to the multidimensional approximate evolution operators (\ref{approx_pwc}),
(\ref{approx_pwb}). These approximate evolution operators are used in a predictor step. In
fact we are predicting the solution at cell interfaces and do not need to use the
hydrostatic reconstruction as it is done by Audusse at el.\ \cite{bouchut1}.

In the following correction step, which is the finite volume update step, the source term is
approximated in the interface-based way. We have proved that the lake at rest, the steady
jet in the rotational frame as well as their combinations are preserved, cf.~Theorems~2.1
and 3.1. Numerical experiments in one and two space dimensions demonstrate correct
resolution of these equilibrium states and of their small perturbations. For smooth
solutions the accuracy of the well-balanced FVEG scheme is superior to that of a recent FV
scheme while the CPU time is comparable.

In future work we want to extend our well-balanced schemes to shallow
water equations with nonlinear friction, which appears in oceanographic as well
as river flow modelling. Another challenge is presented by multi-layere shallow
water models, which are important in oceanology, meteorology and climatology.

\bigskip
{\bf Acknowledgements}
\newline
This research was supported by the Graduate Colleges `Conservation principles in the
modelling and simulation of marine, atmospherical and technical systems' and
`Hierarchie and Symmetry in Mathematical Systems' both sponsored by the Deutsche
Forschungsgemeinschaft, and by the EU-Network Hyke HPRN-CT-2002-00282. The authors
gratefully acknowledge these supports. Furthermore, we would like to thank Normann
Pankratz, RWTH Aachen, for fruitful discussions and for providing the FV results
\cite{bouchut1,noelle1} used in the comparisons in Section 5.4.

\appendix

\section{Derivation of the exact integral equations}
\label{section:appendix_EG}

Applying theory of bicharacteristics to (\ref{sw_linearized}) we can derive exact integral
equations in an analogous way as in \cite{sisc}. In order to keep the presentation
self-contained  we briefly rewrite  main steps of the derivation.

\noindent Let $\umat{R}$ be the matrix of right eigenvectors corresponding to
direction $\uu{n}:=(\cos(\theta),\sin(\theta))$. Its inverse  reads
$$
\umat{R}^{-1} = \frac{1}{2} \left( \begin{array}{ccc}
  -1 &  \frac{\tilde{c}}{g}\cos\theta &      \frac{\tilde{c}}{g}\sin\theta \\
   0 & 2\sin\theta & -2\cos\theta \\
   1 &  \frac{\tilde{c}}{g}\cos\theta &  \frac{\tilde{c}}{g}\sin\theta \\
                                \end{array} \right ).
$$
Let us define the vector of characteristic variables $\uu{v}$ by
$$
 \uu{v} :=  \umat{R}^{-1} \uu{w}.
$$
Multiplying  system (\ref{sw_linearized}) by $ \umat{R}^{-1} $ from the left yields the
following characteristic system
$$
\frac{\partial \uu{v}}{\partial t}
 + \umat{B}_1 (\tilde{\uu{v}})
\frac{\partial \uu{v}}{\partial x} + \umat{B}_2 (\tilde{\uu{v}}) \frac{\partial
\uu{v}}{\partial y} = \uu{r},
$$
where
\begin{eqnarray}
\umat{B}_1  &=&  \displaystyle \left( \begin{array}{ccc}
     \tilde{u} - \tilde{c}\cos\theta &  -\frac{1}{2}\,\tilde{h}\sin\theta & 0 \\
    - g\sin\theta & \tilde{u} & g\sin\theta \\
     0 &  \frac{1}{2}\,\tilde{h}\sin\theta & \tilde{u} + \tilde{c}\cos\theta \\
                                \end{array} \right ),
\nonumber\\
\umat{B}_2  &=&  \displaystyle \left( \begin{array}{ccc}
\tilde{v} - \tilde{c}\sin\theta &  \frac{1}{2}\,\tilde{h}\cos\theta & 0 \\
g \cos\theta & \tilde{v} & - g\cos\theta \\
0 &  -\frac{1}{2}\,\tilde{h}\cos\theta & \tilde{v} + \tilde{c}\sin\theta \\
                                \end{array} \right ),
\nonumber \\
\uu{r}(\uu{n}) &=& \trivek{r_1}{%
r_2}{%
r_3}= \umat{R}^{-1}(\uu{n}) \uu{t} = \trivek{\frac{1}{2}\frac{\tilde{c}}{g}((-g b_x
+ f v) \cos \theta
 - (g b_y + f u) \sin\theta)}%
 { (-g b_x + f v) \sin\theta + (g b_y + f u) \cos \theta }%
 { \frac{1}{2}  \frac{\tilde{c}}{g}( (-g b_x + f v) \cos \theta
 - (g b_y  + f u) \sin\theta)}
\nonumber \\ \phantom{mm}&&
\end{eqnarray}
and the characteristic variables $\uu{v}$ are
\begin{equation}
\nonumber
\uu{v}(\uu{n}) = \trivek{v_1}{%
 v_2}{%
 v_3} = \umat{R}^{-1}(\uu{n}) \uu{u}
= \trivek{\frac{1}{2}(- h +  \frac{\tilde{c}}{g} u \cos\theta
                          +  \frac{\tilde{c}}{g} v \sin\theta)}{%
                        u \sin\theta - v \cos\theta}{%
          \frac{1}{2}(  h + \frac{\tilde{c}}{g} u \cos\theta
                          + \frac{\tilde{c}}{g} v \sin\theta)}.
\end{equation}
The quasi-diagonalised system of the linearized shallow water equations has the
following form
\begin{equation}
\label{so_3} \frac{\partial \uu{v}}{\partial t} + \left( \begin{array}{ccc}
                         \tilde{u} - \tilde{c}\cos\theta & 0 & 0 \\
                          0 & \tilde{u} & 0 \\
                          0 &  0 & \tilde{u} + \tilde{c}\cos\theta \\
                                \end{array} \right )
\frac{\partial \uu{v}}{\partial x}
          + \left( \begin{array}{ccc}
                         \tilde{v} - \tilde{c}\sin\theta &  0 & 0 \\
                          0 & \tilde{v} & 0 \\
                          0 &  0 & \tilde{v} + \tilde{c}\sin\theta \\
                                \end{array} \right )
\frac{\partial \uu{v}}{\partial y}  = \uu{s} + \uu{r}
\end{equation}
with
$$
\uu{s} =  \trivek{s_1}{%
s_2}{%
s_3} = \trivek{\frac{1}{2}\tilde{h}( \sin\theta \frac{\partial v_2}{\partial x}
 + \cos\theta \frac{\partial v_2}{\partial y})}{%
 g \sin\theta(\frac{\partial v_1}{\partial x} - \frac{\partial v_3}{\partial x}) -
 g\cos\theta(\frac{\partial v_1}{\partial y} - \frac{\partial v_3}{\partial y}) }{%
 \frac{1}{2}\tilde{h}(\cos\theta \frac{\partial v_2}{\partial x}
 - \sin\theta \frac{\partial v_2}{\partial y})}.
$$
Let us denote by $\uu{x}_{\ell}$ the $\ell$-th bicharacteristic corresponding to the
$\ell$-th equation of system (\ref{so_3}). The bicharacteristic $\uu{x}_{\ell}$ is
defined in the following way $$ \frac{d\uu{x}_{\ell}(s)}{ds} =
\left(\begin{array}{rr}
                                 b^1_{\ell\ell} \\
                                 b^2_{\ell \ell}
                                \end{array} \right ),
$$
where $b^1_{\ell \ell}, b^2_{\ell \ell}$ are the diagonal entries of the matrices
$\umat{B}_1, \umat{B}_2$, respectively. The bicharacteristics $\uu{x}_{\ell}$ create
the surface of the so-called bicharacteristic cone, see Fig.~1, with the apex $P =
(x, y, t_n + \tau)$ and the footpoints
\begin{eqnarray}
 Q_1(\theta)  & = & (x - (\tilde{u} - \tilde{c} \cos\theta) \tau,
                     y - (\tilde{v} - \tilde{c} \sin\theta) \tau,t_n),
\nonumber\\
 Q_2 & \equiv &  Q_0 =  (x - \tilde{u}\tau, y - \tilde{v}\tau, t_n),
\nonumber\\
 Q_3(\theta)  & = & (x - (\tilde{u} + \tilde{c} \cos\theta) \tau,
                     y - (\tilde{v} + \tilde{c} \sin\theta) \tau,t_n).
\nonumber
\end{eqnarray}
Remember that $\tau = \Delta t/2$ in our case. Integrating each equation of
(\ref{so_3}) along the corresponding bicharacteristic from the apex $P$ down to the
footpoints $Q_{\ell}$ we get
\begin{eqnarray}
\label{w}
 v_{\ell}(P)   = v_{\ell}(Q_{\ell}) + \int_{t_n}^{t_n + \tau}
  s_{\ell}(Q_{\ell}(\tilde{t})) + r_{\ell} (Q_{\ell}(\tilde{t})) \dd\tilde{t}, \quad \ell = 1,2,3.
\end{eqnarray}
\noindent Now multiplying (\ref{w}) with $\umat{R}$ from the left and averaging over all
directions we go back to the original variables $\uu{w}$
\begin{eqnarray}\label{U(P)}
& &\uu{w}(P) = \nonumber \\ \nonumber
& &\frac{1}{2\pi} \int_0^{2\pi}  \trivek{-v_1(Q_1(\theta),\theta) + v_3(Q_3(\theta),\theta)}{%
 \frac{g}{\tilde{c}}\cos\theta v_1(Q_1(\theta),\theta) + \sin\theta v_2(Q_2(\theta),\theta) + \frac{g}{\tilde{c}}\cos\theta v_3(Q_3(\theta),\theta))}{%
 \frac{g}{\tilde{c}}\sin\theta v_1(Q_1(\theta),\theta) - \cos\theta v_2(Q_2(\theta),\theta) + \frac{g}{\tilde{c}}\sin\theta v_3(Q_3(\theta),\theta))}
\dd \theta
\nonumber\\[0.7cm]
& & + \frac{1}{2\pi} \int_0^{2\pi} \trivek{- s_1^{'}(\theta) - r_1^{'}(\theta) + s_3^{'}(\theta) + r_3^{'}(\theta)}{%
 \frac{g}{\tilde{c}}\cos\theta (s_1^{'}(\theta)+ r_1^{'}(\theta)) + \sin\theta (s_2^{'}(\theta) + r_2^{'}(\theta))
 + \frac{g}{\tilde{c}}\cos\theta (s_3^{'}(\theta)+r_3^{'}(\theta))}{%
 \frac{g}{\tilde{c}}\sin\theta ( s_1^{'}(\theta)+ r_1^{'}(\theta)) - \cos\theta (s_2^{'}(\theta) + r_2^{'}(\theta)
 + \frac{g}{\tilde{c}}\sin\theta ( s_3^{'}(\theta) + r_3^{'}(\theta))}
 \dd \theta,
\nonumber\\
\end{eqnarray}
where $ s_{\ell}^{'}(\theta) = \int_{t_n}^{t_n + \tau} {s_{\ell}
(\uu{x}_{\ell}(\tilde{t}, \theta), \tilde{t}, \theta) \dd \tilde{t}} $ is an
integral along the $ \ell $-th bicharacteric and the analogous notation holds for
source terms $r_{\ell}$. It should be noted that the source term $s_\ell$ in
(\ref{w}) arrives from to the multidimensionality of the system, whereas the source
term $r_\ell$ is a physical source term.

Now, we have $\lambda_1 = - \lambda_3$, $Q_1(\theta + \pi) = Q_3(\theta)$ and  the
characteristic variables $v_{\ell}$ are $2 \pi$-periodic. Applying the Gauss
integration, cf.\ (\ref{eq:4}), in order to avoid the derivatives of dependent
variables appearing in $\uu{s}$ we can, after analogous computations as in
\cite{mathcom,sisc}, reformulate the exact integral equations
(\ref{U(P)}) in the following way
\begin{eqnarray}
\label{eq3-ap} h\left( P\right) &=& \frac{1}{2\pi}\int_{0}^{2\pi}
   h\left(Q\right) - \frac{\tilde{c}}{g}u\left(Q\right)\cos\theta-
    \frac{\tilde{c}}{g}v\left(Q\right)\sin\theta \dd\theta
    \phantom{mmmmmmmmmmmmmmmm}
\nonumber \\
&& - \frac{1}{2\pi} \int_{t_n}^{t_n + \tau} \frac{1}{t_n + \tau - \tilde{t}}
     \int_{0}^{2 \pi}\frac{\tilde c}{g}\left( u(\tilde Q) \cos \theta + v(\tilde Q) \sin \theta\right)
     \dd \theta \dd \tilde{t}
\\ \nonumber
&& + \frac{1}{2\pi} \tilde{c} \int_{t_n}^{t_n + \tau}  \int_{0}^{2 \pi} \left(
b_x(\tilde Q) \cos \theta + b_y (\tilde Q ) \sin \theta \right)     \dd \theta \dd
\tilde{t} \nonumber
\\ \nonumber
&& - \frac{1}{2\pi} \frac{\tilde{c} f}{g} \int_{t_n}^{t_n + \tau}  \int_{0}^{2 \pi}
\left( v(\tilde Q) \cos \theta - u (\tilde Q ) \sin \theta \right)     \dd \theta
\dd \tilde{t},
\end{eqnarray}
\begin{eqnarray}
\nonumber \label{eq4-ap}
 u\left(P \right) &=&\frac{1}{2}u\left( Q_0\right) + \frac{1}{2 \pi}\int_{0}^{2\pi}
   - \frac{g}{\tilde{c}}  h\left(Q\right) \cos\theta
   + u \left(Q \right)\cos^2\theta  + v\left(Q \right)\sin\theta\cos\theta \, \dd\theta   \nonumber \\
&& - \frac{g}{2}  \int_{t_n}^{t_n + \tau} \left( h_x(\tilde Q_0) + b_x (\tilde
Q_0)\right)  \dd \tilde{t}
\\ \nonumber
&& - \frac{g}{2\pi}  \int_{t_n}^{t_n + \tau} \int_0^{2 \pi} \left(b_x(\tilde Q)
\cos^2 \theta + b_y(\tilde Q) \sin \theta \cos \theta \right) \dd \theta \dd
\tilde{t}
\\  \nonumber
&&+\frac{1}{2\pi} \int_{t_n}^{t_n + \tau} \frac{1}{t_n + \tau - \tilde{t}}
 \int_{0}^{2 \pi} \left( u(\tilde Q) \cos 2 \theta + v(\tilde Q) \sin  2\theta \right)
 \dd \theta \dd \tilde{t}
\\ \nonumber
&& + \frac{f}{2}  \int_{t_n}^{t_n + \tau}  v(\tilde Q_0) \dd \tilde{t} +
\frac{f}{2\pi} \int_{t_n}^{t_n + \tau} \int_0^{2 \pi} \left(v(\tilde Q) \cos^2
\theta - u(\tilde Q) \sin \theta \cos \theta \right) \dd \theta \dd \tilde{t},
\end{eqnarray}
\begin{eqnarray}
\label{eq5-ap}
\nonumber  v\left(P \right) &=&\frac{1}{2}v\left( Q_0\right) + \frac{1}{2
\pi}\int_{0}^{2\pi}
   - \frac{g}{\tilde{c}} h\left(Q\right) \sin\theta
   + u\left(Q \right)\sin\theta\cos\theta + v \left(Q \right)\sin^2\theta
   \, \dd\theta   \nonumber \\
&& - \frac{g}{2}  \int_{t_n}^{t_n + \tau} \left( h_y(\tilde Q_0) + b_y (\tilde
Q_0)\right)  \dd \tilde{t}
\\ \nonumber
&& - \frac{g}{2\pi}  \int_{t_n}^{t_n + \tau} \int_0^{2 \pi} \left(b_x(\tilde Q) \sin
\theta \cos \theta + b_y(\tilde Q) \sin^2 \theta \right) \dd \theta \dd \tilde{t}
\\  \nonumber
&&+\frac{1}{2\pi} \int_{t_n}^{t_n + \tau} \frac{1}{t_n + \tau - \tilde{t}}
 \int_{0}^{2 \pi} \left( u(\tilde Q) \sin 2 \theta - v(\tilde Q) \cos  2\theta \right)
 \dd \theta \dd \tilde{t}
\\ \nonumber
&& - \frac{f}{2}  \int_{t_n}^{t_n + \tau}  u(\tilde Q_0) \dd \tilde{t} +
\frac{f}{2\pi} \int_{t_n}^{t_n + \tau} \int_0^{2 \pi} \left( v(\tilde Q) \sin \theta
\cos\theta - u(\tilde Q) \sin^2 \theta  \right) \dd \theta \dd \tilde{t}.
\end{eqnarray}
Recall that $\tilde Q_0 = (x- \tilde{u} (t_n + \tau - \tilde t), y- \tilde{v} (t_n +
\tau - \tilde t), \tilde t)$, $\tilde Q  = (x - \tilde{u} (t_n + \tau- \tilde t) + c
(t_n + \tau - \tilde t) \cos \theta, y- \tilde{v}(t_n + \tau - \tilde t) + c (t_n +
\tau - \tilde t )\sin \theta, \tilde t)$ stays for an arbitrary point on the mantle and
$Q = Q(\tilde t)\Big |_{\tilde t = t_n} $ denotes a point at the perimeter of the
sonic circle at time $t_n.$


\section{Proof of Lemma~\ref{lema32}}

\begin{proof}
We show here that the approximate integral equations (\ref{eq6})-(\ref{eq8}) are
consistent with the exact integral equations (\ref{eq3})-(\ref{eq5}),
i.e.~(\ref{eq3-ap})-(\ref{eq5-ap}). In (\ref{eq3}) the integral with bottom
topography terms can be rewritten using the polar-type transformation along the
mantle of the bicharacteristic cone, i.e. $x_{\tilde Q} = x + r (\cos \theta -
\frac{\tilde u}{\tilde c}), \, y_{\tilde Q} = y + r (\sin \theta - \frac{\tilde
v}{\tilde c})$, where $r = \tilde c (t_{n} + \tau - \tilde t)$ is the circle radius
at the time level $\tilde t \in [t_n, t_n + \tau].$ Thus, we have
$$
\frac{d b}{d r} (r, \theta) = b_x (x_{\tilde Q},y_{\tilde Q}) \cos \theta + b_y
(x_{\tilde Q},y_{\tilde Q}) \sin \theta - \frac{1}{\tilde c} \left( \tilde u b_x
(x_{\tilde Q},y_{\tilde Q}) + \tilde v b_y (x_{\tilde Q},y_{\tilde Q})  \right).
$$
Therefore,
\begin{eqnarray}
&&\frac{\tilde c}{2\pi} \int_{t_n}^{t_n + \tau}  \int_{0}^{2 \pi} \left( b_x(\tilde
Q) \cos \theta + b_y (\tilde Q ) \sin \theta \right) \dd \theta \dd \tilde{t}
\label{bb} \\ \nonumber &&=\frac{\tilde c}{2\pi} \int_{\tilde c \tau}^{0} \int_0^{2
\pi} \frac{\dd b(r, \theta)}{\dd r} \dd \theta (- \frac{\dd r}{\tilde c}) +
\frac{1}{2 \pi} \int_{t_n}^{t_n + \tau} \int_0^{2 \pi}  \tilde u b_x (\tilde Q)  +
\tilde v b_y (\tilde Q)
 \dd \theta \dd \tilde t
\\ \nonumber
&&= \int_{0}^{\tilde c \tau} \frac{\dd }{\dd r} \left(\frac{1}{2 \pi} \int_0^{2 \pi}
b \, \dd \theta \right) \dd r  + \frac{1}{2 \pi} \int_{t_n}^{t_n + \tau} \int_0^{2
\pi}  \tilde u b_x(\tilde Q)  +  \tilde v b_y(\tilde Q)   \dd \theta \dd \tilde t
\\ \nonumber
&&= \frac{1}{2 \pi} \int_0^{2 \pi} b(Q) \dd \theta - b(P) +  \frac{1}{2 \pi}
\int_{t_n}^{t_n + \tau} \int_0^{2 \pi}  \tilde u b_x(\tilde Q)   + \tilde v
b_y(\tilde Q)   \dd \theta \dd \tilde t,
\end{eqnarray}
which yields the corresponding terms in (\ref{eq6}).

Further, we show that the integrals in (\ref{eq3}) containing the Coriolis forces are
of order ${\mathcal O}(\Delta t^2)$; note that $\tau = \Delta t/2$. Applying the
rectangle rule in time and the Taylor expansion from Lemma~\ref{aux_lemma:taylor_expansions}
in the center of the sonic circle $Q_0$ yields
\begin{eqnarray}
&&\int_{t_n}^{t_n + \tau} \int_{0}^{2 \pi}  v(\tilde Q) \cos \theta \dd \theta \dd
\tilde t = \tau \int_{0}^{2 \pi} v(Q) \cos \theta \dd \theta \nonumber
\\ \nonumber
&=& \tau \int_{0}^{2 \pi} (v(Q_0) \cos \theta + c \tau \, v_x (Q_0) \cos^2 \theta + c
\tau \, v_y(Q_0) \cos \theta \sin \theta  + O(\Delta t^2)) \dd \theta
\\ \label{CC}
&=& {\mathcal O}(\Delta t^2)
\end{eqnarray}
with an analogous approximation for the Coriolis forces in $y$-direction. Together
with (\ref{bb}) and (\ref{CC}) this yields the approximate integral equation (\ref{eq6}).

In the equation (\ref{eq4}) for velocity $u$  we apply for the mantle integrals
containing the bottom elevation terms the rectangle rule in time and the Taylor
expansion over the center $Q_0$ of the sonic circle $S_0$ at time $t_n$, which lead to
\begin{equation}
\label{TT}
 \frac{1}{2\pi}g \int\limits_{t_n}^{t_n + \tau}\int\limits_{0}^{2\pi}\left(
         b_x(\tilde Q )\cos\theta +   b_y(\tilde Q ) \sin\theta
   \right)\cos\theta\  \dd\theta \,\dd\tilde{t} =
\frac{g \tau}{2}  b_x(Q_0) + {\mathcal O}(\Delta t^2).
\end{equation}
To complete we eliminate the derivative by replacing the term $b_x( Q_0)$ by its
average over the sonic circle $S_0$ and applying the Gauss theorem
\begin{equation}
\label{gauss}
 b_x(Q_0) = \frac{1}{\pi \tilde c^2 \Delta {t}^2} \int_{S_0} b_x (Q) \, \dd x \dd y
+ O(\Delta t^2) = \frac{1}{\pi \tilde c \tau} \int_0^{2 \pi} b (Q) \cos \theta \,
\dd \theta + {\mathcal O}(\Delta t^2),
\end{equation}
which after substitution into (\ref{TT}) yields
\begin{equation}
\frac{1}{2\pi}g\int\limits_{t_n}^{t_n + \tau} \int\limits_{0}^{2\pi}\left(
         b_x(\tilde Q )\cos\theta +   b_y(\tilde Q ) \sin\theta
   \right)\cos\theta\  \dd\theta \, \dd\tilde{t} =
\frac{g}{\tilde c} \frac{1}{2 \pi} \int_0^{2 \pi} b (Q) \cos \theta \, \dd \theta +
{\mathcal O}(\Delta t^2). \label{eqbx}
\end{equation}
Rewriting the Coriolis forces terms using their primitives we obtain analogously to
(\ref{TT}) and (\ref{gauss})
\begin{eqnarray}
&&\frac{f}{2\pi} \int_{t_n}^{t_n + \tau} \int_0^{2 \pi} \left(v(\tilde Q) \cos^2
\theta - u(\tilde Q) \sin \theta \cos \theta \right) \dd \theta \dd \tilde{t}
\\ \nonumber \\
&& = \frac{g}{2\pi}\int_{t_n}^{t_n + \tau} \int_{0}^{2\pi}\left(
         V_x(\tilde Q )\cos\theta  - U_y(\tilde Q ) \sin\theta
   \right)\cos\theta\  \dd\theta \, \dd\tilde{t}
   \nonumber \\ \nonumber
&& =  \frac{g}{\tilde c} \frac{1}{2 \pi} \int_0^{2 \pi} V (Q) \cos \theta \, \dd
\theta + {\mathcal O}(\Delta t^2).
\end{eqnarray}
This balances together with (\ref{eqbx}) the analogous term with $h(Q) \cos \theta$
in (\ref{eq4}). The integral along the middle bicharacteristic
$$
\frac{g}{2}  \int_{t_n}^{t_n + \tau} \left( h_x(\tilde Q_0) + b_x (\tilde Q_0) -
\frac{f}{g} v(\tilde Q_0) \right) \, \dd \tilde{t} = \frac{g}{2} \int_{t_n}^{t_n +
\tau} \left( h_x(\tilde Q_0) + b_x (\tilde Q_0) - V_x (\tilde Q_0) \right) \, \dd
\tilde{t}
$$
can be approximated in a similar way as (\ref{gauss}) applying the Gauss theorem at
each intermediate circular section at $\tilde t$ along the mantle of the
bicharacteristic cone. Substituting into (\ref{eq4})  gives (\ref{eq7}).
Approximation (\ref{eq8}) for the velocity $v$ is obtained in an analogous way as
(\ref{eq7}).
\phantom{mm}\hfill \end{proof}


\end{document}